%% file: main.tex
\documentclass[preprint,12pt,authoryear]{elsarticle}

\usepackage{amssymb}
\usepackage{amsmath}
\allowdisplaybreaks

\usepackage{pdflscape}

\usepackage{tikz}
\usetikzlibrary{arrows.meta}
\usepackage{pgfplots}
\pgfplotsset{compat=1.9}
\usetikzlibrary{patterns}
\usepackage{pgf-pie}

\usepackage{graphicx}
\usepackage{tabularx}
\usepackage{multirow}
\usepackage{subcaption}
\usepackage{hyperref}

\journal{European Journal of Operational Research}

\begin{document}

\begin{frontmatter}



\title{Literature Survey on the Container Stowage\\ Planning Problem}


\author[inst1]{Jaike van Twiller}

\affiliation[inst1]{organization={IT University of Copenhagen},
            addressline={Rued Langgaards Vej 7}, 
            city={Copenhagen},
            postcode={2300}, 
            country={Denmark}}

\author[inst2]{Agnieszka Sivertsen}
\author[inst3]{Dario Pacino\corref{cor1}}
\ead{darpa@dtu.dk}
\cortext[cor1]{}
\author[inst1]{Rune Møller Jensen}

\affiliation[inst2]{organization={Roskilde University},
            addressline={Universitetsvej 1}, 
            city={Roskilde},
            postcode={4000}, 
            country={Denmark}}

\affiliation[inst3]{organization={Technical University of Denmark},
            addressline={Bygningstorvet 116B}, 
            city={Kgs. Lyngby},
            postcode={2800}, 
            country={Denmark}}

\begin{abstract}
Container shipping drives the global economy and is an eco-friendly mode of transportation. A key objective is to maximize the utilization of vessels, which is challenging due to the NP-hardness of stowage planning. This article surveys the literature on the Container Stowage Planning Problem (CSPP). We introduce a classification scheme to analyze single-port and multi-port CSPPs, as well as the hierarchical decomposition of CSPPs into the master and slot planning problem. Our survey shows that the area has a relatively small number of publications and that it is hard to evaluate the industrial applicability of many of the proposed solution methods due to the oversimplification of problem formulations. To address this issue, we propose a research agenda with directions for future work, including establishing a representative problem definition and providing new benchmark instances where needed. 

\end{abstract}


\begin{highlights}

\item Literature review and classification scheme for the Container Stowage Planning Problem
\item Identification of four significant groups of research: single-port, multi-port, master planning, and slot planning
\item Comparison of problem formulations and solution approaches
\item Description of minimal representative problem definition
\item Publication of benchmarks instances
\item Research agenda and discussion on the state of the art

\end{highlights}


\end{frontmatter}


\section{Introduction}
\label{sec:introduction}
\input{chapters/01_introduction.tex}

\section{The container stowage planning problem}
\label{sec:stowagePlanning}
\input{chapters/02_complexity_problem.tex}

\section{Classification scheme}
\label{sec:literatureSurvey}
\input{chapters/03_survey_results.tex}


\section{Research Agenda}
\label{sec:mathModels}
\input{chapters/05_research_agenda.tex}
\section{Conclusion}
\label{sec:conclusion}
This paper provides a review of the literature that studies the Container Stowage Planning Problem. The studies are summarized according to a classification scheme that outlines the fundamental characteristics of the problem and the applied solution approaches. As there is a lack of a common understanding of the problem characteristics, this paper provided a description of a representative problem definition based on several years of academic and industrial collaborations. In light of this definition, a research agenda is proposed for each of the major branches of research in the Container Stowage Planning Problem (single-port planning, multi-port planning, master planning, and slot planning). Moreover, this paper identifies, and in one case, provides publicly available benchmark sets in the hope that future research will make use of them as a reference point and a way to compare results. Where possible, these benchmarks have been used to compare recent research results, and provide some computational comparison. It is our hope that this survey will help improve the field and acts as inspiration for future developments.

\section*{Acknowledgements}
This work is partially funded by the Danish Maritime Fund (grant nr. 2021-069) and the Innovation Fund Denmark (grant no. 1044-00145A).

\clearpage
\appendix
\section{Classification tables}
\label{app:classification_tables}
\input{appendices/02_full_tables.tex}
\clearpage
\section{Collected result tables}
\label{app:result_tables}
\input{appendices/01_appendixA.tex}
\clearpage
\section{Master planning formulations}
\input{chapters/04_pacino11model}

\def\urlprefix{}
\def\url#1{}

\bibliographystyle{formatting/elsarticle-harv}
\bibliography{references/stowage_papers}

\end{document}

%% file: chapters/01_introduction.tex
Container shipping is an underappreciated business. Most people know little about container vessels and the media coverage often focuses on negative aspects. The truth is that container shipping is the most environmentally friendly mode of transportation with the least CO2 emissions per metric ton of goods shipped per kilometer \citep{InternationalChamberofShipping2023EnvironmentalTransport}.
Economically, it runs the supply chains of the world, and is in fact believed to have been more important for globalization than freer markets \citep{TheEconomist2013FreeHero}.

From an operations research (OR) point of view, the overall objective of container shipping is to maximize the utilization of vessels while minimizing operational costs. Container vessels, however, are challenging. To stow containers on them according to this objective is a combinatorial optimization problem with an unusually wide range of complex constraints and objectives including seaworthiness requirements, stacking rules, crane utilization, and fuel consumption. For instance, minimizing the number of containers that block containers in lower stacks tiers is NP-hard \citep{Avriel2000ContainerGraphs}. 

Research on the container stowage planning problem (CSPP) is unfortunately scarce compared to other areas of OR. In this survey article, we were only able to find 54 key contributions in our literature search covering the 67 years that have passed since the first container vessel sailed in 1956. CSPP studies are challenged in several ways. First, container shipping has been deregulated about thirty years after the airline industry when the conference system was outlawed in Europe in 2008. For that reason, there has been less focus on advanced capacity management systems. Second, as mentioned above, the domain is unknown to most researchers and it is only recently that a comprehensive description \citep{Jensen2018ContainerPlanning} and benchmark suite (\cite{Larsen2021AProblem}) was published. Finally, the problem is highly complex and important aspects are subtle and hard to model. It is not clear how to study it in a reduced representative form suitable for scientific research.

In this article, we survey the CSPP literature and propose a research agenda with directions for future work. Our main conclusion is that while several algorithmic frameworks and problem decompositions have been investigated, it is still a maturing research area with relatively few publications and a lack of benchmark suites and problem definition consensus. As a large and sustainable mode of transportation, this state of affairs is important to change.   

The remainder of the article is organized as follows. Section~\ref{sec:stowagePlanning} gives an overview of the CSPP focusing on its key combinatorial components. Section~\ref{sec:literatureSurvey} defines our literature classification scheme. Section~\ref{sec:contributions} reviews the literature and divides the contributions into single-port planning, multi-port planning, master planning, slot planning, computational complexity, and other relevant work. Sections~\ref{sec:mathModels} presents a comprehensive research agenda with suggested directions for future work in each area. Finally, Section~\ref{sec:conclusion} draws conclusions from the survey.

%% file: chapters/02_complexity_problem.tex
The CSPP is extensive and it is beyond the scope of this article to describe it in detail. For a full introduction, we refer the reader to \cite{Jensen2018ContainerPlanning}. Instead, this section focuses on the combinatorial structure of the problem and identifies key combinatorial aspects that are important to model in academic studies in order to ensure that they are representative of the real-world problem.

To this end, we first briefly describe container shipping. Shipping lines are similar to bus lines but on water. Their fleet of vessels is assigned to closed-loop services with fixed schedules. 
The CSPP is to decide where on the vessel the booked cargo to load is stowed. It is an operational problem that is solved by a stowage team. Even though stowage plans are made one port at a time, the CSPP is multi-port in nature as the cargo placed in the current port affects the vessel condition and free capacity in future ports. The input to the CSPP in this multi-port form is the arrival condition of the vessel in the first part (the so-called {\em remain-on-board} (ROB) condition), the {\em load-list} of cargo to load for each port call, and the vessel and terminal data. The result of the CSPP is a stowage plan for the first or all port calls. 
The primary objective of the CSPP is to load all booked cargo by maximizing the available capacity of the vessel. The secondary objective is to minimize terminal fees and the operational costs of the vessel. Finally, since cargo bookings are uncertain, stowage plans must be robust and allow many different cargo compositions in future ports.  

The complexity of stowage planning is due to the large size of container vessels that today can be more than 24,000 {\em twenty-foot equivalent units} (TEU) and a myriad of interacting seaworthiness requirements. To understand the essence of these combinatorial aspects, we need some physical insight into the problem. 

\begin{figure}[h!]
	\centering
		\includegraphics[scale=0.4]{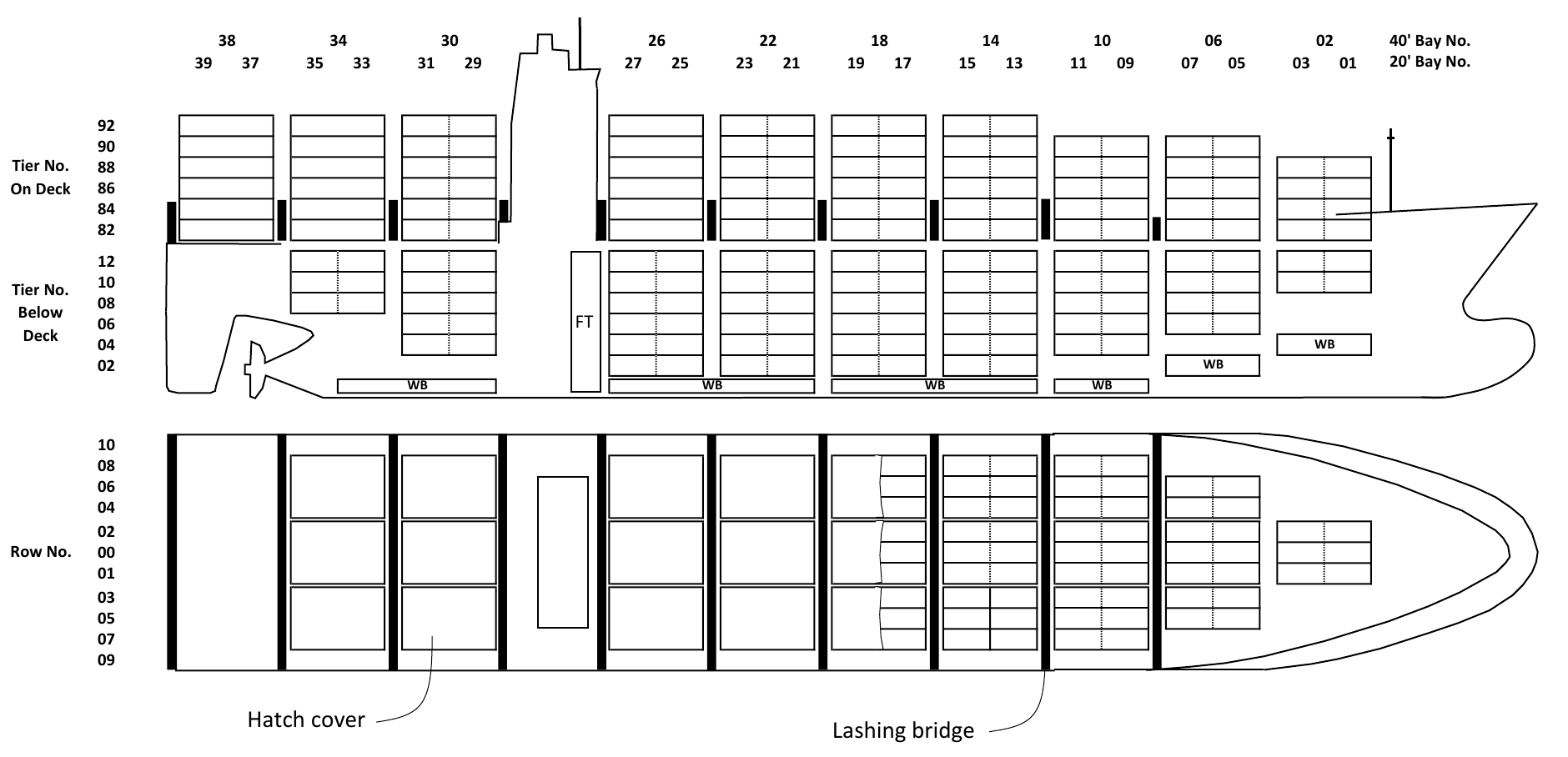} 
	\caption{Vessel side and top view.}
	\label{fig:vesselside}
\end{figure}

\begin{figure}[h!]
	\centering
		\includegraphics[scale=0.6]{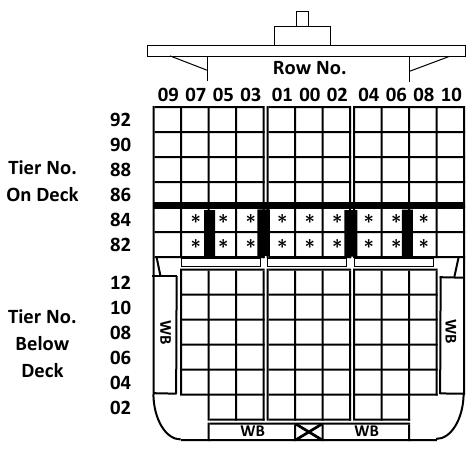} 
	\caption{Vessel bay front view.}
	\label{fig:vesselfront}
\end{figure}

Figure~\ref{fig:vesselside} and \ref{fig:vesselfront} show the cellular design of container vessels. The storage area consists of {\em bays} (02 - 38) with {\em stacks} or {\em rows} of {\em cells} that normally either can hold one 40' container or two 20' containers. The securing system {\em below deck} consists of {\em cell guides} that hold the stacks in place. Each hold is sealed with {\em hatch covers}. The stacks {\em on deck} rest on the hatch covers or the deck of the ship. The stacks are kept in place by {\em twist locks} that bind containers together and {\em lashing rods} that tie container corners to the deck or lashing bridges that are raised to increase stability. All stacks have weight limits. Below deck and fore on deck, they also have height limits. Both are essential to the model to get the volume and weight capacity of the vessel right. Some cells have power plugs for refrigerated containers ({\em reefers}). They are indicated with stars in Figure~\ref{fig:vesselside}. 

Most containers are 20', 40' and 45' long, 8' wide, and 8'6" high. Some of these (mostly 40' and 45') are 9'6" high and are called {\em highcubes}. The weight ranges from about four tons for empty containers to about 30 tons for heavy containers. 
The port, where a container is loaded, is referred to as the  {\em port of load} (POL), while the port where it is unloaded is called the {\em port of discharge} (POD). {\em Specials} include reefers, IMDG, and OOG containers. IMDGs carry dangerous cargo and may be required to be segregated. OOGs (Out-of-Gauges) are over-dimensioned, and reefers are refrigerated containers and require electric power. 

An important combinatorial aspect of stowage planning is {\em restow} minimization. A restow happens when a container is unloaded and loaded again before its POD. Restows can be {\em voluntary} or {\em mandatory}. Mandatory restows are well-studied and are caused by {\em overstowage}. {\em Stack overstowage} happens when a container in a stack is stowed on top of a container to an earlier port. The top container must be restowed in order for the crane to reach the container below it.
Restow fees are high relative to the profit margin of each container. For that reason, the number of restows must be minimized. Unfortunately, minimization of mandatory restows has been shown to be an NP-hard problem \citep{Avriel2000ContainerGraphs}. Even when overstowage is simplified to {\em hatch overstowage}, which is overstowage between containers stowed above and below a hatch cover \citep{Tierney2014OnProblems}. Voluntary restows are also subject to restow fees, but these are made on purpose to increase the capacity of the vessel. For instance, a non-reefer container in a reefer slot can be restowed to fit an extra reefer. Since restow minimization is NP-hard and plays a central role in stowage planning, a proper optimization model should include the POL and POD of containers and represent voluntary and mandatory restows.


Another combinatorial element is stowage rules. 20' containers can not be stowed on top of 40' containers and IMDG containers must be segregated depending on their content. Since reefer containers are spark generators, they must be placed away from most IMDGs. Moreover, since reefer plugs typically are at bottom of stacks, a 40' reefer kills 20' capacity. Wrt. length, the interactions above can be modeled by just 20' and 40' containers. 45' containers must be placed above lashing bridges or in 45' bays, but this could be a minor combinatorial issue.

It is an open question whether IMDG segregation is NP-hard. In practice, it is challenging to deal with and should for that reason be modeled on IMDG heavy services.~\footnote{Carriers also can have their own rules, so-called {\em handling instructions}. They are easier to deal with than IMDGs and can be ignored} Reefer containers should always be modeled as they have limited positions available and affect other containers. Highcubes are frequent and often take up more than half of the volume. Due to height limitations below deck and in fore bays on deck, it is important to mix highcubes and normal containers right to utilize all the stack volume. It is unclear to us, though, how critical this combinatorial aspect is. 

Container weights have a high impact on stowage conditions. They can vary considerably and affect many seaworthiness requirements that are associated with the weight distribution of the vessel. The vertical centre of gravity (VCG) must be sufficiently low to avoid that the vessel capsizes. This transversal stability is measured by the {\em metacentric height} (GM) that must be above a certain minimum. A high GM, on the other hand, makes the vessel stiff and may cause the {\em lashing forces} in lashing rods to be exceeded in particular if heavy containers are stowed high in stacks on deck. The longitudinal centre of gravity (LCG) determines the difference between the fore and aft draft of the vessel called the {\em trim}. A vessel has a maximum allowed draft (maximum {\em loadline}) which may be further reduced due to port draft restrictions. The trim affects the {\em line of sight} (LOS) from the bridge which must be sufficient. The trim also affects the energy efficiency of the engine. The longitudinal weight distribution further causes stress forces on the vessel that all must be within limits.

The question is to what level of granularity we need to model these weight-related constraints of stowage planning. With respect to container weight, homogeneous weight is out of the question as it is unrealistic and compromises a reasonable representation of most constraints above. However, even just three or four weight classes enable us to model the combinatorial interactions of these constraints. Wrt. the weight-related constraints, it is key to model GM as it both affects transversal stability and lashing forces. Lashing forces are computed using complex mechanical simulation that makes them hard to embed in optimization models. Nevertheless, they should not be ignored as they easily make one or more top tiers on deck impossible to use (i.e., more than 5\% of the total volume capacity). Trim and list requirements can be translated into total weight and box constraints on LCG and TCG. Stress forces also impact capacity. In particular, bending moments (BM) can reduce the weight capacity in the fore and aft of the vessel and for that reason should not be ignored. Shear forces (SF) can be high if mixing full and empty bays. Torsion moments (TM) are mainly a problem on large vessels with wide hulls. We assess that their combinatorial impact is secondary to BM in most cases. With respect to the vessel structure, the hull can be represented by a sequence of box-shaped sections such that hydrostatic equilibrium can be linearly defined. This representation is fairly accurate for draft, trim, and stress forces down to just six sections \cite{Jensen2018TheCapacity}. 
Ballast water is also important to model as it provides a flexible weight buffer that can be used to ease most weight-related constraints. 

Another combinatorial aspect of stowage planning is terminal constraints and objectives. Among the constraints, we have draft restrictions and crane work height. Normally, the impact of these restrictions is limited. The hard combinatorial aspects have to do with minimizing the port stay. This is important since a short port stay may allow the vessel to catch up on the schedule or save fuel by reducing the speed between ports. There mainly are two ways to minimize the port stay: 1) minimize the total number of quay crane moves, and 2) minimize the makespan of the cranes. With respect to  1), we need to minimize the number of restows which we already saw is NP-hard. With respect to 2), it is important to understand the relation between the terminal and the shipping line. Typically, the terminal guarantees a certain number of moves within a time window without specifying the number of cranes working on the vessel during this period. The shipping line can observe the average number of assigned cranes and distribute the moves along the vessel such that these cranes can work in parallel.  A complicating matter is that two quay cranes due to their width are unable to work simultaneously on two adjacent bays. The total work time is therefore given by the {\em long crane} which is the largest number of moves $m$ of any pair of adjacent bays. Let $M$ denote the total number of moves of the vessel. The {\em crane intensity} (CI) is then defined as $M / m$. It is an upper bound on the number of cranes that can work in parallel on the vessel without any of them being idle. The stowage plan should have a CI that is larger than the average number of cranes assigned to the vessel by the terminal.  

Finally, we turn to a rather hidden combinatorial aspect of stowage planning that has to do with the robustness of plans. To achieve high flexibility for the kind of containers that can be loaded in future ports, stowage planners use certain stowage patterns. To define these patterns, let a {\em block} denote a storage space either above or below a hatch cover. Hence, if a bay has three hatch covers, it has six blocks in total: one center and two wing blocks on and below deck. A basic stowage pattern is to avoid mixing PODs in a block. This ensures by design that the block has no stack restows and empties the whole block in the port of discharge such that the POD choice of the containers to load back into the block is as free as possible. To avoid hatch restows by design as well, we need to require that the two blocks above and below a hatch cover hold the same POD. We call this pattern {\em block stowage}. To avoid only changing weight on one side of the vessel and risk excessive torsion moments, the block stowage pattern is often extended to {\em paired block stowage}, where blocks in the wing hold the same POD, while the center may hold another POD. 

Despite the fact that stowage patterns reduce the space of possible plans, they seem to increase the problem complexity substantially \cite{Christensen2017AStowage}. Due to this impact and the fact that the patterns are industry standard, they should be included to some degree in a representative model of the problem.

%% file: chapters/03_survey_results.tex
Before diving into the survey and the proposed classification scheme, let us clarify the search strategy used to collect and select the relevant literature. The articles were found by using Google Scholar using \textit{Stowage Planning} as the keyword. The title of the resulting articles has then been evaluated and publications that were clearly not relevant have been removed. The remaining publications were further filtered by reading the abstracts. It is here that publications focusing on, e.g., the packing of cargo into a single container were removed. All the references of the remaining publications have been analyzed and missing contributions have been added to the list. A number of publications presenting minor incremental work have been removed from the analysis. It is a key finding that few if any studies fulfill even the minimal representation requirements discussed in  Section~\ref{sec:stowagePlanning}. Hence, the 54 publications selected for this survey either treat most of them or present an original solution approach.

Articles that do not focus directly on solving the CSPP (or one of its sub-problems) were also not considered in the survey, e.g., various descriptions of visualization tools \citep{Aye2010VisualizationSystem, Song2010ResearchLogistics}, loading computers \citep{Nugroho2021RegulatoryConsiderations, Wu2021ResearchShips}, container data sharing systems \citep{Conca2018AutomationPlanning}, crane scheduling tools \citep{Hsu2021JointApproach} or loading sequence planners \citep{Serban2017AStrategy}. Survey papers \citep{Zhang2008}, decision support tools \citep{Martin1988ComputerizedEvaluation, Saginaw1989DecisionPlanning} and instance generators \citep{Cruz-reyes2013ConstructivePlanning} are not discussed in this classification either.

\begin{figure}[!h]
	\centering
		\includegraphics[scale=0.46]{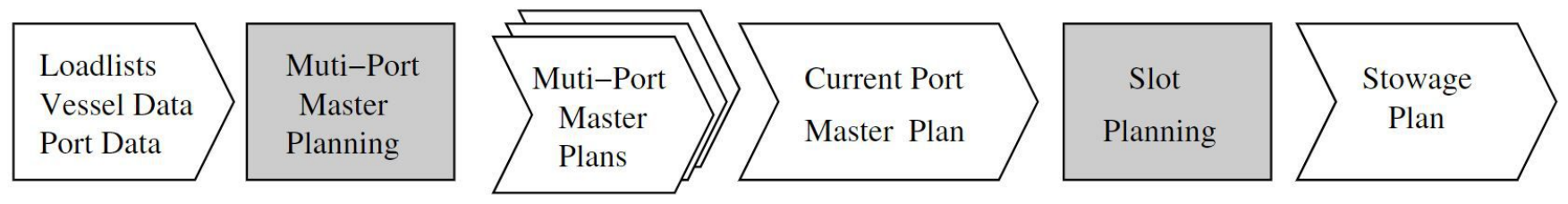} 
	\caption{Hierarchical decomposition of stowage planning into master and slot planning from \cite{Pacino2011FastVessels}.}
	\label{fig:decomposition}
\end{figure}

To reduce the complexity of the problem, the seminal work of \cite{Wilson2000ContainerSolutions} introduced a hierarchical decomposition to obtain a sequence of tractable interdependent sub-problems. It quickly became popular and widely used in several studies (e.g., \cite{Kaisar2006ATRANSPORTATION, Pacino2018CranePlanning}). For that reason, the authors propose a division of the problem into two sub-problems as shown in Figure \ref{fig:decomposition}. The first sub-problem is named the {\em Master Planning Problem} (MPP) and its solution (a {\em master plan}) is the assignment of groups of containers to storage areas of the vessel (blocks). The idea is to address high-level constraints and objectives, such as overall weight distribution, crane utilization, hatch cover moves, and cargo consolidation. Containers in the MPP are often grouped by their weights and types, and a master plan is created for all the ports in the voyage. The second sub-problem is called the {\em Slot Planning Problem} (SPP), which uses the MPP as input. Given the assignment of groups of containers to each block, the SPP assigns individual containers to slots in the block for every port separately. The assignment fulfills only low-level constraints, such as stacking rules, capacity constraints, and overstowage constraints. This results in a complete (multi-port) stowage plan.

Considering the above decomposition, we grouped the publications as follows: publications treating the full CSPP and publications focusing only on one of the sub-problems. Among the works solving the full CSPP, we distinguish between multi-port and single-port stowage planning, where the latter approach only takes a loadlist for a single port into account. For each category, we provide a qualitative comparison and a quantitative comparison when possible. Analogously, we aggregated publications for the MPP and the SPP and analyzed them independently.

\begin{table}
\centering
\caption{\label{tab:classification_Scheme}Classification scheme.}
\resizebox{\linewidth}{!}{%
\begin{tabular}{lll}
\textbf{Label}            & \textbf{Value}     & \textbf{Description}                                                                                                                                     \\ 
\hline
\textit{Cargo}            &                    & Cargo characteristics                                                                                                                                    \\ 
\hline
                          & \textit{Uni}       & Uniform weight containers                                                                                                                                \\ 
\cline{2-3}
                          & \textit{Class}     & Container grouped by weight classes                                                                                                                      \\ 
\cline{2-3}
                          & \textit{Mix}       & Mixed weight containers                                                                                                                                  \\ 
\hline
\textit{Hydro}            &                    & Hydrostatics                                                                                                                                             \\ 
\hline
                          & \textit{Rich}      & Stability and stress force constraints                                                                                                                   \\ 
\cline{2-3}
                          & \textit{Stab}      & Stability constraints only                                                                                                                               \\ 
\cline{2-3}
                          & \textit{Equi}      & Longitudinal, vertical and/or transversal equilibriums                                                                                                     \\ 
\cline{2-3}
                          & \textit{None}      & No hydrostatics constraints                                                                                                                              \\ 
\hline
\textit{CSPP aspects}     &                    & CSPP aspects present in the problem formulations                                                                                                         \\ 
\hline
                          & \textit{MinRe}     & Involuntary container restows are allowed as well as minimized                                                                                           \\ 
\cline{2-3}
                          & \textit{VolRe}     & Voluntary and involuntary container restows are allowed as well as
  minmized                                                                            \\ 
\cline{2-3}
                          & \textit{NARe}      & Container restows are not allowed                                                                                                                        \\ 
\cline{2-3}
                          & \textit{HR}        & Hatch restows created by hatch cover lifts                                                                                                               \\ 
\cline{2-3}
                          & \textit{RF}        & Refrigerated containers                                                                                                                                            \\ 
\cline{2-3}
                          & \textit{DG}        & Dangerous cargo                                                                                                                                           \\ 
\cline{2-3}
                          & \textit{BW}        & Ballast water                                                                                                                                            \\ 
\cline{2-3}
                          & \textit{La}        & Lashing forces                                                                                                                                           \\ 
\cline{2-3}
                          & \textit{CO}        & Crane operations                                                                                                                                             \\ 
\cline{2-3}
                          & \textit{BS}        & Block stowage                                                                                                                                            \\ 
\hline
\textit{Obj}              &                    & Elements of the objective function                                                                                                                       \\ 
\hline
                          & \textit{PS}        & Minimize port stay by minimizng overstowage and optimizing cranes work                                                                                   \\ 
\cline{2-3}
                          & \textit{VU}        & Maximize vessel utilization, consolidation                                                                                                               \\ 
\cline{2-3}
                          & \textit{H}         & Minimize fuel consulption, improve hydrostatics                                                                                                          \\ 
\hline
\textit{HD}               &                    & Hierarchical decomposition                                                                                                                               \\ 
\hline
\textit{Sc}               &                    & Size of computational study                                                                                                                              \\ 
\hline
                          & \textit{S}         & Small. Vessels with a capacity below 2,500 TEU, for slot planner blocks below 75 TEU                  \\ 
\cline{2-3}
                          & \textit{M}         & \begin{tabular}[c]{@{}l@{}}Medium. Vessels with a capacity between 2,500 and 15,000 TEU,\\~for slot planner blocks between 75 and 150 TEU\end{tabular}  \\ 
\cline{2-3}
                          & \textit{L}         & Large. Vessels with a capacity above 15,000 TEU, for slot planner blocks above 150 TEU                       \\ 
\hline
\textit{Solution methods} &                    & Applied optimization techniques                                                                                                                          \\ 
\hline
                          & \textit{Greedy}    & Greedy approach                                                                                                                                          \\ 
\cline{2-3}
                          & \textit{Exact}     & Exact algorithms                                                                                                                                         \\ 
\cline{2-3}
                          & \textit{Method 1/Method 2}    & Hybrid of two methods, e.g. exact and heuristic                                                                                                                               \\ 
\cline{2-3}
                          & \textit{MatHeu}    & Matheuristic                                                                                                                                             \\ 
\cline{2-3}
                          & \textit{TreeB}     & Tree-based approach                                                                                                                                      \\ 
\cline{2-3}
                          & \textit{NeighMeta} & Neighborhood metaheuristic                                                                                                                               \\ 
\cline{2-3}
                          & \textit{PopulMeta} & Population metaheuristic                                                                                                                                 \\ 
\cline{2-3}
                          & \textit{ML}        & Machine Learning                                                                                                                                         \\
\cline{2-3}
\end{tabular}
}
\end{table}

To have an overview of the classified publications and to be able to easily compare the problem formulations presented in them, we developed the classification scheme that is shown in Table \ref{tab:classification_Scheme}. The \textit{Cargo} attribute determines what types of containers the problem formulation includes. They are divided into three categories: there can be models where all containers have the same weight (\textit{Uni}) or varying weights (\textit{Mix}). Weight classes (\textit{Class}) can also be used, where each weight class corresponds to a weight range.

The attribute \textit{Hydro} specifies the level of hydrostatics constraints included in the problem formulation. The \textit{Rich} level indicates that both stability and stress force constraints are part of the formulation, whereas, \textit{Stab} means that only stability constraints like GM, trim and/or list are included. \textit{Equi} is related to longitudinal, vertical and/or transversal equilibrium. There is also a possibility that the problem formulation doesn't contain any hydrostatic considerations (\textit{None}).

We note that basic capacity limitations of a container vessel, e.g., constrained height and weight of stacks, are part of every problem formulation included in our classification.

The group called \textit{CSPP aspects} specifies which elements of the problem are included in the problem formulation. The first group is restow handling which can also vary in the problem formulations. A common approach is to minimize the number of restows in the stowage plan (\textit{MinRe}), whereas in some cases they are completely forbidden (\textit{NARe}). In some models, an attempt to create voluntary restows to be able to stow more containers is made while minimizing the number of involuntary restows (\textit{VolRe}). A few formulations include hatch restows (\textit{HR}) caused by hatch overstowage. Crane operations (\textit{CO}) indicate that crane work optimization (for more than one crane) is incorporated in the problem formulation. It can be considered either in the objective function or as a constraint. The inclusion of (\textit{BS}) indicates whether block stowage best practices are modeled. \textit{BW} indicates that is it allowed to use ballast water to fix hydrostatics and \textit{La} means that lashing forces are part of the problem formulation. Special containers might be included in the load list. It is indicated by \textit{RF} (refrigerated containers) and \textit{DG} (dangerous cargo).

The objective function (\textit{Obj}) might include different aspects of the CSPP. The focus might be on the port stay (\textit{PS}), vessel volume utilization (\textit{VU}), hydrostatics (\textit{H}) to for example minimize fuel consumption or a combination of these.

The problem might be divided into two or more sub-problems and solved using hierarchical decomposition (\textit{HD}). The attribute, \textit{SC}, indicates the scale of the computational study in terms of problem size and is grouped into three categories: small (\textit{S}), medium (\textit{M}) and large (\textit{L}). 

Notice that elements of the table might be written in parentheses. This indicates that this aspect is only partially incorporated into the problem formulation.

\section{Literature review} \label{sec:contributions}
In the following, we will compare some of the most influential or recent contributions in what we believe are the two main scientific areas of interest: single-port and multi-port stowage planning. The first two subsections describe publications in these areas. Several of the multi-port contributions solve the problem using hierarchical decomposition and include MPP and SPP models. The third and fourth subsections describe papers that exclusively address the MPP and SPP sub-problems. Finally, the last two subsections discuss studies of computational complexity and other relevant contributions.       

\subsection{Single-port container stowage planning}\label{sec:solMethods1}
The high impact that a stowage plan of an earlier port can have on later ports is what dictates the inclusion of cargo forecasts and, as a consequence, the solution of multi-port versions of the problem. Single-port versions of the problem are, however, still interesting as they can be seen as more lightweight operational plans or as sub-problems \citep{Delgado2012ABays}. 

Table \ref{tab:full_SPClassification} shows the classification of the single-port studies in this review. 
From the table, it is easy to see a general consensus that contributions must include the modeling of several container types and some aspects of vessel stability. A noticeable exception is the work of \cite{Delgado2012AStowage, Shen2017AProblem}, where a weight distribution of the cargo is assumed to be an input to the algorithms, and the work of \cite{Zhao2018ContainerSearch}, where no explanation is given for this omission.
In the studies of \cite{Sciomachen2003TheProblem,Ambrosino2004StowingProblem,Ambrosino2006AProblem,Sciomachen2007AProductivity,Ambrosino2010AnProblem, Cruz-Reyes2015,Li2020OptimizingShipping}, vessel stability is only considered as balanced weights on the four sections of the vessel (bow, stern, port, and starboard). Works including more accurate measures are more recent \citep{Cho1981DevelopmentPlanning,Hu2012CombinatorialTerminal,Zhu2020IntegerProblem,Larsen2021AProblem,ElYaagoubi2022Multi-objectiveSystem}. 

The vast majority of the literature focuses on the minimization of time at port, either in the form of the time spent moving containers or in the minimization of the number of restows. Only a subset of the studies, though, includes the modeling of restows due to hatch covers \citep{Cho1981DevelopmentPlanning,Delgado2012ABays,Zhu2020IntegerProblem,Larsen2021AProblem}, only \cite{Larsen2021AProblem} and \cite{Sciomachen2007AProductivity} includes workload distribution of the quay cranes, and a block stowage strategy is proposed in \cite{Larsen2021AProblem}.

All the mathematical models proposed for the single-port stowage planning problem use, or adapt, the formulation introduced by \cite{Ambrosino2004StowingProblem}. It is a  four-index formulation indicating whether a container $c$ is assigned to a slot in bay $b$, row/stack $r$ and tier $t$. Exceptions are the work of \cite{Larsen2021AProblem}, where symmetries in the container index are broken by the modeling of container classes, and  \cite{Zhu2020IntegerProblem}, where the decision variable is split in two. The first variable assigns containers to blocks and the second variable assigns containers to tiers within blocks, thereby effectively abstracting away the stack/row position. This formulation was able to solve problems up to 1000 TEUs, compared to the  198 TEUs of the original formulation \citep{Ambrosino2004StowingProblem}. The model by \cite{Larsen2021AProblem}, is deemed intractable even for medium-sized vessels (7300 TEUs; no data is available for smaller vessels). The formulation has also been extended to integrate other problems: the blocks relocation problem \citep{Li2020OptimizingShipping}, and barge assignment \citep{ElYaagoubi2022Multi-objectiveSystem}.

No efficient mathematical formulation or exact method has yet been found that can solve the single-port container stowage problem for real-size vessels, which explains the focus of the literature on heuristic approaches. Most solution methods rely on metaheuristics. Local search procedures that exchange containers with the aim to solve vessel stability are used by \cite{Cho1981DevelopmentPlanning,Ambrosino2006AProblem,Li2020OptimizingShipping}. Ant colony optimization is proposed by \cite{Ambrosino2010AnProblem}, two genetic algorithms are introduced in \cite{Hu2012CombinatorialTerminal, ElYaagoubi2022Multi-objectiveSystem}, and an adaptive large neighborhood search is presented in \cite{Larsen2021AProblem}. 

Construction heuristics are proposed by \cite{Sciomachen2003TheProblem} and \cite{Delgado2012AStowage}, where the latter is based on a cargo distribution obtained with a linear program. Other approaches include a heuristics branching procedure \citep{Sciomachen2007AProductivity}, a tree-search-based heuristics \citep{Zhao2018ContainerSearch}, and machine learning \citep{Shen2017AProblem}. 

Though several approaches have been proposed, the lack of a common benchmark and problem definition makes it hard to compare their performance. \cite{Larsen2021AProblem} has recently published a benchmark (available at \url{https://doi.org/10.11583/DTU.9916760}) in order to address this issue. Given that most approaches have been tested on rather small instances, further research on their performance on larger instances is valuable.

\begin{table}[h!]
\centering
\caption{\label{tab:full_SPClassification}Classification of single-port container stowage planning.}
\scalebox{.72}{
\scriptsize
\begin{tabular}{llllllll}
\textbf{Paper} & \textbf{Cargo} & \textbf{Hydro} & \textbf{CSPP aspects} & \textbf{Obj} & \textbf{Sc} & \textbf{HD} & \textbf{Solution methods}            \\ 
\hline
\cite{Cho1981DevelopmentPlanning}        & Mix            & Stab             & NARe, RF         & PS, VU       & S           &           & Exact/NeighMeta         \\ 
\hline
\cite{Sciomachen2003TheProblem} & Class          & Equi            & NARe, RF                 & PS           & S           &           & Exact/Greedy                 \\ 
\hline
\cite{Ambrosino2004StowingProblem}  & Class          & Equi            & NARe             & PS           & S           &           & Exact                  \\ 
\hline
\cite{Ambrosino2006AProblem}  & Class          & Equi            & NARe             & PS           & S           &           & Exact/NeighMeta         \\ 
\hline
\cite{Sciomachen2007AProductivity} & Class          & Equi            & NARe, CO             & PS           & S           &           & Exact/NeighMeta                  \\ 
\hline
\cite{Ambrosino2010AnProblem}  & Class          & Equi            & NARe                 & PS           & S           &           & NeighMeta  \\ 
\hline
\cite{Delgado2012AStowage}    & Mix            & None                  & NARe, RF                    & VU           & L           & \checkmark           & Exact/Greedy         \\ 
\hline
\cite{Hu2012CombinatorialTerminal}         & Class          & Stab             & MinRe             & PS, H        & S           &           & PopulMeta         \\ 
\hline
\cite{Cruz-Reyes2015} & Class          & Equi            & NARe                 & PS           & S           &            & Exact/Greedy         \\ 
\hline
\cite{Shen2017AProblem}       & Mix            & None                  & MinRe                 & PS           & S           &             & ML                     \\ 
\hline
\cite{Zhao2018ContainerSearch}       & Mix            & None                  & MinRe                 & PS           & S           &             & TreeB                  \\ 
\hline
\cite{Li2020OptimizingShipping}         & Mix            & Equi            & MinRe, NARe                 & PS           & S           &           & NeighMeta, Exact                 \\ 
\hline
\cite{Zhu2020IntegerProblem}        & Mix            & Rich                  & MinRe, HR, RF         & PS           & S           &             & Exact                  \\ 
\hline
\cite{Larsen2021AProblem}     & Class          & Rich                  & MinRe, HR, CO, BS, RF & PS, VU, H    & L           &             & NeighMeta                 \\ 
\hline
\cite{ElYaagoubi2022Multi-objectiveSystem} & Mix            & Stab             & MinRe, La             & PS, H        & S           &             & PopulMeta, Exact          \\
\hline
\end{tabular}
}
\end{table}

\subsection{Multi-port container stowage planning}\label{sec:solMethods2}
In addition to single-port plans, we recognize the importance of considering vessel conditions and cargo forecasts at future ports in multi-port planning \citep{Delgado2012ABays}. The goal is to find robust plans that maximize vessel utilization and minimize operational costs during the voyage.

Table \ref{tab:full_MPClassification} shows the classification of multi-port work, from which can be derived that little consensus exists on the modeling of cargo weights and stability. 
Most contributions combine stability with varying cargo types, e.g., varying cargo weights \citep{ Liu2011RandomizedPlans, Hu2017Multi-objectiveRoute, Parreno-Torres2021SolvingAlgorithm}, or special cargo as reefers and dangerous goods \citep{Hamedi2011CONTAINERSHIPOPERATIONS, Liu2011RandomizedPlans, Chang2022SolvingMode}. These hydrostatics can be extended with shear forces \citep{Shields1984, Pacino2011FastVessels} as well as bending moments \citep{Botter1992StowageSolution}.
The contributions by \cite{Kaisar2006ATRANSPORTATION, Azevedo2014SolvingMeta-heuristics, Azevedo2018SolvingAlgorithm, Li2018Multi-PortUncertainties} approximate stability as \cite{Ambrosino2004StowingProblem} by balancing cargo weight with respect to the vertical, longitudinal and transversal dimensions.
There is, however, also work that disregards hydrostatics in order to focus on other combinatorial aspects (e.g., \cite{Avriel1998StowageShifts, Wilson2000ContainerSolutions, Pacino2018CranePlanning}).

Most studies in Table \ref{tab:full_MPClassification} aim to minimize load and discharge moves or mandatory restows. There are, however, two exceptions, namely, voluntary restows (e.g., \cite{Avriel1998StowageShifts, Ding2015StowageShifts, Roberti2018APlans}), and not allowing restows  \citep{Li2018Multi-PortUncertainties, Roberti2018APlans, Pacino2018CranePlanning}. Furthermore, earlier work considers hatch restows (e.g., \cite{Shields1984, Kang2002StowageTransportation, Pacino2011FastVessels}), while stowage plans with crane operations are more frequently recurring (e.g., \cite{Hamedi2011CONTAINERSHIPOPERATIONS, Azevedo2018SolvingAlgorithm, Chang2022SolvingMode}). In less frequent cases, block stowage strategies are proposed by \cite{Wilson2000ContainerSolutions, Liu2011RandomizedPlans, Pacino2018CranePlanning}, and ballast water is modeled in \cite{Shields1984, Botter1992StowageSolution}. Lashing forces are only considered in \cite{Shields1984}.

The origin of most mathematical models in Table \ref{tab:full_MPClassification} is traced back to \cite{Botter1992StowageSolution}. The formulation contains two decision variables to load or remove container $c$ into/from a slot with bay $b$, row $r$ and tier $t$ at any port $p$ between loading port $l$ and discharge port $d$. This enables the modeling of voluntary restows. 
In the $k$-shift problem, both variables are integrated to stow uniform cargo on a box-formed vessel with a single bay \citep{AvrielMordecaiPenn1993, Avriel1998StowageShifts}. The decision variable indicates whether slot $(r,t)$ is occupied by a container at loading port $l$ with destination $d$ and planned discharging port $\upsilon$. Note that $\upsilon<d$ causes a restow to an arbitrary slot $(r',t')$ at port $\upsilon$, whereas $\upsilon=d$ results in zero restows. Recently, the problem is extended to multiple cargo weights and lengths for medium-sized vessels \citep{Parreno-Torres2020ImprovingProblems, Parreno-Torres2021SolvingAlgorithm}.
The work by \cite{Kang2002StowageTransportation} assigns cargo to blocks and subsequently creates stacks ordered by destination and weight, thereby reducing mandatory restows. This assignment formulation often scales to medium-sized vessels (e.g., \cite{Kaisar2006ATRANSPORTATION, Hamedi2011CONTAINERSHIPOPERATIONS, Pacino2011FastVessels}). It is extended to a stochastic program with uncertain container weights \citep{Li2018Multi-PortUncertainties} and a block stowage problem with crane intensity \citep{Pacino2018CranePlanning}.

Despite the many implementations in Table \ref{tab:full_MPClassification}, it is challenging to compare their performance adequately. On the one hand, little consensus exists on a general problem definition or benchmark instances. Hence, the problem is attempted to be solved from different angles. On the other hand, the $k$-shift problem is defined clearly and includes multiple benchmarks, but needs to be extended further to become representative.

An optimal solution to a representative multi-port problem is yet to be found. The early work of \cite{Botter1992StowageSolution} demonstrates this for their formulation, and therefore suggests a hierarchical decomposition that combines a 0-1 IP with an enumeration tree. Consequently, hybrids of exact and neighborhood-based methods \citep{Wilson2000ContainerSolutions, Kaisar2006ATRANSPORTATION, Pacino2011FastVessels} or greedy and tree-search heuristics \cite{Kang2002StowageTransportation} are proposed to sequentially solve master and slot planning. 
More recently, decompositions are substituted by heuristic frameworks, such as genetic algorithms \citep{Hamedi2011CONTAINERSHIPOPERATIONS, Hu2017Multi-objectiveRoute, Chang2022SolvingMode}, (large) neighborhood search  \citep{Pacino2018CranePlanning, Li2018Multi-PortUncertainties}, or a framework of greedy and tabu-search heuristics \citep{Liu2011RandomizedPlans}.

\begin{figure}[htp]
    \footnotesize
    \centering \input{plots/multi_port_solve_times.tex}
    \caption{Execution time over the number of ports in the tested instance \citep{Pacino2011FastVessels,Kang2002StowageTransportation}. In red is the linear trend line.}
    \label{fig:plot_multi_port_time}
\end{figure}
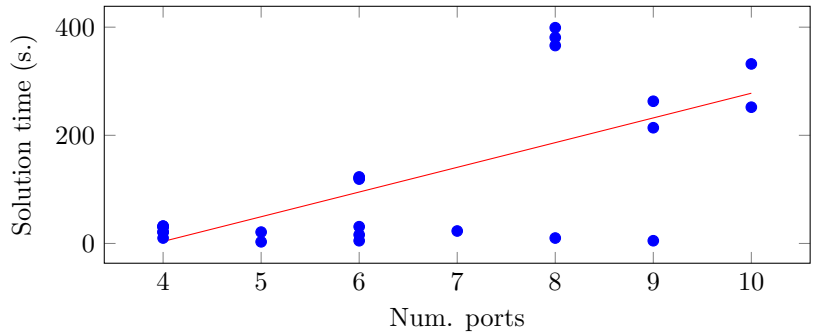

Early computational experiments have been published by \cite{Kang2002StowageTransportation}, who iteratively solve the problem at each port and add constraints to the MP to minimize restows. In \cite{Pacino2011FastVessels} we instead see a focus on larger problem instances and an increased complexity of SPP constraints. The decomposition is similar to the one proposed in \cite{Kang2002StowageTransportation}, but without iterative interaction between MBPP and SPP. 
Though the experiments are based on very different hardware, it is still possible to analyze the impact of vessel size and the number of planned ports on the computation time. Figure~\ref{fig:plot_multi_port_time} shows a plot of the reported execution time from the two papers. An analysis of the results shows a weak negative correlation between the vessel size and the execution time (-0.29). In contrast, a stronger correlation appears between the execution time and the number of ports (0.67). The red line in the figure shows the increasing linear trend of the results\footnote{The analysis was done discarding results where the algorithm timed out. Also, 2 results with extremely long runtimes were considered outliers.}. The collected data is shown in Table~\ref{tab:mutli_port_full_problem_results} of \ref{app:result_tables}. Due to single-instance computational studies, we could not include more experiments in this comparison.

Regarding the $k$-shift problem, construction heuristics are used in \cite{Avriel1998StowageShifts, Ding2015StowageShifts} to find upper bounds. Moreover, \cite{Roberti2018APlans} suggest a branch-and-price framework to solve the problem for medium-sized vessels, while \cite{Parreno-Torres2019SolutionProblem} formulates a 0-1 IP with valid inequalities to improve the linear relaxation. A compact formulation of this 0-1 IP is suggested in \cite{Parreno-Torres2020ImprovingProblems, Parreno-Torres2021SolvingAlgorithm}.
Despite these efforts, exact methods struggle with solving instances of large vessels. This motivates the use of heuristics, for which promising results are found by GRASP and matheuristics \citep{Parreno-Torres2019SolutionProblem, Parreno-Torres2020ImprovingProblems, Parreno-Torres2021SolvingAlgorithm}. Other heuristics such as genetic algorithms, beam search and simulated annealing can solve small vessel instances \citep{Dubrovsky2002AProblem, Azevedo2014SolvingMeta-heuristics, Azevedo2018SolvingAlgorithm}.

The first benchmarks are proposed by \cite{Avriel1998StowageShifts}, i.e. {\em Mixed}, {\em Long} and {\em Short} instances that refer to the expected time cargo is on-board. The {\em Authentic} instances by \cite{Ding2015StowageShifts} have ensured fully loaded vessels at each port, while the {\em Required} instances by \cite{Roberti2018APlans} have ensured the presence of shifts in optimal solutions. 

Table~\ref{tab:k-shift-results} summarizes the computational results of the exact methods for the $k$-shift problem, which contains formulations of \cite{Avriel1998StowageShifts}, \cite{Roberti2018APlans}, \cite{Parreno-Torres2019SolutionProblem} and \cite{Parreno-Torres2021SolvingAlgorithm}. The number of optimal solutions and upper bounds (feasible solutions) are presented for each method/formulation and instance group. All formulations except PAP can find a feasible upper bound for almost every instance. Moreover, the contributions of PCAR and RP find optimal solutions to almost each Long, Mixed, Short and Authentic instance, while APSW and PAP achieve optimal solutions in most of those instances. Nevertheless, only RP can find optimal solutions in 61 of the Required instances.

\begin{table}[h!]
    \footnotesize
    \caption{Results summary of exact methods for the $k$-shift problem. Column {Inst. group} indicates the instance group, while column {\em \#} shows the number of instances in each group. For each method, two columns are presented: {\em Opt.} indicating the number of optimal solutions found and {\em UB} indicating feasible solutions. The four methods are \cite{Avriel1998StowageShifts} (\em APSW), \cite{Parreno-Torres2019SolutionProblem} ({\em PAP}), \cite{Parreno-Torres2021SolvingAlgorithm} ({\em PCAR}), and \cite{Roberti2018APlans} ({\em RP}).}
    \label{tab:k-shift-results}
    \centering
    \input{tables/research_theory_table1.tex}
\end{table}

\begin{table}[h!]
\centering
\footnotesize
\caption{\label{tab:full_MPClassification}Classification of multi-port container stowage planning problems}
\scalebox{.72}{
\scriptsize
\begin{tabular}{llllllll}
\textbf{Paper} & \textbf{Cargo} & \textbf{Hydro} & \textbf{CSPP aspects}     & \textbf{Obj} & \textbf{Sc} & \textbf{HD} & \textbf{Solution methods}       \\ 
\hline
\cite{Shields1984}    & Class          & Rich           & MinRe, HR, BW, La, RF     & PS, VU, H    & M           &             & Greedy                          \\ 
\hline
\cite{Botter1992StowageSolution}     & Mix            & Rich           & VolRe, HR, BW         & PS           & S           & \checkmark           & Exact/TreeB                  \\ 
\hline
\cite{AvrielMordecaiPenn1993}     & Uni            & None           & VolRe                     & PS           & S           &             & Greedy                   \\ 
\hline
\cite{Avriel1998StowageShifts}     & Uni            & None           & VolRe                     & PS           & S           &             & Greedy                          \\ 
\hline
\cite{Wilson2000ContainerSolutions}     & Class          & None           & MinRe, HR, CO, BS, DG, RF     & PS, VU       & S           & \checkmark           & Exact/NeighMeta                  \\ 
\hline
\cite{Dubrovsky2002AProblem}  & Uni            & Stab           & VolRe                     & PS           & S           &             & PopulMeta                  \\ 
\hline
\cite{Kang2002StowageTransportation}       & Class          & Stab           & MinRe, HR                 & PS           & M           & \checkmark           & Greedy/TreeB                  \\ 
\hline
\cite{Kaisar2006ATRANSPORTATION}     & Class          & Equi           & MinRe, DG, RF         & PS           & M           & \checkmark           & Exact/NeighMeta                   \\ 
\hline
\cite{Hamedi2011CONTAINERSHIPOPERATIONS}     & Mix            & Stab           & MinRe, CO, DG, RF         & PS           & M           &             & PopulMeta                          \\ 
\hline
\cite{Liu2011RandomizedPlans}        & Class            & Stab           & MinRe, HR, CO, BS, DG, RF     & PS, VU, H    & M           &            & Greedy/NeighMeta                  \\ 
\hline
\cite{Pacino2011FastVessels}     & Mix            & Rich           & MinRe, HR, CO, RF         & PS           & M           & \checkmark           & Exact/NeighMeta                   \\ 
\hline
\cite{Azevedo2014SolvingMeta-heuristics}    & Uni            & Equi           & VolRe                     & PS, H        & S           &             & NeighMeta, TreeB, PopulMeta           \\ 
\hline
\cite{Ding2015StowageShifts}       & Uni            & None           & VolRe                     & PS           & M           &             & Greedy                          \\ 
\hline
\cite{Hu2017Multi-objectiveRoute}         & Mix            & Stab           & MinRe                     & PS, H        & S           &             & PopulMeta                  \\ 
\hline
\cite{Azevedo2018SolvingAlgorithm}    & Uni            & Equi           & VolRe, CO                 & PS           & S           &             & PopulMeta                          \\ 
\hline
\cite{Li2018Multi-PortUncertainties}         & Class          & Equi           & NARe                      & VU       & S           &             & NeighMeta, Exact                   \\ 
\hline
\cite{Pacino2018CranePlanning}     & Class          & None           & NARe, CO, BS              & PS           & L           &             & NeighMeta                          \\ 
\hline
\cite{Roberti2018APlans}    & Uni            & None           & VolRe, NARe                     & VU           & M           &             & Exact                           \\ 
\hline
\cite{Parreno-Torres2019SolutionProblem}    & Uni            & None           & VolRe                     & PS           & L           &             & Exact, NeighMeta                   \\ 
\hline
\cite{Parreno-Torres2020ImprovingProblems}     & Class          & Stab           & VolRe                     & PS           & L           &             & Exact, MatHeu, NeighMeta  \\ 
\hline
\cite{Parreno-Torres2021SolvingAlgorithm}    & Class          & Stab           & VolRe                     & PS           & M           &             & Exact, MatHeu  \\ 
\hline
\cite{Chang2022SolvingMode}      & Mix            & Stab           & MinRe, CO, (DG), RF                 & PS           & S           &             & PopulMeta                          \\
\hline
\end{tabular}
}
\end{table}

\subsection{Master planning}
\label{sec:master_planning}
\begin{table}
\centering
\caption{\label{tab:MBPP_Classification}Classification of master planning problems.}
\scalebox{.9}{
\scriptsize
\begin{tabular}{lllllll}
\textbf{Paper}    & \textbf{Cargo} & \textbf{Hydro} & \textbf{CSPP aspects} & \textbf{Obj} & \textbf{Sc} & \textbf{Solution methods}  \\ 
\hline
\cite{Pacino2012AnTanks}      & Class          & Rich           & BW, RF            & H            & L           & Exact                      \\ 
\hline
\cite{Pacino2013AnPlanning}        & Class          & Rich           & (MinRe), HR, CO, RF     & PS           & M           & NeighMeta                  \\ 
\hline
\cite{Ambrosino2015APlanning}     & Class          & Equi           & HR, CO, RF     & PS           & L           & MatHeu              \\ 
\hline
\cite{Ambrosino2015ComputationalProblems} & Class          & Equi           & HR, CO, RF     & PS           & L           & Exact                \\ 
\hline
\cite{Ambrosino2015ExperimentalProblem}  & Class          & Equi           & HR, CO         & PS, VU       & M           & Exact, MatHeu              \\ 
\hline
\cite{Ambrosino2018ShippingApproach}     & Class          & Rich           & HR, CO, RF     & PS           & L           & Exact, MatHeu              \\ 
\hline
\cite{Kebedow2018IncludingProblem}       & Class          & Rich           & HR, CO, DG, RF        & PS, VU       & L           & Exact                      \\ 
\hline
\cite{Bilican2020AParameters}       & Class          & Rich           & MinRe, HR, BW         & PS, H        & L           & Exact/NeighMeta                     \\ 
\hline
\cite{Chao2021}          & Class          & None           & HR             & PS, VU       & M           & Exact                      \\
\hline
\end{tabular}
}
\end{table}

Though it is an important part of the hierarchical decomposition that many researchers use, the SPP (when solved heuristically) has a minor impact on the runtime of the solution approach and its general objectives. For that reason, a number of scholars find it legitimate to focus their studies to solve the Master Planning Problem (MPP).

Table \ref{tab:MBPP_Classification} compares publications and shows trends in master planning research. In the MPP, cargo is grouped by weight and sorted into corresponding weight classes, since the knowledge of the true container weights is less important when building a master plan. While reefer containers are included in most test instances, there is only one study that handles dangerous cargo as well \citep{Kebedow2018IncludingProblem}. 

The problem formulations tend to include hydrostatic calculations of varying degrees, often more than just trim and GM estimations. \cite{Ambrosino2015ExperimentalProblem, Ambrosino2015APlanning, Ambrosino2015ComputationalProblems} based hydrostatics calculations on a bit older and less accurate model including equilibrium consideration that entails balancing weights on the vessel \citep{Ambrosino2004StowingProblem}. \cite{Ambrosino2018ShippingApproach} uses a richer stability model that includes trim, GM, and shear forces. \cite{Chao2021} did not consider stability constraints in their problem definition at all, but focused instead on a new IP formulation based on a minimum-cost flow problem with a multi-commodity network structure. The objective function contains the cost of the assignment of a container and the cost of using extra bays.

There are only a few studies that allow using ballast water to fix possible instabilities \citep{Pacino2012AnTanks, Bilican2020AParameters}. \cite{Pacino2012AnTanks} proposed an IP with an approximation of the displacement and a linearization of the centre of gravity calculations to include the effect of the ballast water on the hydrostatic values. 

The objective function mostly included balancing crane work and minimizing hatch restows, but some of the formulations focused also on maximizing vessel utilization in addition to port stay optimization \citep{Ambrosino2015ExperimentalProblem, Kebedow2018IncludingProblem, Chao2021}. While accounting for unnecessary hatch cover moves is considered in the MPP publications, there is rarely a focus on minimizing restows within blocks. Only \cite{Bilican2020AParameters, Pacino2013AnPlanning} presented related constraints.

From a mathematical modeling point of view, most models (\cite{Pacino2013AnPlanning,Ambrosino2018ShippingApproach, Kebedow2018IncludingProblem, Bilican2020AParameters}) are inspired by or extensions of the formulation proposed by \cite{Ambrosino2004StowingProblem} and \cite{Pacino2011FastVessels}. The only scientific work that has proposed a different formulation is that of \cite{Chao2021}, where the problem is modeled using a network-flow representation. 

A summary of the solution methods is shown in Table \ref{tab:MBPP_Classification}. \cite{Ambrosino2015APlanning, Ambrosino2018ShippingApproach} introduced the use of matheuristics that decompose the MIP. Firstly, an assignment of container destinations to blocks was found by solving the relaxed MIP with linearly relaxed variables indicating the assignment of containers to blocks. Secondly, a heuristic called the progressive random fixing procedure was used to obtain the feasible solution, where the assignment of container destinations to blocks is given from the previous phase and there is no relaxation of the MIP. \cite{Ambrosino2015ComputationalProblems, Ambrosino2015ExperimentalProblem} introduced two MIP formulations: the first one was a binary representation of the problem, and the second one expressed the number of containers in TEU and resulted in a more compact formulation. For the solution approach, the authors proposed to primarily solve a relaxation of the second model and then use its solution as input for solving the first model. The use of large neighborhood search was explored by \cite{Pacino2013AnPlanning} using the results from a relaxed IP as a warm start. \cite{Bilican2020AParameters} proposed a heuristic, where firstly a problem was solved with relaxed trim constraints and secondly a local search was performed to fix eventual instabilities. 

By looking at Table \ref{tab:MBPP_Classification}, we can observe that most of the proposed approaches were tested against large instances (above 15,000 containers). 

Though the problems and the formulations are very similar, the lack of a common benchmark has made it impossible to compare results. To remedy this situation, we propose a new set of publicly available benchmark instances based on the vessel data from \cite{Larsen2021AProblem}. For each of the 3 available vessels, 2 random instances are generated for each combination of ports $\{5,7,10\}$ percentage of cargo already on board $\{0\%,15\%,30\%\}$ and vessel utilization $\{60\%,70\%,80\%\}$. The cargo list simulates an ocean-going service where the long-haul leg (sailing from one region to another) guarantees the provided vessel utilization. The benchmark has a total of 162 instances, all of which are available at \url{https://doi.org/10.11583/DTU.22293412}. 

Given this new set of benchmark instances, it is now possible to compare the efficiency of the assignment-based  \citep{Pacino2011FastVessels} and the network-flow formulation \citep{Chao2021}. The network-flow formulation does not include any stability constraints as mentioned before and it enforces zero hatch overstowage. The assignment-based formulation has thus been adjusted to follow this problem definition. \cite{Chao2021} proposes an objective function that diverges from the makespan minimization from the literature. Unfortunately, its description is not accurate enough to be reproduced, hence, we modified the formulation to minimize the makespan as described in \cite{Pacino2011FastVessels}. We have also corrected the capacity constraints, since in the original work they were posted per container type and not for the block as a whole (the two formulations are presented in ~\ref{app:result_tables}).

\begin{table}[]
    \centering
    \scriptsize
    \caption{Runtime comparison of relaxation of the network-flow and the assignment-based formulation for the MPP without stability constraints.}
    \label{tab:MP_simple_comparison_relaxed}
    \input{tables/master_plan_simple_relaxed.tex}
\end{table}

\begin{table}[]
    \centering
    \scriptsize
    \caption{Runtime comparison of the network-flow and the assignment-based formulation for the MPP without stability constraints.}
    \label{tab:MP_simple_comparison}
    \input{tables/master_plan_simple.tex}
\end{table}

Table~\ref{tab:MP_simple_comparison_relaxed} shows the solution time comparison between the two formulations. The first three columns indicate the vessel's size, the number of ports, and the number of instances solved. Six instances are solved for each of these combinations, which represent all the instances in the benchmark for which the ROB condition is empty (needed due to the hatch overstow constraint). The next two columns are the optimal solution and the time needed to find it. Notice that the problem is solved in the relaxed version, where the decision variables are continuous. As can be seen, both formulations can find optimal solutions to this relaxed problem within a reasonable time for small and medium-sized vessels. The network-flow formulation, however, is clearly underperforming, possibly due to the increasing number of arcs needed as the number of ports increases. The network-flow formulation was not able to find feasible solutions to any of the largest instances (large vessel and 10 ports), while the assignment-based formulation timed out (3600 sec.) only for two of those instances. 

Table~\ref{tab:MP_simple_comparison} presents the results of the two formulations without the relaxation of the decision variables. The table includes an extra column representing the gap between the returned solution and the lower bound (MIP Gap). The results resemble those of Table~\ref{tab:MP_simple_comparison_relaxed}, where the formulations struggle as the instances increase in size.

Despite the fact that the network-flow formulation performs worse than the assignment-based formulation, further studying its application can be interesting as such formulations are well-studied for decomposition methods.

\begin{table}[]
    \centering
    \scriptsize
    \caption{Aggregated results of the assignment-based formulation to the complete problem with relaxed decision variables.}
    \label{tab:mp_full_relatex}
    \input{tables/master_plan_full_relaxed.tex}
\end{table}

\begin{table}[]
    \centering
    \scriptsize
    \caption{Aggregated results of the assignment-based formulation to the complete problem.}
    \label{tab:mp_full}
    \input{tables/master_plan_full.tex}
\end{table}

Tables~\ref{tab:mp_full_relatex} and \ref{tab:mp_full} report the results of the assignment-based formulation on a rich version of the CSPP, where stability and shear-force constraints are included. The formulation minimizes hatch overstowage and crane makespan. Both tables report aggregated results grouped by vessel size and the number of planned ports (columns 1 and 2). The tables also report the number of instances in each group (N. Inst.), the number of instances solved and which of those are optimal (Sol/Opt), mean values of the optimality gap reported by the solver (MIP Gap), the mean of the objective function (Obj.), the mean number of hatch overstows (HO), the mean of the makespan (MK), and the average runtime (Time).

It can be easily seen from the tables that as the vessel size and the number of ports increase, so does the time required by the solver to find solutions. For the version of the problem where the assignment variables are relaxed, it is possible to find optimal solutions for all the instances within the 3600 sec. time limit used in the experiments. Without the variable relaxation (Table~\ref{tab:mp_full}), it is hard to find feasible solutions even for instances with a small vessel. The full set of results can be found in ~\ref{app:result_tables}.

\subsection{Slot planning}\label{sec:solMethods4}

\begin{table}
\centering
\caption{\label{tab:SPP_Classification}Classification of slot planning problems.}
\scalebox{.9}{
\scriptsize
\begin{tabular}{lllllll}
\textbf{Paper} & \textbf{Cargo} & \textbf{Hydro} & \textbf{CSPP aspects} & \textbf{Obj} & \textbf{Sc} & \textbf{Solution methods}  \\ 
\hline
\cite{Pacino2010ABays}     & Mix            & None           & MinRe, RF             & PS, VU       & L           & NeighMeta                  \\ 
\hline
\cite{Delgado2012ABays}    & Mix            & None           & MinRe, RF         & PS, VU       & L           & Exact                      \\ 
\hline
\cite{Parreno2016AProblem}    & Mix            & None           & MinRe, DG, RF         & PS, VU       & M           & NeighMeta                  \\ 
\hline
\cite{Yifan2016Group-BayShip}      & Mix            & None           & MinRe                 & PS        & S           & Greedy/PopulMeta                  \\ 
\hline
\cite{Jin2019AnBay}        & Mix            & None           & MinRe                 & PS           & S           & Exact                      \\ 
\hline
\cite{Kebedow2019IncludingProblem}    & Mix            & None           & MinRe, DG, RF         & PS, VU       & M           & Exact                      \\ 
\hline
\cite{Korach2020MatheuristicsBays}     & Mix            & None           & MinRe, RF         & PS, VU       & L           & MatHeu                     \\ 
\hline
\cite{Rashed2021AVessels}     & Mix            & None           & MinRe, RF         & PS, VU       & L           & NeighMeta                  \\
\hline
\end{tabular}
}
\end{table}

Table \ref{tab:SPP_Classification} provides a summary of the literature concerning exclusively the SPP. As the table shows, some of the important constraints, e.g., hydrostatics or hatch restows are not present in the problem formulations. Since the output of the MPP is the input for the SPP and these constraints are already fulfilled in the MPP part of the problem, they can be ignored in this phase.

There is more consensus reached on the SPP than on the MPP. A plausible reason is that a common definition of the problem and a set of benchmark instances have been available since the publication of \cite{Delgado2009}. The proposed formulation was an inspiration to several researchers in this area. This is seen in Table \ref{tab:SPP_Classification}, where the definitions of the problem are quite uniform. \cite{Delgado2009} considered a wide spectrum of container types including reefers and highcubes, and important low-level constraints including stacking rules and capacity constraints. The objective function contained several aspects of the CSPP: minimizing port stay by avoiding unnecessary crane moves, consolidation by minimizing the number of used stacks and preserving reefer slots for reefer containers. The assumption is that all the containers from the load list can be stowed in the block. \cite{Parreno2016AProblem} modified the objective function based on the fact that not all containers can be stowed in their proposed formulation, so the aim is to load most of them. Since we need to solve the SPP for every block on the vessel separately, the computation time of the proposed approach has to be low, such that the whole process of creating a slot plan is finished in a reasonable time. A slot plan for one block should be created ideally in less than one second \citep{Delgado2012ABays}.

What is important to underline is that the model proposed by \cite{Delgado2009} considered only creating stowage plans for below deck stacks. This made it possible to ignore constraints related to lashing, line of sight and 45' containers. The study of \cite{Kebedow2019IncludingProblem, Jin2019AnBay} took the on deck section of the vessel into account by using respectively arbitrary height or weight limits to mimic, in a simple way, lashing constraints. 

\cite{Parreno2016AProblem, Kebedow2019IncludingProblem} introduced handling of hazardous cargo and segregation rules in the SPP. \cite{Jin2019AnBay, Yifan2016Group-BayShip} included port operations objectives by minimizing restow containers on the yard and crane moves while loading and discharging the vessel.

Table \ref{tab:SPP_Classification} shows the solution methods suggested for the SPP. \cite{Delgado2009, Delgado2012ABays} proposed constraint programming (CP) with the usage of, among others, symmetry-breaking constraints, and branching strategies to achieve better computations time. Constraint-based local search (CBLS) explored by \cite{Pacino2010ABays} and \cite{Parreno2016AProblem} proposed Greedy Randomized Adaptive Search Procedure (GRASP), with a randomized construction of initial solutions and a local search to make improvements. A hybrid method involving A* and a Genetic Algorithm (GA) was developed by \cite{Yifan2016Group-BayShip}. A* was used to find a feasible loading sequence and GA to find a feasible allocation of containers to slots. A fuzzy logic algorithm with a rule-based search was presented in \cite{Rashed2021AVessels}. Additionally, a matheuristic was developed by \cite{Korach2020MatheuristicsBays}. It combined a large neighborhood search with a mathematical solver to iteratively destroy and rebuilt parts of the solution.

Table~\ref{tab:slot_planning_cmp} shows a comparison between all the slot planning approaches that have adopted the benchmark taken from the work of  \cite{Delgado2009}. The first column indicated the instance group (we refer the reader to \cite{Delgado2012ABays} for a detailed description), and the second the number of instances in that group. Next, the table is divided into 6 sections each representing the results of a publication: CBLS is the constraint-based local search of \cite{Pacino2013FastSearch}, IP is the integer programming formulation of \cite{Delgado2012ABays} (with 10 seconds time limit), CP is the constraint programming model of \cite{Delgado2009} (with 10 seconds time limit), Fuzzy is the fuzzy logic approach of \cite{Rashed2021AVessels}, Matheuristic is the matheuristic approach of \cite{Korach2020MatheuristicsBays}, and GRASP is the GRASP approach of \cite{Parreno2016AProblem} (run for 1 second). For each of the publications, the table reports the percentage of feasible solutions (Sol) and optimal solutions (Opt) found, plus the time used to compute all the instances in the group. The best results are highlighted in bold.

\begin{table}[h]
    \caption{Slot planning methods summary}
    \label{tab:slot_planning_cmp}
    \footnotesize
    \centering
    \input{tables/slot_plnning_cmp.tex}
\end{table}

The results in Table~\ref{tab:slot_planning_cmp} cannot be fully compared as experiments have been run on different hardware. That said, the CPUs used in \cite{Korach2020MatheuristicsBays} and \cite{Rashed2021AVessels} are comparable, and even though the hardware used in \cite{Parreno2016AProblem} is older, it can be assumed that some improvement can be expected if run on modern machines. With this in mind, the table shows a clear improvement from the original work of \cite{Delgado2009,Pacino2013FastSearch}. Feasible solutions have been found for all instances, and only for a few instances, optimal solutions do not exist. Given these results, this set of benchmarks seems to have achieved its purpose. In \cite{Parreno2016AProblem} it was pointed out that the benchmark contains a very limited set of discharge ports, which reduces drastically the complexity of restows, and hence propose a more challenging set of instances and a revised version of the problem including a load maximization objective and the handling of dangerous goods. On the one hand, this new set of benchmark instances brings new challenges to the problem. On the other hand, it is less representative of the kind of instances that a slot planning problem will face when being part of a decomposition algorithm. In the latter case, it is to be expected that the master plan will ensure to have as many containers as possible with the same discharge port. The new benchmark from \cite{Parreno2016AProblem} (including the set of instances from \cite{Delgado2009}) can be found at \url{https://doi.org/10.11583/DTU.22284475}.

A different direction is taken by \cite{Kebedow2019IncludingStowage}, where the focus is on the modeling of a large set (compared to that of \cite{Parreno2016AProblem}) of rules for dangerous goods. Unfortunately, the publication did not present the mathematical formulation and did not present results on the original benchmark. The new instances generated in \cite{Kebedow2019IncludingProblem} can, however, be found at \url{https://doi.org/10.11583/DTU.22293991} and be used for future comparison.

\subsection{Computational complexity} \label{sec:computational_complexity}
Relatively little work has focused on the study of the computational complexity of the CSPP. The first study focused on the complexity that stability constraints such as metacentric limits had on single-stack (GM-OSOP) and multi-stack overstowage problems (GM-MSOP) with uniform cargo \citep{Aslidis1989CombinatorialProblems}. A polynomial time algorithm is proved to exist for the GM-OSOP (with a time complexity of $\mathcal{O}(m^2n^3)$, where $m$ and $n$ refer to ports and containers respectively), while the computational complexity of the GM-MSOP is conjectured to be NP-Complete. An extension of this work is provided by \cite{Avriel2000ContainerGraphs}, which presented an NP-Completeness proof based on a reduction from the $C$-coloring problem of circle graphs, where $C$ represents the number of uncapacitated stacks (or the colors of the graph). The authors also proved that a polynomial time algorithm exists for $C<4$, and provided an algorithm to calculate upper and lower bounds on the number of stacks needed to find a solution with zero shifts. Further research \citep{Tierney2014OnProblems} showed that the capacitated version of the {\em k}-shift problem described in \cite{Avriel2000ContainerGraphs} can be solved in polynomial time for a fixed-sized vessel. The exponent in the polynomial is too large for any practical use, but the proof can be used to demonstrate that conclusions over experimental results conducted on a single vessel are not representative of the problem's complexity. Furthermore, \cite{Tierney2014OnProblems} study the computational complexity of the Hatch Overstow Problem (HOP) and show that the assignment of containers over and below and hatch cover with at most $k$ hatch overstows is NP-Complete by reduction from the set covering problem.  

\subsection{Other relevant publications}
\label{sec:misc}
Some studies have focused on other aspects of stowage optimization than solving the CSPP. They are presented in the following subsection.

An extension of the MPP to a selection problem was introduced by \cite{Christensen2017AStowage, Christensen2019AProblem, Kebedow2019IncludingStowage}. They considered a revenue management problem called cargo mix, where the goal was to select which bookings to accept in each port of call to maximize profit. A matheuristic was proposed in \cite{Christensen2017AStowage} composed of 3 stages: generating schedules of discharge ports by solving the longest path problem in an acyclic-directed graph, solving relaxed MIP where stability constraints are dropped and the final stage was fixing hydrostatics by possibly removing cargo. Stochastic programming was proposed by \cite{Christensen2019AProblem} considering uncertain demand per port, rolling horizon heuristic was introduced that decomposed the problem into sub-problems with shorter planning horizons.

The work of \cite{Jensen2018TheCapacity, Ajspur2019AModels, Jensen2022RevenueChallenge} introduced the Standard Capacity Model (SCM). It is a polyhedron model derived from MPP models and contributes the first linear approximations of hydrostatic equilibriums and restows. The purpose of the SCM is to increase the accuracy of cargo network-flow models such as \cite{Zurheide2015RevenueIndustry}, while maintaining their scalability. \cite{Jensen2022RevenueChallenge} applied the SCM to a yield optimization problem over 90 days in 2018 of Maersk's Asia - Europe service network with over 250 port calls. Optimal results could be computed in less than 30 minutes and showed that simple fixed capacity models used by carriers today can overestimate revenue by more than 20\%.   

The study presented in \cite{Lee2020LashingContainerships} proposed a multimodal deep learning model to predict the expected lashing forces for container stowage plans. Calculations for lashing forces are tedious and for this reason are hard to be incorporated into models of the CSPP. With the use of machine learning, this process could be faster and the presented results were promising, i.e., the average gap between predicted and true values was 0.66\%.

Two interactive decision support tools were presented with the usage of Binary Decision Diagrams (BDDs) \cite{Jensen2012FastDiagrams} and boolean satisfiability (SAT) \cite{Kroer2012SATPlanning}. The comparison of both methods was presented in \cite{Kroer2016SymbolicPlanning}. The software allowed for marking infeasible areas in a bay, but also suggestions of slots in which containers could be placed, and vice versa for containers and potential slots. The BDDs performed well in real-life instances.

%% file: plots/multi_port_solve_times.tex
\begin{tikzpicture}
\begin{axis}[
    ylabel={Solution time (s.)},
    xlabel={Num. ports},
    height=5cm,
    width=.8\textwidth
    ]
    \addplot [red, domain=4:10 ] {45.688 * x - 179.05};
    \addplot [blue, only marks] table [x=ports, y=time]{plots/multi_port_solve_times.csv};
\end{axis}
\end{tikzpicture}

%% file: tables/research_theory_table1.tex
\begin{tabular}{lrrrrrrrrrr}
            &      & \multicolumn{2}{c}{APSW}  &\multicolumn{2}{c}{PAP}  &\multicolumn{2}{c}{PCAR}  & \multicolumn{2}{c}{RP}\\\hline
Inst. Group & \#   & Opt.  & UB                & Opt.  & UB               & Opt.  & UB               & Opt.  & UB  \\\hline
Short       & 81   & 78    & 81                & 79    & 79               & 81    & 81               & 81    & 81 \\
Mixed       & 81   & 78    & 81                & 77    & 77               & 81    & 81               & 81    & 81 \\
Authentic   & 81   & 77    & 81                & 72    & 81               & 79    & 81               & 81    & 81 \\
Long        & 81   & 76    & 81                & 76    & 81               & 81    & 81               & 81    & 81 \\ 
Required    & 81   &  1    & 81                & 13    & 60               & 14    & 81               & 61    & 79 \\\hline
\end{tabular}

%% file: tables/master_plan_simple_relaxed.tex
\begin{tabular}{lc|rr|rr}
\multicolumn{2}{c|}{\textbf{}}           & \multicolumn{2}{c|}{\textbf{\cite{Chao2021}}} & \multicolumn{2}{c}{\textbf{\cite{Pacino2011FastVessels}}} \\ \hline
\textbf{Vessel}         & \textbf{Ports} & \textbf{Objective} & \textbf{Time} & \textbf{Objective}  & \textbf{Time} \\ \hline
\multirow{3}{*}{Small}  & 5              & 1544.33            & 5.41          & 1544.33             & 2.48          \\
                        & 7              & 1922.61            & 21.15         & 1922.61             & 9.77          \\
                        & 10             & 2326.74            & 588.16        & 2326.74             & 47.04         \\ \hline
\multirow{3}{*}{Medium} & 5              & 1936.94            & 6.16          & 1936.94             & 3.25          \\
                        & 7              & 2425.49            & 27.86         & 2425.49             & 12.78         \\
                        & 10             & 2326.74            & 588.16        & 2326.74             & 47.04         \\ \hline
\multirow{3}{*}{Large}  & 5              & 2713.70            & 53.16         & 2713.68             & 16.62         \\
                        & 7              & 3376.71            & 1406.08       & 3376.60             & 450.90        \\
                        & 10             & -                  & -             & 3821.90             & 756.61        \\ \hline
\end{tabular}

%% file: tables/master_plan_simple.tex
\begin{tabular}{lc|rrr|rrr}
\multicolumn{2}{c|}{\textbf{}}           & \multicolumn{3}{c|}{\textbf{\cite{Chao2021}}}                    & \multicolumn{3}{c}{\textbf{\cite{Pacino2011FastVessels}}}   \\ \hline
\textbf{Vessel}         & \textbf{Ports} & \textbf{Objective} & \textbf{MIP Gap} & \textbf{Time} & \textbf{Objective} & MIP Gap & Time    \\ \hline
\multirow{3}{*}{Small}  & 5              & 1546.17            & 0\%             & 22.51         & 1546.17            & 0\%    & 22.81   \\
                        & 7              & 1925.83            & 0\%             & 344.28        & 1925.83            & 0\%    & 77.97   \\
                        & 10             & 2328.83            & 0\%             & 2295.48       & 2329.00            & 0\%    & 1254.14 \\ \hline
\multirow{3}{*}{Medium} & 5              & 1939.33            & 0\%             & 65.80         & 1939.33            & 0\%    & 13.54   \\
                        & 7              & 2428.50            & 0\%             & 117.70        & 2428.50            & 0\%    & 70.60   \\
                        & 10             & 2328.83            & 0\%             & 2295.48       & 2329.00            & 0\%    & 1254.14 \\ \hline
\multirow{3}{*}{Large}  & 5              & 3063.00            & 9\%             & 3202.49       & 2726.00            & 0\%    & 3023.26 \\
                        & 7              & 4750.00            & 21\%             & 3604.85       & 3445.17            & 2\%    & 3600.04 \\
                        & 10             & 6422.00            & 38\%             & 3600.14       & 4210.00            & 5\%    & 3600.09 \\ \hline
\end{tabular}

%% file: tables/master_plan_full_relaxed.tex
\begin{tabular}{cccrrrrrr}
\textbf{Vessel}         & \textbf{Ports} & \textbf{N. Inst.} & \textbf{Sol/Opt} & \textbf{MIP Gap} & \textbf{Obj.} & \textbf{HO} & \textbf{MK} & \textbf{Time} \\ \hline
\multirow{3}{*}{Small}  & 5              & 18                & 18(18)           & 0\%              & 3497.69       & 4.00        & 1575.56     & 12.04         \\
                        & 7              & 18                & 18(18)           & 0\%              & 4752.69       & 7.61        & 1936.18     & 42.82         \\
                        & 10             & 18                & 18(18)           & 0\%              & 4751.59       & 5.87        & 2340.71     & 192.68        \\ \hline
\multirow{3}{*}{Medium} & 5              & 18                & 18(18)           & 0\%              & 4820.40       & 6.06        & 1964.94     & 17.15         \\
                        & 7              & 18                & 18(18)           & 0\%              & 6260.96       & 8.61        & 2433.34     & 70.09         \\
                        & 10             & 18                & 18(18)           & 0\%              & 5843.55       & 8.16        & 2943.79     & 312.38        \\ \hline
\multirow{3}{*}{Large}  & 5              & 18                & 18(18)           & 0\%              & 6598.71       & 7.39        & 2782.16     & 28.93         \\
                        & 7              & 18                & 18(18)           & 0\%              & 8440.38       & 7.44        & 3434.99     & 163.23        \\
                        & 10             & 18                & 18(18)           & 0\%              & 6862.35       & 6.33        & 4151.48     & 1301.05       \\ \hline
\end{tabular}

%% file: tables/master_plan_full.tex
\begin{tabular}{llcrlllll}
\multicolumn{1}{c}{\textbf{Vessel}} & \multicolumn{1}{c}{\textbf{Ports}} & \textbf{N. Inst.} & \textbf{Sol/Opt} & \multicolumn{1}{r}{\textbf{MIP Gap}} & \multicolumn{1}{r}{\textbf{Obj.}} & \multicolumn{1}{r}{\textbf{HO}} & \multicolumn{1}{r}{\textbf{MK}} & \multicolumn{1}{r}{\textbf{Time}} \\ \hline
\multirow{3}{*}{Small}              & 5                                  & 18                & 18(14)           & 0\%                                  & 3501.44                           & 6.39                            & 1579.22                         & 1228.37                           \\
                                    & 7                                  & 18                & 17(8)            & 1\%                                  & 4532.12                           & 8.18                            & 1943.88                         & 2923.58                           \\
                                    & 10                                 & 18                & 0(0)             & \multicolumn{1}{c}{-}                & \multicolumn{1}{c}{-}             & \multicolumn{1}{c}{-}           & \multicolumn{1}{c}{-}           & \multicolumn{1}{c}{-}             \\ \hline
\multirow{3}{*}{Medium}             & 5                                  & 18                & 15(14)           & 0\%                                  & 4691.73                           & 5.80                            & 1905.07                         & 1243.34                           \\
                                    & 7                                  & 18                & 10(4)            & 25\%                                 & 17796.50                          & 37.60                           & 2626.50                         & 3076.16                           \\
                                    & 10                                 & 18                & 2(0)             & 19\%                                 & 6936.50                           & 10.00                           & 2736.50                         & 3600.16                           \\ \hline
\multirow{3}{*}{Large}              & 5                                  & 18                & 9(0)             & 1\%                                  & 6550.44                           & 8.44                            & 2583.78                         & 3600.14                           \\
                                    & 7                                  & 18                & 0(0)             & \multicolumn{1}{c}{-}                & \multicolumn{1}{c}{-}             & \multicolumn{1}{c}{-}           & \multicolumn{1}{c}{-}           & \multicolumn{1}{c}{-}             \\
                                    & 10                                 & 18                & 7(0)             & 82\%                                 & 108633.00                         & 262.71                          & 5290.14                         & 3600.27                           \\ \hline
\end{tabular}

%% file: tables/slot_plnning_cmp.tex
\begin{tabular}{cr|rrr|rrr|rrr}
\multicolumn{2}{c}{}  & \multicolumn{3}{c}{CBLS}             & \multicolumn{3}{c}{IP (10s)}         & \multicolumn{3}{c}{CP(10s)}           \\\hline
Group & Inst & Sol           & Opt          & Time  & Sol          & Opt          & Time   & Sol           & Opt           & Time  \\\hline
1     & 13   & \textbf{100}  & 59           & 0.10  & \textbf{100} & \textbf{100} & 1.80   & \textbf{100}  & \textbf{100}  & 0.10  \\
2     & 22   & \textbf{100}  & 77           & 3.60  & 95           & 91           & 50.40  & 31            & 91            & 21.60 \\
3     & 13   & \textbf{100}  & 92           & 0.50  & 92           & 85           & 35.30  & \textbf{100}  & \textbf{100}  & 0.50  \\
4     & 78   & \textbf{100}  & 92           & 6.00  & 96           & 94           & 87.00  & 99            & 99            & 19.70 \\
5     & 36   & 97            & 58           & 7.10  & 72           & 56           & 192.00 & 92            & 92            & 39.00 \\
6     & 15   & 93            & 80           & 1.20  & \textbf{100} & 93           & 13.00  & \textbf{100}  & \textbf{100}  & 5.40  \\
7     & 14   & 93            & 79           & 2.30  & 64           & 29           & 102.80 & 64            & 64            & 53.50 \\
8     & 14   & 93            & 43           & 1.50  & 79           & 64           & 74.10  & 93            & 93            & 10.50 \\
9     & 17   & 94            & 47           & 5.20  & 53           & 41           & 112.30 & 88            & \textbf{88}   & 36.50 \\
10    & 8    & \textbf{100}  & 88           & 0.70  & 88           & 62           & 31.50  & \textbf{100}  & \textbf{100}  & 0.70  \\
11    & 6    & 50            & 17           & 1.30  & 67           & 50           & 30.50  & 83            & \textbf{83}   & 10.30 \\\hline
 \multicolumn{11}{c}{}\\
\multicolumn{2}{c}{}      & \multicolumn{3}{c}{Fuzzy}            & \multicolumn{3}{c}{Matheuristic}     & \multicolumn{3}{c}{GRASP (1s)}        \\\hline
Group & Inst & Sol           & Opt          & Time  & Sol          & Opt          & Time   & Sol           & Opt           & Time  \\\hline
1     & 13   & \textbf{100}  & \textbf{100} & 4.52  & \textbf{100} & \textbf{100} & 1.30   & \textbf{100}  & \textbf{100}  & 6.30  \\
2     & 22   & \textbf{100}  & 95           & 8.12  & 95           & 86           & 12.10  & \textbf{100}  & \textbf{100}  & 15.40 \\
3     & 13   & \textbf{100}  & \textbf{100} & 5.10  & 92           & 92           & 8.00   & \textbf{100}  & \textbf{100}  & 7.80  \\
4     & 78   & \textbf{100}  & \textbf{100} & 28.32 & \textbf{100} & \textbf{100} & 17.70  & \textbf{100}  & 99            & 37.30 \\
5     & 36   & \textbf{100}  & \textbf{94}  & 16.41 & 89           & 83           & 29.40  & \textbf{100}  & \textbf{94}   & 22.30 \\
6     & 15   & \textbf{100}  & 93           & 4.68  & \textbf{100} & \textbf{100} & 2.50   & \textbf{100}  & \textbf{100}  & 6.00  \\
7     & 14   & \textbf{100}  & 71           & 6.66  & 86           & 71           & 11.00  & \textbf{100}  & \textbf{93}   & 9.00  \\
8     & 14   & \textbf{100}  & \textbf{100} & 6.40  & \textbf{100} & 79           & 5.40   & \textbf{100}  & \textbf{100}  & 6.10  \\
9     & 17   & \textbf{100}  & 82           & 8.69  & 94           & 76           & 13.00  & 94            & 82            & 12.30 \\
10    & 8    & 94            & 88           & 3.33  & \textbf{100} & 88           & 2.50   & \textbf{100}  & \textbf{100}  & 4.80  \\
11    & 6    & \textbf{100}  & 67           & 2.56  & 83           & \textbf{83}  & 3.30   & \textbf{100}  & \textbf{83}   & 3.20  \\\hline
\end{tabular}

%% file: chapters/05_research_agenda.tex
In light of this review, we present our conclusions on the state-of-the-art and propose possible areas of future research. We will do so by starting to describe the challenges with respect to problem representation.
As mentioned in Section \ref{sec:stowagePlanning}, the included combinatorial aspects should be representative of the real-world problem. From Section~\ref{sec:contributions} it can be derived that a subset of these aspects has not been studied sufficiently. For instance, we believe that lashing is only modeled by dynamic stack capacity in \cite{ElYaagoubi2022Multi-objectiveSystem}, while \cite{Shields1984} only commented on their approach. In practice, the proper use of lashing rods can significantly increase on-deck capacity, while the effects of different lashing models are yet to be investigated. In addition, a substantial body of work implements voluntary restows to reduce restows at future ports (e.g., \cite{Botter1992StowageSolution, Avriel1998StowageShifts, Roberti2018APlans}). Nevertheless, it remains unclear to what extent these impact port stay and vessel utilization. Similarly, block stowage patterns are limited to block purity in \cite{Wilson2000ContainerSolutions, Liu2011RandomizedPlans, Pacino2018CranePlanning, Larsen2021AProblem}, even though more sophisticated patterns (e.g., paired block stowage) have been adopted by the industry. Thus, future work should also investigate these best practices.

Despite the individual cases, we assess that the interaction between key combinatorial aspects is studied insufficiently. Overall, each additional constraint reduces the capacity at ports or vessels, but how these interact and jointly impact the objectives should be investigated further. To do so, a fully representative model is necessary, which is yet to be modeled for the CSPP and MBPP. With respect to the SPP, most combinatorial aspects, except lashing, have been modeled adequately.

\begin{figure}[h]
    \centering
    \input{figures/tcg_error.tex}
    \caption{Example calculation of the transversal center of gravity. According to the balance constraints from \cite{Ambrosino2004StowingProblem} (where the total weight at each side of the center should be equal) the presented example is assumed balances ($tcg=0$) while it clearly is not the case.}
    \label{fig:CG_error}
\end{figure}
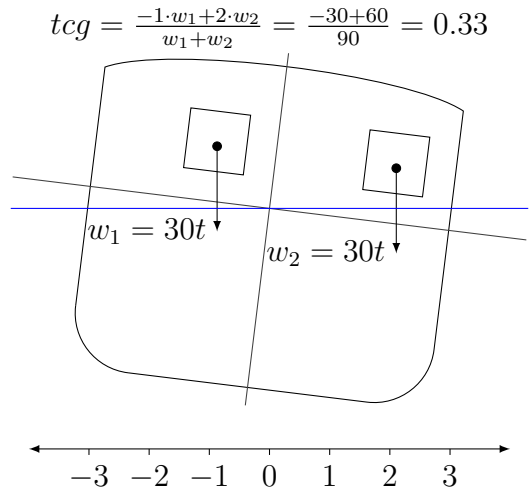

We suggest the following minimum requirements for a representative problem. The cargo model should consider 20/40 ft. lengths, standard and highcube heights, but also special cargo such as reefers and IMDG. In addition, voluntary and hatch restows represent reality well, whereas a combination of GM, trim, list and stress forces must model vessel stability. Any future work worth publishing should be aware of the issue raised from modeling stability constraints as simple balancing of weight (e.g., \cite{Ambrosino2004StowingProblem}). 
Figure~\ref{fig:CG_error} shows a simple example of how such simplifications are too far from reality and cannot be used. Since then, many scientific studies reverted to the use of a center of gravity calculation (i.e. \cite{Ambrosino2015APlanning, Pacino2011FastVessels, Christensen2017AStowage, Zhu2020IntegerProblem}), while unfortunately, some remain oblivious to this mistake (e.g., disregarding container weights \cite{Azevedo2014SolvingMeta-heuristics, Azevedo2018SolvingAlgorithm}, or disregarding the position of the weights \cite{Li2018Multi-PortUncertainties,Li2020OptimizingShipping,Cruz-Reyes2015,Kebedow2019IncludingProblem}).
Furthermore, lashing forces should be included as they impact on-deck stack capacity, while incorporating crane operations and block stowage enables to evaluate and enhance (un)loading efficiency. The main objectives are to minimize port stay and maximize utilization on at least 15,000 TEU vessels.

\subsection{Solution methods}
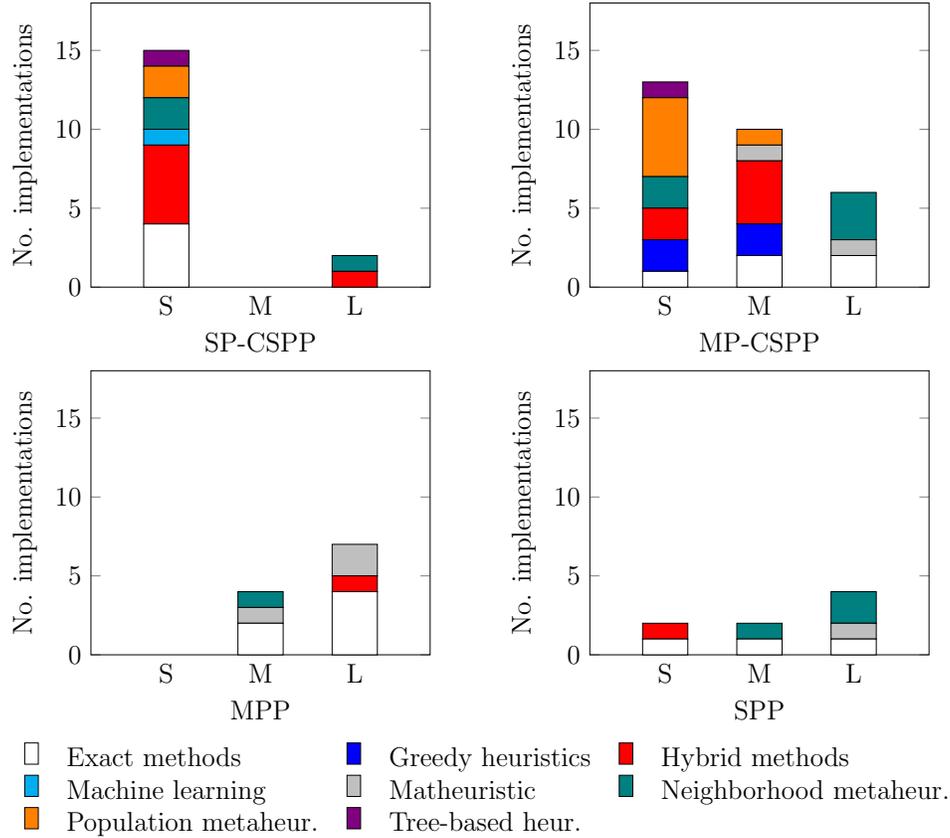
\begin{figure}[htp]
    \centering
    \resizebox{.95\linewidth}{!}{\input{plots/solution_methods_barchart.tex}}
    \caption{\label{fig:sol_methods_bar}Number of implemented solution methods for single-port (SP-CSPP), multi-port (MP-CSPP), master bay planning (MPP) and slot planning problems (SPP) with varying vessel sizes as defined in Table \ref{tab:scheme} (S=Small size, M=Medium size, L=Large size).}
\end{figure}

As representative problems are scarce, the underlying problem can vary greatly. Hence, we should tread carefully before drawing any conclusions from this comparison. Even though plenty of solution methods are proposed in Figure \ref{fig:sol_methods_bar}, their experiments are often limited. In order to verify their generalizability, implementations should strive for computational studies with multiple realistic instances (e.g., \cite{Kang2002StowageTransportation, Pacino2011FastVessels, Parreno-Torres2019SolutionProblem}). The use or extension of benchmark instances enables such comparative studies (e.g., \cite{Avriel1998StowageShifts, Ding2015StowageShifts, Larsen2021AProblem}). Consequently, these studies will help us to find adequate solution methods.  
 
As in many operations research studies, articles focus on providing new mathematical formulations to a problem (e.g., \cite{Botter1992StowageSolution, Ambrosino2004StowingProblem, Delgado2012ABays}), or using those formulations to evaluate the efficiency of heuristic solution methods (e.g., \cite{Ambrosino2010AnProblem, Pacino2013AnPlanning, Korach2020MatheuristicsBays}). Unfortunately, we have seen articles on the container stowage planning problem that either do not properly cite the origin of a mathematical formulation, or even present it as their own with only minor changes (if any). We encourage future authors and reviewers to be more critical, so that the literature is not overwhelmed with minor contributions that do now enhance the state-of-the-art.


\subsection{Future work}
In contrast to the work on the $k$-shift problem, very little consensus can be found on a common definition of the CSPP, a set of benchmark instances or the existence of a research road map. Most of the issues related to this lack of coordination can be attributed to the lack of publicly available data, and the high knowledge entry lever required to truly understand the calculations behind the vessel stability constraints. It is only recently that a textbook detailing the CSPP has been published \citep{Jensen2018ContainerPlanning}. Researchers that were lucky enough to collaborate with the industry, were constrained by non-confidentiality agreements from publishing details of their results (e.g., \cite{Wilson2001ContainerStudy}), or from making available the benchmark data (e.g., \cite{Pacino2011FastVessels}).

In this subsection, we suggest future work for each of the areas of research.

\subsubsection{Single-Port Container Stowage Planning}
The computational results of the single-port stowage planning (e.g., \cite{Cho1981DevelopmentPlanning, Ambrosino2010AnProblem}) are positive in terms of solution quality and computational efficiency. The impact the procedure has on today's large vessels, however, needs to be better evaluated. Most approaches are tested on small vessels for which the repositioning of a single container can have a significant effect on stability, which is no longer the case for the large vessels the industry now uses. 

The KPIs mentioned by \cite{Larsen2021AProblem} are interesting when compared to some of the model enhancements presented in \cite{Zhu2020IntegerProblem}, where it was argued that containers should be stowed tier-wise rather than having tall stacks. In contrast, \cite{Pacino2010ABays, Delgado2012AStowage, Larsen2021AProblem} argued that leaving free stacks provides a flexible stowage plan for future ports. As proposed by \cite{Larsen2021AProblem}, validation of such KPIs using simulation approaches is necessary. The vessel data mentioned in Subsection \ref{sec:solMethods1}, though simplified by the authors, still presents itself with a high learning curve. Thus, it is advisable that papers studying a specific version of the CSPP derive simplified data instances. An example can be seen in Subsection~\ref{sec:master_planning}, where a benchmark for the multi-port master planning problem is provided. 

Research on the single-port container stowage planning problem is far from finished, and we see the following as important future research directions: the design of exact methods for the identification of optimal solutions, the evaluation of the validity and usefulness of the proposed KPIs, and the evaluation of the use of other heuristic methods, e.g., based on container exchanges as proposed in \cite{Cho1981DevelopmentPlanning,Ambrosino2006AProblem}. Given the currently available data and experimental results, new research that does not include a full set of stability constraints is no longer of scientific interest.

\subsubsection{Multi-Port Container Stowage Planning}
With respect to the benchmarks for multi-port container stowage planning, given that the solution approach of \cite{Roberti2018APlans} provides optimal solutions within a minute, it can be concluded that the Long, Mixed, Short and Authentic instances of Table \ref{tab:k-shift-results} are now closed. As the formulation of \cite{Parreno-Torres2021SolvingAlgorithm} is able to find feasible solutions to all instances, further research on heuristics for the $k$-shift problem does not seem to be a valuable future direction any longer.

\cite{Parreno-Torres2021SolvingAlgorithm} arrives at the same conclusion and hence proposes to extend the $k$-shift problem with variable cargo sizes and simple stability constraints. It is shown that the additional complexity has a negative impact on the IP formulation and hence a matheuristic approach is proposed. It is likely that similar results could be obtained by extending the work of \cite{Roberti2018APlans} as the stability constraints would increase the number of constraints posted across the generated columns in the formulation and hence are likely to worsen the quality of the lower-bound found by the column generation within the approach.

As future research directions, we propose the study of exact and heuristic methods for the $k$-shift CSPP with simple stability constraints. This problem corresponds to the definition provided by \cite{Parreno-Torres2021SolvingAlgorithm}, where container types, weights and simple stability constraints are added to the original $k$-shift problem. It is unclear from the results presented in \cite{Parreno-Torres2021SolvingAlgorithm}, whether instances based on the Short, Mixed and Long transport matrices will result in any mandatory shifts, hence merit can be given also to future research  that studies, or leverages, the special case of the zero-shift problem.

To the best of the author's knowledge, limited progress has been made on the identification of a single-phase heuristic procedure for the multi-port container stowage planning problem. 
The vessel data provided by \cite{Larsen2021AProblem}, in combination with the cargo lists which we will describe in the next section, could be used as a common benchmark for future research. It is also the authors' opinion that future studies on exact methods are better suited as extensions of the $k$-shift problem (see \cite{Parreno-Torres2021SolvingAlgorithm}). 

\subsubsection{Master Planning}
Research on master planning is far from concluded. As a part of a hierarchical decomposition, master planning is most often solved using a relaxation of a mixed-integer programming formulation (e.g., \cite{Pacino2011FastVessels,Chou2021ApplyingStudy}). Though this has positive outcomes, the method is far from infallible, and its performance is heavily dependent on the features of the specific instance and on the combinatorial aspects included in the problem. 

From a problem representation point of view, combinatorial aspects such as block stowage and paired block stowage should be better studied. Only a few works have studied the impact of such stowage patterns on the achieved solutions and the performance of the solution methods (e.g., \cite{Wilson2000ContainerSolutions, Liu2011RandomizedPlans, Pacino2018CranePlanning, Larsen2021AProblem}). Though the use of mathematical modeling has the flexibility of easily allowing additional side constraints to the problem, research on heuristic methods with more stable performance should also be carried out. 

It is our hope that the new set of benchmark instances provided in this article (see Section~\ref{sec:master_planning}) will increase the quality and quantity of research on this problem.

\subsubsection{Slot Planning}
Thanks to the publicly available benchmarks, slot planning reached a high level of quality, and the problem, as currently defined, is (at least from an industrial point of view) solved. The benchmark, however, focuses on the SPP specific to below deck blocks. Aspects such as lashing forces have not been explored yet. 

Lashing forces are particularly interesting as little knowledge is currently available. The position of the container on deck, not only depends on its weight and the general load condition, but also on the type of lashing equipment available on the vessel. To which degree the mechanical calculation of the lashing forces can be simplified, and which assumption can be made to better implement solution algorithms is a field yet unexplored. The inclusion of lashing constraints is an important part of stowage planning, as a miscalculation might disallow an entire tier of a container from being loaded.

Being part of a hierarchical decomposition, slot planning has dependencies on the solution of the master planning problem. As of now, it is assumed that a master planning solution always generates feasible slot planning problems. In reality, this is not true (as shown by \cite{Pacino2011FastVessels}). Hence, slot planning could be extended to include the entire vessel, thereby allowing for the flexible assignment of containers to exchange between blocks and thus improving the solution quality. 

Another interesting extension of the slot planning problem is the integration with terminal operations. Some researchers have already realized this potential \citep{Monaco2014TheProblem, Iris2018FlexibleScheduling}, where the individual assignment of containers to container types is optimized with respect to the position of the cargo in the terminal. Other possible integrations include i.e. quay crane assignment and scheduling and container sequencing.

%% file: figures/tcg_error.tex
\begin{tikzpicture}
\newcommand{\figTcgAng}{7}
\begin{scope}[scale=0.8]
\begin{scope}[rotate around ={-\figTcgAng:(0,2)}]
\draw (-3,0) -- (-3,4) .. controls (-2,4.5) and (2,4.5) .. (3,4) -- (3,0) arc (0:-90:1) -- (-2,-1);
\draw (-3,0) arc (0:90:-1);
\draw [darkgray] (0,-1.3) -- (0,4.6);
\draw [darkgray] (-4.3,2) -- (4.3,2);
\draw (-1.5,3.5) rectangle (-0.5,2.5);
\draw (1.5,3.5) rectangle (2.5,2.5);
\draw [Circle-latex, rotate around={\figTcgAng:(-1,3)}] (-1,3) -- (-1,1.5) node [anchor=east] {$w_1=30t$};
\draw [Circle-latex, rotate around={\figTcgAng:(2,3)}] (2,3) -- (2,1.5) node [anchor=east] {$w_2=30t$};
\end{scope}
\draw [blue] (-4.3,2) -- (4.3,2);
\draw[latex-latex] (-4,-2) -- (4,-2);
\foreach \x in {-3,-2,-1,0,1,2,3}
  \draw (\x,-2) -- (\x,-2.1) node[anchor=north] {$\x$};
\node at (0,5) {$tcg = \frac{-1\cdot w_1 + 2\cdot w_2}{w_1+w_2} = \frac{-30 + 60}{90} = 0.33$};
\end{scope}
\end{tikzpicture}

%% file: plots/solution_methods_barchart.tex
\pgfplotstableread{
0    4	0	5	1	0	2	2	1
1    0	0	0	0	0	0	0	0
2    0	0	1	0	0	1	0	0
3    1	2	2	0	0	2	5	1
4    2	2	4	0	1	0	1	0
5    2	0	0	0	1	3	0	0
6    0	0	0	0	0	0	0	0
7    2	0	0	0	1	1	0	0
8    4	0	1	0	2	0	0	0
9    1	0	1	0	0	0	0	0
10   1	0	0	0	0	1	0	0
11   1	0	0	0	1	2	0	0
}\datatable	

\begin{tabular}{@{}cc@{}}
\begin{tikzpicture}
\begin{axis}[
    height=6cm,
    width=.5\textwidth,
    ylabel=No. implementations,
    xtick=data,
    xticklabels={S,M,L},
    x tick style={draw opacity=0}, 
    xlabel={SP-CSPP},
    enlarge y limits=false,
    enlarge x limits=0.4,
    ymin=0,ymax=18,
    ybar stacked,
    bar width=20pt,
    legend columns=3,
    legend entries={Exact methods, Greedy heuristics, Hybrid methods, Machine learning, Matheuristic, Neighborhood metaheur., Population metaheur., Tree-based heur.},
    legend to name=named,
    legend style={cells={anchor=west}, align=left, draw=none, column sep=2ex}
]
\addplot [fill=white, restrict x to domain=0:2] table[x index=0,y index=1] \datatable;
\addplot [fill=blue, restrict x to domain=0:2] table[x index=0,y index=2] \datatable;
\addplot [fill=red, restrict x to domain=0:2] table[x index=0,y index=3] \datatable;
\addplot [fill=cyan, restrict x to domain=0:2] table[x index=0,y index=4] \datatable;
\addplot [fill=lightgray, restrict x to domain=0:2] table[x index=0,y index=5] \datatable;
\addplot [fill=teal, restrict x to domain=0:2] table[x index=0,y index=6] \datatable;
\addplot [fill=orange, restrict x to domain=0:2] table[x index=0,y index=7] \datatable;
\addplot [fill=violet, restrict x to domain=0:2] table[x index=0,y index=8] \datatable;
\end{axis}
\end{tikzpicture}
&
\begin{tikzpicture}
\begin{axis}[
    height=6cm,
    width=.5\textwidth,
    ylabel=No. implementations,
    xtick=data,
    xticklabels={S,M,L},
    x tick style={draw opacity=0}, 
    xlabel={MP-CSPP},
    enlarge y limits=false,
    enlarge x limits=0.4,
    ymin=0,ymax=18,
    ybar stacked,
    bar width=20pt,
    legend columns=3,
    legend entries={Exact methods, Greedy heuristics, Hybrid methods, Machine learning, Matheuristic, Neighborhood metaheur., Population metaheur., Tree-based heur.},
    legend to name=named,
    legend style={cells={anchor=west}, align=left, draw=none, column sep=2ex}
]
\addplot [fill=white, restrict x to domain=3:5] table[x index=0,y index=1] \datatable;
\addplot [fill=blue, restrict x to domain=3:5] table[x index=0,y index=2] \datatable;
\addplot [fill=red, restrict x to domain=3:5] table[x index=0,y index=3] \datatable;
\addplot [fill=cyan, restrict x to domain=3:5] table[x index=0,y index=4] \datatable;
\addplot [fill=lightgray, restrict x to domain=3:5] table[x index=0,y index=5] \datatable;
\addplot [fill=teal, restrict x to domain=3:5] table[x index=0,y index=6] \datatable;
\addplot [fill=orange, restrict x to domain=3:5] table[x index=0,y index=7] \datatable;
\addplot [fill=violet, restrict x to domain=3:5] table[x index=0,y index=8] \datatable;
\end{axis}
\end{tikzpicture}
\\
\begin{tikzpicture}
\begin{axis}[
    height=6cm,
    width=.5\textwidth,
    ylabel=No. implementations,
    xtick=data,
    xticklabels={S,M,L},
    x tick style={draw opacity=0}, 
    xlabel={MPP},
    enlarge y limits=false,
    enlarge x limits=0.4,
    ymin=0,ymax=18,
    ybar stacked,
    bar width=20pt,
]
\addplot [fill=white, restrict x to domain=6:8] table[x index=0,y index=1] \datatable;
\addplot [fill=blue, restrict x to domain=6:8] table[x index=0,y index=2] \datatable;
\addplot [fill=red, restrict x to domain=6:8] table[x index=0,y index=3] \datatable;
\addplot [fill=cyan, restrict x to domain=6:8] table[x index=0,y index=4] \datatable;
\addplot [fill=lightgray, restrict x to domain=6:8] table[x index=0,y index=5] \datatable;
\addplot [fill=teal, restrict x to domain=6:8] table[x index=0,y index=6] \datatable;
\addplot [fill=orange, restrict x to domain=6:8] table[x index=0,y index=7] \datatable;
\addplot [fill=violet, restrict x to domain=6:8] table[x index=0,y index=8] \datatable;
\end{axis}
\end{tikzpicture}
&
\begin{tikzpicture}
\begin{axis}[
    height=6cm,
    width=.5\textwidth,
    ylabel=No. implementations,
    xtick=data,
    xticklabels={S,M,L},
    x tick style={draw opacity=0}, 
    xlabel={SPP},
    enlarge y limits=false,
    enlarge x limits=0.4,
    ymin=0,ymax=18,
    ybar stacked,
    bar width=20pt,
]
\addplot [fill=white, restrict x to domain=9:11] table[x index=0,y index=1] \datatable;
\addplot [fill=blue, restrict x to domain=9:11] table[x index=0,y index=2] \datatable;
\addplot [fill=red, restrict x to domain=9:11] table[x index=0,y index=3] \datatable;
\addplot [fill=cyan, restrict x to domain=9:11] table[x index=0,y index=4] \datatable;
\addplot [fill=lightgray, restrict x to domain=9:11] table[x index=0,y index=5] \datatable;
\addplot [fill=teal, restrict x to domain=9:11] table[x index=0,y index=6] \datatable;
\addplot [fill=orange, restrict x to domain=9:11] table[x index=0,y index=7] \datatable;
\addplot [fill=violet, restrict x to domain=9:11] table[x index=0,y index=8] \datatable;
\end{axis}
\end{tikzpicture}
\\
\multicolumn{2}{c}{
\ref{named}}

\end{tabular}



%% file: appendices/02_full_tables.tex
\begin{table}[h!]
\centering
\caption{\label{tab:scheme}Classification scheme}
\resizebox{\textwidth}{!}{%
\begin{tabular}{lll}
\hline
\textbf{Label}   & \textbf{Value}               & \textbf{Description}                                                                                    \\ \hline
\multicolumn{2}{l}{\textit{Port}}               & Port mode                                                                                               \\ \hline
\textit{}        & \textit{SP}                  & Single port                                                                                             \\ \hline
\textit{}        & \textit{MP}                  & Multi port                                                                                              \\ \hline
\multicolumn{2}{l}{\textit{Cargo. Length}}              & Length of containers                                                                     \\ \hline
\textit{}        & \textit{20'}             & 20' containers                                                                            \\ \hline
\textit{}        & \textit{40'}         & 40' containers                                                                    \\ \hline
\textit{}        & \textit{45'}               & 45' containers                                                                          \\ \hline
\multicolumn{2}{l}{\textit{Cargo. Height}}              & Height of containers                                                                        \\ \hline
\textit{}        & \textit{DC}             & Standard height containers                                                                         \\ \hline
\textit{}        & \textit{HC}         & High cube containers                                                                  \\ \hline
\multicolumn{2}{l}{\textit{Cargo. Weight}}              & Weight of containers                                                                     \\ \hline
\textit{}        & \textit{Mix}             & Mixed weight containers                                                                           \\ \hline
\textit{}        & \textit{WC or WCx}         & Containers grouped in weight classes, if x given, it is the number of weight classes considered                                                                  \\ \hline
\textit{}        & \textit{Uni}               & Uniform weight containers                                                                        \\ \hline
\multicolumn{2}{l}{\textit{Cargo. Specials}}              & Special containers                                                                    \\ \hline
\textit{}        & \textit{DG}             & Dangerous cargo                                                                          \\ \hline
\textit{}        & \textit{RF}         & Refrigerated containers                                                                  \\ \hline
\textit{}        & \textit{OT}                & Open-top containers                          \\ \hline
\multicolumn{2}{l}{\textit{Re}}         & Restow containers handling                                                                                       \\ \hline
\textit{}        & \textit{Min}            &
Involuntary container restows are allowed as well as minimized                        \\ \hline
\textit{}        & \textit{Vol}            & Voluntary and involuntary container restows are allowed as well as minimized                           \\ \hline
\textit{}        & \textit{NA}         & Container restows are not allowed                                                         \\ \hline
\textit{}        & \textit{None}                & Container restows are disregarded                                                                              \\ \hline
\multicolumn{2}{l}{\textit{HR}}                 & Hatch restows created by hatch cover lifts                                                                                       \\ \hline
\multicolumn{2}{l}{\textit{Hydrostatics. Stability}}       & Stability constraints                                                                                   \\ \hline
\textit{}        & \textit{GM}                  & Metacentric height                                                                                      \\ \hline
\textit{}        & \textit{Trim}                & Trim                                                                                                    \\ \hline
\textit{}        & \textit{List}                  & List                                                                                  \\ \hline
\textit{}        & \textit{LE}                  & Longitudinal equilibrium                                                                                 \\ \hline
\textit{}        & \textit{VE}                  & Vertical equilibrium                                                                                     \\ \hline
\textit{}        & \textit{TE}                  & Transversal equilibrium                                                                                     \\ \hline
\multicolumn{2}{l}{\textit{Hydrostatics. Stress forces}}       & Stress forces constraints                                                                                   \\ \hline
\textit{}        & \textit{SF}                  & Shear force                                                                                             \\ \hline
\textit{}        & \textit{BM}                  & Bending moment                                                                                          \\ \hline
\textit{}        & \textit{TM}                  & Torsion moment                                                                                          \\ \hline
\multicolumn{2}{l}{\textit{BW}}                 & Ballast water                                                                                           \\ \hline
\multicolumn{2}{l}{\textit{La}}                 & Lashing forces                                                                                                \\ \hline
\multicolumn{2}{l}{\textit{CO}}                 & Crane operations                                                                                             \\ \hline
\multicolumn{2}{l}{\textit{BS}}                 & Block stowage                                                                                             \\ \hline
\multicolumn{2}{l}{\textit{Obj}}          & Elements of the objective function                                                                      \\ \hline
                 & \textit{PS}           & Minimize port stay by minimizing overstowage and optimizing cranes work                                 \\ \hline
                 & \textit{VU}  & Maximizing vessel utilization, consolidation                                               \\ \hline
                 & \textit{H}        & Minimizing fuel consumption, improving hydrostatics                                                \\ \hline
\multicolumn{2}{l}{\textit{HD}}                 & Hierarchical decomposition                                                                              \\ \hline
\multicolumn{2}{l}{\textit{Sc}}              & Size of computational study                                                                             \\ \hline
\textit{}        & \textit{S}               & \begin{tabular}[c]{@{}l@{}}Small. Vessels with a capacity below 2,500 TEU, for slot planner blocks below 75 TEU \end{tabular}                     \\ \hline
\textit{}        & \textit{M}              & \begin{tabular}[c]{@{}l@{}}Medium. Vessels with a capacity between 2,500 and 15,000 TEU,\\~for slot planner blocks between 75 and 150 TEU\end{tabular}  \\ \hline
\textit{}        & \textit{L}               & \begin{tabular}[c]{@{}l@{}}Large. Vessels with a capacity above 15,000 TEU, for slot planner blocks above 150 TEU \end{tabular}            \\ \hline
\multicolumn{2}{l}{\textit{Solution   methods}} & Applied optimization approaches                                                                                      \\ \hline
\end{tabular}%
}
\end{table}

\begin{landscape}
\begin{table}[h!]
\centering
\caption{Full classification of single-port container stowage planning problems}
\resizebox{\columnwidth}{!}{%
\begin{tabular}{|l|l|l|l|l|l|l|l|l|l|l|l|l|l|l|l|l|l|} 
\hline
\multicolumn{1}{|c|}{{\textbf{Paper}}} & \multicolumn{1}{c|}{{\textbf{Port}}} & \multicolumn{4}{c|}{\textbf{Cargo}}                                                                               & \multicolumn{1}{c|}{{\textbf{Re}}} & \multicolumn{1}{c|}{{\textbf{HR}}} & \multicolumn{2}{c|}{\textbf{Hydrostatics}}                       & \multicolumn{1}{c|}{{\textbf{BW}}} & \multicolumn{1}{c|}{{\textbf{La}}} & \multicolumn{1}{c|}{{\textbf{CO}}} & \multicolumn{1}{c|}{{\textbf{BS}}} & \multicolumn{1}{c|}{{\textbf{Obj}}} & \multicolumn{1}{c|}{{\textbf{HD}}} & \multicolumn{1}{c|}{{\textbf{Sc}}} & \multicolumn{1}{c|}{{\textbf{Solution methods}}}  \\ 
\cline{3-6}\cline{9-10}
\multicolumn{1}{|c|}{}                                & \multicolumn{1}{c|}{}                               & \multicolumn{1}{c|}{\textbf{Length}} & \multicolumn{1}{c|}{\textbf{Height}} & \textbf{Weight} & \textbf{Specials} & \multicolumn{1}{c|}{}                             & \multicolumn{1}{c|}{}                             & \multicolumn{1}{c|}{\textbf{Stability}} & \textbf{Stress forces} & \multicolumn{1}{c|}{}                             & \multicolumn{1}{c|}{}                             & \multicolumn{1}{c|}{}                             & \multicolumn{1}{c|}{}                              & \multicolumn{1}{c|}{}                             & \multicolumn{1}{c|}{}                             & \multicolumn{1}{c|}{}  & \multicolumn{1}{c|}{}                                            \\ 
\hline
\cite{Cho1981DevelopmentPlanning}         & SP & 40'      & DC     & Mix & RF & NA & & GM, trim    &        & ~ & ~                         & ~           & ~              & PS, VU    & & S & Dynamic programming, neighborhood based  \\ 
\hline
\cite{Sciomachen2003TheProblem}   & SP & 10', 20', 30', 40' & DC, HC     & WC3 &  RF, OT      & NA  &     & LE, TE           &        &     &  &     & & PS     &  & S & Hybrid exact and greedy \\ \hline
\cite{Ambrosino2004StowingProblem}   & SP & 20', 40' &   DC   &  WC3 &       & NA  &     & LE, TE           &        &     &  &    & & PS     &  & S & Exact                       \\ \hline
\cite{Ambrosino2006AProblem}   & SP & 20', 40' & DC, HC & WC3 & ~  & NA & ~                         & LE, TE       &        & ~ & ~                         &  & ~ & PS        &  & S & Hybrid exact and neighborhood based        \\ 
\hline
\cite{Sciomachen2007AProductivity}   & SP & 20', 40' & DC     & WC3 &        & NA  &     & LE, TE           &        &     &  & \checkmark & & PS    &  & S & Hybrid exact and neighborhood based                               \\ \hline
\cite{Ambrosino2010AnProblem}   & SP & 20', 40' &   DC   &  WC3 &       & NA  &     & LE, TE           &        &     &  &    &  & PS   &   & S & Neighborhood based                       \\ \hline
\cite{Delgado2012AStowage}     & SP & 20', 40' & DC, HC & Mix & RF & NA  & &              &        & ~ & ~         & ~                & ~                         & VU        & \checkmark & L & Hybrid exact and greedy                      \\ 
\hline
\cite{Hu2012CombinatorialTerminal}          & SP & 40'      & DC     & WC  & ~  & Min & ~                         & GM, trim &        & ~ & ~                         &  & ~ & PS, H     &  & S & Population based               \\ 
\hline
\cite{Cruz-Reyes2015}    & SP & 20', 40' & DC     & WC3 &        & NA  &     & LE, TE           &        &     &  &  &   & PS   &   & S & Hybrid exact and greedy                       \\ \hline
\cite{Shen2017AProblem}      & SP & 40'           & DC     & Mix &        & Min &     &  &  &  &     &   &  & PS  &   & S & Machine learning                   \\ \hline
\cite{Zhao2018ContainerSearch}         & SP & 20', 40' & DC, HC & Mix & ~  & Min & ~                         &              &        & ~ & ~                         & ~            & ~             & PS        & ~                         & S & Tree based                         \\ 
\hline
\cite{Li2020OptimizingShipping}          & SP & 20'      & DC     & Mix & ~  & Min, NA & ~                         & LE, TE   &        & ~ & ~                         & ~             & ~            & PS        &  & S & Neighborhood based, exact                       \\ 
\hline
\cite{Zhu2020IntegerProblem}          & SP & 20', 40' & DC     & Mix & RF & Min & \checkmark & GM, trim     & SF     &   &                           &               & ~            & PS        &                           & S & Exact                              \\ 
\hline
\cite{Larsen2021AProblem}      & SP & 20', 40' & DC, HC & WC6 & RF & Min & \checkmark & GM, trim, list  & SF, BM & ~ & ~                         & \checkmark & \checkmark & PS, VU, H & ~                         & L & Neighborhood based                       \\ 
\hline
\cite{ElYaagoubi2022Multi-objectiveSystem} & SP & 20', 40' & DC     & Mix & ~  & Min & ~                         & GM           &        & ~ & \checkmark & ~                    & ~      & PS, H     & ~                         & S & Population based, exact                \\
\hline
\end{tabular}
}
\end{table}

\begin{table}[h!]
\centering
\caption{Full classification of multi-port container stowage planning problems}
\resizebox{\columnwidth}{!}{%
\begin{tabular}{|l|l|l|l|l|l|l|l|l|l|l|l|l|l|l|l|l|l|} 
\hline
\multicolumn{1}{|c|}{{\textbf{Paper}}} & \multicolumn{1}{c|}{{\textbf{Port}}} & \multicolumn{4}{c|}{\textbf{Cargo}}                                                                               & \multicolumn{1}{c|}{{\textbf{Re}}} & \multicolumn{1}{c|}{{\textbf{HR}}} & \multicolumn{2}{c|}{\textbf{Hydrostatics}}                       & \multicolumn{1}{c|}{{\textbf{BW}}} & \multicolumn{1}{c|}{{\textbf{La}}} & \multicolumn{1}{c|}{{\textbf{CO}}} &
\multicolumn{1}{c|}{{\textbf{BS}}} &
\multicolumn{1}{c|}{{\textbf{Obj}}} & \multicolumn{1}{c|}{{\textbf{HD}}} & \multicolumn{1}{c|}{{\textbf{Sc}}} & \multicolumn{1}{c|}{{\textbf{Solution methods}}}  \\ 
\cline{3-6}\cline{9-10}
\multicolumn{1}{|c|}{}                                & \multicolumn{1}{c|}{}                               & \multicolumn{1}{c|}{\textbf{Length}} & \multicolumn{1}{c|}{\textbf{Height}} & \textbf{Weight} & \textbf{Specials} & \multicolumn{1}{c|}{}                             & \multicolumn{1}{c|}{}                             & \multicolumn{1}{c|}{\textbf{Stability}} & \textbf{Stress forces} & \multicolumn{1}{c|}{}                             & \multicolumn{1}{c|}{}                             & \multicolumn{1}{c|}{}                             & \multicolumn{1}{c|}{}                              & \multicolumn{1}{c|}{}                             & \multicolumn{1}{c|}{}                             & \multicolumn{1}{c|}{}  & \multicolumn{1}{c|}{}                                            \\ 
\hline
\cite{Shields1984}                                           & MP                                                  & 20', 40'                             & DC                                   & WC              & RF                & Min                                               & \checkmark                         & GM, trim, list                    & SF                     & \checkmark                         & \checkmark                         &                                                & ~    & PS, VU, H                                          &                                                   & M                                                 & Greedy                                                           \\ 
\hline
\cite{Botter1992StowageSolution}                                             & MP                                                  & 20', 40'                             & DC                                   & Mix             &                   & Vol                                               & \checkmark                         & GM, trim, list                       & SF, BM                 & \checkmark                         &                                                   & ~                        & ~ & PS                                                 & \checkmark                         & S                                                 & Hybrid exact and tree based                                                 \\ 
\hline
\cite{AvrielMordecaiPenn1993}                                           & MP                                                  & 40'                                  & DC                                   & Uni             &                   & Vol                                               &                                                   &                                         &                        &                                                   &                                                   &                                            & ~       & PS                                                 &                                                   & S                                                 & Greedy                                                    \\ 
\hline
\cite{Avriel1998StowageShifts}                                             & MP                                                  & 40'                                  & DC                                   & Uni             &                   & Vol                                               &                                                   &                                         &                        &                                                   &                                                   &                                           & ~        & PS                                                 &                                                   & S                                                 & Greedy                                                           \\ 
\hline
\cite{Wilson2000ContainerSolutions}                                             & MP                                                  & 20', 40'                             & DC                                   & WC              & DG, RF            & Min                                               & \checkmark                         &                                         &                        &                                                   &                                                   & \checkmark            & \checkmark               & PS, VU                                             & \checkmark                         & S                                                 & Hybrid exact and neighborhood based                                              \\ 
\hline
\cite{Dubrovsky2002AProblem}                                        & MP                                                  & 40'                                  & DC                                   & Uni             &                   & Vol                                               &                                                   & List                                      &                        &                                                   &                                                   &                                           & ~        & PS                                                 &                                                   & S                                                 & Population based                                             \\ 
\hline
\cite{Kang2002StowageTransportation}                                                & MP                                                  & 40'                                  & DC                                   & WC              &                   & Min                                               & \checkmark                         & GM, trim, list                      &                        &                                                   &                                                   &       & ~                                            & PS                                                 & \checkmark                         & M                                                 & Hybrid greedy and tree based                                                \\ 
\hline
\cite{Kaisar2006ATRANSPORTATION}                                               & MP                                                  & 20', 40'                             & DC                                   & WC              & DG, RF            & Min                                               &                          & VE, LE, TE                               &                        &                                                   &                                                   &                                           & ~        & PS                                                 & \checkmark                         & M                                                 & Hybrid exact and neighborhood based                                              \\ 
\hline
\cite{Hamedi2011CONTAINERSHIPOPERATIONS}                                             & MP                                                  & 20', 40'                             & DC                                   & Mix             & DG, RF            & Min                                               &                                                   & Trim, list                             &                        &                                                   &                                                   & \checkmark        & ~                 & PS                                                 &                                                   & M                                                 & Population based                                                     \\ 
\hline
\cite{Liu2011RandomizedPlans}                                               & MP                                                  & 20', 40'                             & DC, HC                               & Class             & DG, RF                & Min                                               & \checkmark                         & GM, trim list                           &                        &                                                   &                                                   & \checkmark            & \checkmark             & PS, VU, H                                          &                         & M                                                 & Hybrid greedy and neighborhood based                                             \\ 
\hline
\cite{Pacino2011FastVessels}                                            & MP                                                  & 20', 40'                             & DC                                   & Mix             & RF                & Min                                               & \checkmark                         & GM, trim                   & SF                     &                                                   &                                                   & \checkmark             & ~            & PS                                                 & \checkmark                         & M                                                 & Hybrid exact and neighborhood based                      \\ 
\hline
\cite{Azevedo2014SolvingMeta-heuristics}                                            & MP                                                  & 40'                                  & DC                                   & Uni             &                   & Vol                                               &                                                   & LE, TE                                     &                        &                                                   &                                                   &                                           & ~        & PS, H                                              &                                                   & S                                                 & Neighborhood based, tree based, population based                           \\ 
\hline
\cite{Ding2015StowageShifts}                                             & MP                                                  & 40'                                  & DC                                   & Uni             &                   & Vol                                               &                                                   &                                         &                        &                                                   &                                                   &                                           & ~        & PS                                                 &                                                   & M                                                 & Greedy                                                           \\ 
\hline
\cite{Hu2017Multi-objectiveRoute}                                                 & MP                                                  & 20', 40'                             & DC                                   & Mix             &                   & Min                                               &                                                   & GM, trim, list                             &                        &                                                   &                                                   &                                           & ~        & PS, H                                              &                         & S                                                 & Population based                                            \\ 
\hline
\cite{Azevedo2018SolvingAlgorithm}                                           & MP                                                  & 40'                                  & DC                                   & Uni             &                   & Vol                                               &                                                   & LE, TE                                      &                        &                                                   &                                                   & \checkmark                        & ~ & PS                                                 &                                                   & S                                                 & Population based                                                     \\ 
\hline
\cite{Li2018Multi-PortUncertainties}                                                & MP                                                  & 40'                                  & DC                                   & WC              &                   & NA                                                &                                                   & VE, LE, TE                              &                        &                                                   &                                                   &                                           & ~        & VU                                             &                                                   & S                                                 & Neighborhood based, exact                                     \\ 
\hline
\cite{Pacino2018CranePlanning}                                             & MP                                                  & 40'                                  & DC                                   & WC              &                   & NA                                                &                                                   &                                         &                        &                                                   &                                                   & \checkmark        & \checkmark                 & PS                                                 &                         & L                                                 & Neighborhood based                                                     \\ 
\hline
\cite{Roberti2018APlans}                                           & MP                                                  & 40'                                  & DC                                   & Uni             &                   & Vol, NA                                               &                                                   &                                         &                        &                                                   &                                                   &                                           & ~        & VU                                                 &                                                   & M                                                 & Exact                                                            \\ 
\hline
\cite{Parreno-Torres2019SolutionProblem}                                            & MP                                                  & 40'                                  & DC                                   & Uni             &                   & Vol                                               &                                                   &                                    &                        &                                                   &                                                   &                                           & ~        & PS                                                 &                                                   & L                                                 & Exact, neighborhood based                                                    \\ 
\hline
\cite{Parreno-Torres2020ImprovingProblems}                                             & MP                                                  & 20', 40'                             & DC                                   & WC3             &                   & Vol                                               &                                                   & Trim, list                                  &                        &                                                   &                                                   &                                          & ~         & PS                                                 &                                                   & L                                                 & Exact, matheuristic, neighborhood based                     \\ 
\hline
\cite{Parreno-Torres2021SolvingAlgorithm}                                           & MP                                                  & 20', 40'                             & DC                                   & WC3             &                   & Vol                                               &                                                   & GM, trim                                  &                        &                                                   &                                                   &                                         & ~          & PS                                                 &                                                   & M                                                 & Exact, matheuristic                        \\ 
\hline
\cite{Chang2022SolvingMode}                                              & MP                                                  & 20', 40'                             & DC                                   & Mix             &  (DG), RF                 & Min                                               &                                                   & Trim, list                            &                        &                                                   &                                                   & \checkmark        & ~                 & PS                                                 &                                                   & S                                                 & Population based                                                     \\
\hline
\end{tabular}
}
\end{table}
\end{landscape}

\begin{landscape}
    
\begin{table}[h!]
\centering
\caption{Full classification of master bay planning problems}
\resizebox{\columnwidth}{!}{%
\begin{tabular}{|l|l|l|l|l|l|l|l|l|l|l|l|l|l|l|l|} 
\hline
\multicolumn{1}{|c|}{{\textbf{Paper}}} & \multicolumn{1}{c|}{{\textbf{Port}}} & \multicolumn{4}{c|}{\textbf{Cargo}}                                                                               & \multicolumn{1}{c|}{{\textbf{Re}}} & \multicolumn{1}{c|}{{\textbf{HR}}} & \multicolumn{2}{c|}{\textbf{Hydrostatics}}                       & \multicolumn{1}{c|}{{\textbf{BW}}} & \multicolumn{1}{c|}{{\textbf{La}}} & \multicolumn{1}{c|}{{\textbf{CO}}} & \multicolumn{1}{c|}{{\textbf{Obj}}}  & \multicolumn{1}{c|}{{\textbf{Sc}}} & \multicolumn{1}{c|}{{\textbf{Solution methods}}}  \\ 
\cline{3-6}\cline{9-10}
\multicolumn{1}{|c|}{}                                & \multicolumn{1}{c|}{}                               & \multicolumn{1}{c|}{\textbf{Length}} & \multicolumn{1}{c|}{\textbf{Height}} & \textbf{Weight} & \textbf{Specials} & \multicolumn{1}{c|}{}                             & \multicolumn{1}{c|}{}                             & \multicolumn{1}{c|}{\textbf{Stability}} & \textbf{Stress forces} & \multicolumn{1}{c|}{}                             & \multicolumn{1}{c|}{}                                       & \multicolumn{1}{c|}{}                              & \multicolumn{1}{c|}{}                             & \multicolumn{1}{c|}{}                             & \multicolumn{1}{c|}{}         \\        
\hline

\cite{Pacino2012AnTanks}        & MP & 20', 40' & DC     & WC2 & RF     & None & \checkmark & GM, trim & SF, BM & \checkmark &  &     & H      & L & Exact                               \\ 
\hline
\cite{Pacino2013AnPlanning}          & MP & 20', 40' & DC     & WC2 & RF     & (Min)  & \checkmark & GM, trim           & SF, BM &     &  & \checkmark & PS     & M & Neighborhood metaheuristic          \\ 
\hline
\cite{Ambrosino2015APlanning}      & MP & 20', 40' & DC     & WC3 & RF, OT & Min  & \checkmark & LE, TE           &        &     &  & \checkmark & PS     & L & Matheuristic                 \\ 
\hline
\cite{Ambrosino2015ComputationalProblems} & MP & 20', 40' & DC     & WC3 & RF, OT & Min  & \checkmark & LE, TE           &        &     &  & \checkmark & PS     & L & Exact                               \\ 
\hline
\cite{Ambrosino2015ExperimentalProblem}   & MP & 20', 40' & DC     & WC3 &        & Min  & \checkmark & LE, TE           &        &     &  & \checkmark & PS, VU & M & Exact, matheuristic                 \\ 
\hline
\cite{Ambrosino2018ShippingApproach}      & MP & 20', 40' & DC     & WC3 & RF, OT & Min  & \checkmark & GM, trim           &  SF      &     &  & \checkmark & PS     & L & Exact, matheuristic                 \\ 
\hline
\cite{Kebedow2018IncludingProblem}        & MP & 20', 40' & DC     & WC2 & DG, RF & None & \checkmark & GM, trim, list      & SF     &     &  & \checkmark & PS, VU & L & Exact                               \\ 
\hline
\cite{Bilican2020AParameters}       & MP & 20', 40' & DC     & WC6 &        & Min  & \checkmark & Trim             & SF, BM & \checkmark &  &     & PS, H  & L & Hybrid exact and neighborhood based  \\ 
\hline
\cite{Chao2021}            & MP & 20', 40' & DC, HC & WC  &        & Min  & \checkmark &                  &        &     &  &     & PS, VU & M & Exact                               \\
\hline

\end{tabular}
}
\end{table}

\begin{table}[h!]
\centering
\caption{Full classification of slot planning problems}
\resizebox{\columnwidth}{!}{%
\begin{tabular}{|l|l|l|l|l|l|l|l|l|l|l|l|l|l|l|l|} 
\hline
\multicolumn{1}{|c|}{{\textbf{Paper}}} & \multicolumn{1}{c|}{{\textbf{Port}}} & \multicolumn{4}{c|}{\textbf{Cargo}}                                                                               & \multicolumn{1}{c|}{{\textbf{Re}}} & \multicolumn{1}{c|}{{\textbf{HR}}} & \multicolumn{2}{c|}{\textbf{Hydrostatics}}                       & \multicolumn{1}{c|}{{\textbf{BW}}} & \multicolumn{1}{c|}{{\textbf{La}}} & \multicolumn{1}{c|}{{\textbf{CO}}} & \multicolumn{1}{c|}{{\textbf{Obj}}}  & \multicolumn{1}{c|}{{\textbf{Sc}}} & \multicolumn{1}{c|}{{\textbf{Solution methods}}}  \\ 
\cline{3-6}\cline{9-10}
\multicolumn{1}{|c|}{}                                & \multicolumn{1}{c|}{}                               & \multicolumn{1}{c|}{\textbf{Length}} & \multicolumn{1}{c|}{\textbf{Height}} & \textbf{Weight} & \textbf{Specials} & \multicolumn{1}{c|}{}                             & \multicolumn{1}{c|}{}                             & \multicolumn{1}{c|}{\textbf{Stability}} & \textbf{Stress forces} & \multicolumn{1}{c|}{}                             & \multicolumn{1}{c|}{}                             & \multicolumn{1}{c|}{}                             & \multicolumn{1}{c|}{}                              & \multicolumn{1}{c|}{}                              & \multicolumn{1}{c|}{}         \\        
\hline
\cite{Pacino2010ABays}    & SP & 20', 40'      & DC, HC & Mix & RF     & Min &     &  &  &  &     &     & PS, VU & L & Neighborhood based        \\ 
\hline
\cite{Delgado2012ABays}    & SP & 20', 40'      & DC, HC & Mix & RF     & Min &  &  &  &  &     &     & PS, VU & L & Exact      \\ 
\hline
\cite{Parreno2016AProblem}   & SP & 20', 40'      & DC, HC & Mix & DG, RF & Min &     &  &  &  &     &     & PS, VU & M & Neighborhood based        \\ 
\hline
\cite{Yifan2016Group-BayShip}     & SP & 20', 40'      & DC, HC & Mix &        & Min &     &  &  &  &     &     & PS  & S & Hybrid greedy and population based   \\ 
\hline
\cite{Jin2019AnBay}        & SP & 40'           & DC     & Mix &        & Min &     &  &  &  &     &     & PS     & S & Exact                              \\ 
\hline
\cite{Kebedow2019IncludingProblem}    & SP & 20', 40'      & DC     & Mix & DG, RF & Min &     &  &  &  &     &     & PS, VU & M & Exact                              \\ 
\hline
\cite{Korach2020MatheuristicsBays}     & SP & 20', 40'      & DC, HC & Mix & RF     & Min &     &  &  &  &     &     & PS, VU & L & Matheuristic                       \\ 
\hline
\cite{Rashed2021AVessels}     & SP & 20', 40'      & DC, HC & Mix & RF     & Min &     &  &  &  &     &     & PS, VU & L & Neighborhood metaheuristic                              \\ 
\hline
\end{tabular}
}
\end{table}

\end{landscape}

%% file: appendices/01_appendixA.tex
\begin{table}[!h]
    \centering
    \small
    \input{tables/multi_port_full_problem_results.tex}
    \caption{Reported results from \cite{Pacino2011FastVessels} and \cite{Kang2002StowageTransportation}}
    \label{tab:mutli_port_full_problem_results}
\end{table}

%% file: tables/multi_port_full_problem_results.tex
\begin{tabular}{rrr}
         TEUs	&	Num. ports	&	Solve time (sec.)	\\\hline\\
        \multicolumn{3}{c}{Results from \cite{Pacino2011FastVessels}}\\\hline
        4755	&	4	&	10	\\
        9618	&	4	&	30	\\
        9984	&	4	&	21	\\
        2584	&	5	&	3	\\
        4755	&	5	&	21	\\
        4456	&	6	&	16	\\
        6545	&	6	&	5	\\
        8490	&	6	&	31	\\
        6717	&	7	&	23	\\
        7490	&	8	&	10	\\
        4478	&	9	&	5	\\
        5047	&	9	&	263	\\
        5052	&	9	&	214	\\
        4478	&	10	&	332	\\
        7344	&	10	&	252	\\
        9160	&	10	&	2079	\\
        9118	&	11	&	3711	\\
        9118	&	11	&	{\em timeout}	\\
        5044	&	14	&	{\em timeout}	\\
        9160	&	16	&	69	\\\hline\\
        \multicolumn{3}{c}{Results from \cite{Kang2002StowageTransportation}}\\\hline
        2500	&	4	&	32	\\
        3000	&	4	&	21	\\
        4000	&	4	&	32	\\
        2500	&	6	&	119	\\
        3000	&	6	&	121	\\
        4000	&	6	&	123	\\
        2500	&	8	&	366	\\
        3000	&	8	&	381	\\
        4000	&	8	&	399	\\\hline
    \end{tabular}

%% file: chapters/04_pacino11model.tex

The master planning problem aims to allocate cargo to subsections of bays. Those are often called locations \citep{Pacino2011FastVessels} or blocks \citep{Chou2021ApplyingStudy}. A block can either be a logical grouping of containers or be defined by the position of the hatch covers. The assignment of containers to blocks must ensure that the vessel is seaworthy while minimizing the handling time of the vessel. This appendix briefly presents the master planning formulations used for the computational comparison of Section~\ref{sec:master_planning}. We refer the reader to the original publications for a more in-depth description. Tables~\ref{tab:generalNotationSets} and \ref{tab:generalNotationPar} show the sets and parameters shared between the \cite{Pacino2011FastVessels} and \cite{Chou2021ApplyingStudy} formulations.

\begin{table}[h]
    \centering
    \begin{tabularx}{\textwidth}{lX}
        \multicolumn{2}{l}{{\bf Sets}}\\
        $B$& The set of bays  \\
        $BL$ & The set of blocks \\
        $BL_b$ & The set of blocks in bay $b\in B$\\
        $BL^O$ & The set of blocks over deck\\
        $BL^U_l$ & The set of blocks below deck under block $l\in L^O$\\
        $Bin$ & The set of adjacent bays $(b_1,b_2) \in B\times B$\\
        $T$ & The set of container types\\
        $T^{\{20,40,R\}}$ & The set of 20-, 40-foot and Reefer container types\\ 
        $P$ & The set of ports\\
        $TR$ & The set of transports (port pairs $(o,d) \in P\times P,o<p$)\\
        $TR^{ON}_p$ & The set of transports $(o,d) \in TR$ where $o<p$ and $d>p$\\
        $TR^{A}_p$ & The set of active transports at port $p\in P$ where $(o,d) \in TR$ where $o=p$ or $d=p$\\
        $TR^{OV}_p$ & The set of overstowing transports at port $p\in P$ where $(o,d) \in TR$ where $o<p$ and $d>p$\\
    \end{tabularx}
    \caption{Common Sets}
    \label{tab:generalNotationSets}
\end{table}
\begin{table}[!h]
    \centering
    \begin{tabularx}{\textwidth}{lX}
        \multicolumn{2}{l}{{\bf Parameters}}\\
        $K^{\{20,40,R\}}_l$ & The 20-,40-foot, and Reefer capacity of block $l\in L$\\
        $TEU^\tau$ & The Twenty-Foot Equivalent units of container type $\tau \in T$\\
        $W_\tau$ & The weight of container type $\tau \in T$\\
        $L_\tau$ & The length of container type $\tau \in T$\\
        $LD^{\tau}_t$ & The number of container of type $\tau \in T$ to be loaded for transport $t\in TR$\\
        $R^\tau_{p,l}$ & The number of containers of type $\tau \in T$, already on board the vessel with destination $p\in P$.\\
        $W^K_b$ & The lightship weight at bay $b\in B$\\
        $W^{\max}_l$ & The maximum weight limit of block $l\in L$\\
        $D_p$ & The total displacement of the vessel leaving port $p\in P$\\
        $CG^{\alpha}_l$ & The centre of gravity components $\alpha \in {L=LCG, V=VCG, T=TCG}$ of block $l \in BL$\\
        $\overline{CG}^{\alpha}_b$ & The centre of gravity components $\alpha \in {L=LCG, V=VCG, T=TCG}$ of bay $b \in B$\\
        $LCG^{\{Min,Max\}}$ & The limits for the vessel's longitudinal centre of gravity\\
        $VCG^{\{ax\}}$ & The maximum vertical centre of gravity of the the vessel\\
        $TCG^{\{Min,Max\}}$ & The limits for the vessel's transversal centre of gravity\\
        $Shear^{\{Min,Max\}}_b~\quad~\quad$ & The shear limits at bay $b\in B$\\
    \end{tabularx}
    \caption{Common Paremeters}
    \label{tab:generalNotationPar}
\end{table}
\clearpage
\subsection{\cite{Pacino2011FastVessels} formulation}
The formulation is based on a main set of decision variables indicating the number of containers of a specific type to be stowed on a block during a transport leg. Other indicator variables are used for the calculation of the objective value. A description of the variables follows.
\\\\
\begin{tabularx}{\linewidth}{lX}
     $x^\tau_{tl}~\in~\mathbb{Z}^+~\quad $ & The number of containers of type $\tau \in T$ stowed in block $l\in L$ during transport $t \in TR$.\\
     $\delta_{pl}~\in~\mathbb{B}~\quad $ & Stowage indicator equal to 1 if any load or discharge containers are present in block $l \in L$ at port $p\in P$.\\
     $y^O_{pl}~\in~\mathbb{R}^+$ & The number of hatch-overstow containers in block $l\in L$ at port $p\in P$\\
     $y^T_p~\in~\mathbb{R}^+$ & The long-crane (or makespan lower bound) at port $p\in P$.\\
\end{tabularx}
\\\\

Following is the mathematical formulation and its description.
\begin{align}
    & \text{\bf min} \sum_{p\in P} \left(\mathcal{C}^Ty^T_p + \sum_{l\in L} \mathcal{C}^Oy^O_pl \right) \label{pacino:obj}\\
    & \text{\bf s.t.}\nonumber \\
    & \sum_{l\in L} x^\tau_{tl} = LD^\tau_t \quad \forall \tau \in T, t \in TR \label{pacino:loadall}\\
    & x^\tau_{(1,p)l} = R^\tau_p \quad \forall \tau \in T, p\in P \label{pacino:release}\\
    & \sum_{t \in TR^{ON}_p}\sum_{\tau \in T} TEU^\tau x^\tau_{tl} \leq K^{20}_l \quad \forall p\in P, l\in L\label{pacino:capTEU}\\
    & \sum_{t \in TR^{ON}_p}\sum_{\tau \in T^\alpha} x^\tau_{tl} \leq K^{\alpha}_l \quad \forall p\in P, l\in L, \alpha\in \{20,40,R\}\label{pacino:cap}\\
    &\sum_{t \in TR^{ON}_p}W_\tau x^\tau_{tl} \leq W^{\max}_l \quad \forall p\in P, l\in L \label{pacino:weight}\\
    &\sum_{b\in B}W^K_b\overline{CG}^L_b +\!\! \sum_{t \in TR^{ON}_p}\sum_{l\in L}CG^L_lW_\tau x^\tau_{tl} \geq LCG^{Min}D_p \quad \forall p\in P\label{pacino:minLCG}\\
    &\sum_{b\in B}W^K_b\overline{CG}^L_b +\!\! \sum_{t \in TR^{ON}_p}\sum_{l\in L}CG^L_lW_\tau x^\tau_{tl} \leq LCG^{Max}D_p \quad \forall p\in P\label{pacino:maxLCG}\\
    &\sum_{b\in B}W^K_b\overline{CG}^T_b +\!\! \sum_{t \in TR^{ON}_p}\sum_{l\in L}CG^T_lW_\tau x^\tau_{tl} \geq TCG^{Min}D_p \quad \forall p\in P\label{pacino:minTCG}\\
    &\sum_{b\in B}W^K_b\overline{CG}^T_b +\!\! \sum_{t \in TR^{ON}_p}\sum_{l\in L}CG^T_lW_\tau x^\tau_{tl} \leq TCG^{Max}D_p \quad \forall p\in P\label{pacino:maxTCG}\\
    &\sum_{b\in B}W^K_b\overline{CG}^V_b +\!\! \sum_{t \in TR^{ON}_p}\sum_{l\in L}CG^V_lW_\tau x^\tau_{tl} \leq VCG^{Max}D_p \quad \forall p\in P\label{pacino:maxVCG}\\
    & \sum_{b'=1}^b\left( W^K_{b'} +\!\! \sum_{l\in BL_{b'}}\sum_{\tau\in T}\sum_{t \in TR^{ON}_p}\!\!W_\tau x^\tau_{tl} \right)\geq Shear^{Min}_b \quad \forall p \in P, b \in B\label{pacino:shearMin}\\
    & \sum_{b'=1}^b\left( W^K_{b'} +\!\! \sum_{l\in BL_{b'}}\sum_{\tau\in T}\sum_{t \in TR^{ON}_p}\!\!W_\tau x^\tau_{tl} \right)\leq Shear^{Max}_b \quad \forall p \in P, b \in B\label{pacino:shearMax}\\
    &\sum_{\tau \in T}\sum_{t\in TR^A_p}\sum_{l' \in BL^U_l}x^\tau_{tl'} M\delta_{pl} \quad \forall p\in P, l \in BL^O \label{pacino:delta}\\
    & \sum_{\tau \in T}\sum_{t\in TR^{OV}_p}x^\tau_{tl} - M(1-\delta_{pl})\leq y^O_p \quad \forall p\in P, l \in BL^O\label{pacino:hatch}\\
    &\sum_{b \in N}\sum_{\tau \in T}\sum_{t \in TR^{ON}_p}\sum_{l\in B} x^\tau_{tl} \leq y^T_p\quad \forall p\in P, N \in Bin \label{pacino:makespan}
\end{align}

The formulation minimizes the makespan and the hatch overstowage at each port~\eqref{pacino:obj}. Constraint~\eqref{pacino:loadall} ensures that all cargo must be stowed on the vessel, while constraint~\eqref{pacino:release} enforces that cargo already on board does not change position. The total block capacity, the type-specific block capacity, and the block weight capacity are constrained by~\eqref{pacino:capTEU},~\eqref{pacino:cap},~\eqref{pacino:weight}, respectively. To ensure vessel stability, the longitudinal, transversal, and vertical centers of gravity are constrained within the given limits by constraints~\eqref{pacino:minLCG}~-~\eqref{pacino:maxVCG}.Constraints~\eqref{pacino:shearMin} and~\eqref{pacino:shearMax} impose minimum and maximum levels for the shear forces that act on the vessel.

Using an indicator variable, constraint~\eqref{pacino:delta} identifies, for each port, blocks below deck that require container moves. Should those container moves be blocked by containers on-deck, they will be captured in constraint~\eqref{pacino:hatch} as overstowing. Finally, the bay pair with the maximum number of movements at each port is identified with constraint~\eqref{pacino:makespan}, representing the makespan. 

\subsection{\citep{Chao2021} formulation}
The model is based on a network-flow formulation. For a detailed description of the network structure, we refer the reader to the original publication ~\citep{Chao2021}. Following is an extension of the common sets and parameters needed for this formulation.
\\\\
\begin{tabularx}{\linewidth}{lX}
    \multicolumn{2}{l}{{\bf Sets}}\\
    $N$ & The set of all nodes\\
    $P$ & The set of ports\\
    $\mathcal{A}$ & The set of all arcs\\
    $\mathcal{A}^-_i, \mathcal{A}^+_i$ & The set of outgoing and incoming arcs of node $i\in N$\\
    $OD$ & The set of transports (origin/destination port pairs) $(o,d) \in P \times P$ \\
    $OD^{ON}$ & The set of transports $t=(o,d) \in OD$ for which $o\leq p$ and $d>p$\\
    $OD^{A}$ & The set of transports $t=(o,d) \in OD$ for which $o=p$ or $d=p$\\
    $A^{TR}_{t}$ & The set of all arcs belonging to the transport $t \in OD$\\
    $A^{\tau}_t$ & The set of arcs connecting the transport nodes and the container type nodes\\
    $A^{L}_{tl}$ & The set of arcs connecting the container type nodes and the block nodes\\
    $K_t$ & The number of containers in transport $t\in OD$\\
    $E_{ij}$ & The number of containers flowing through the type nodes\\
    \multicolumn{2}{l}{{\bf Parameters}}\\
    $S$ & The source node\\
    $T$ & The sink node\\
    $q_i$ & Is equal to $q$ for $i=S$, $-q$ for $i=T$ and zero for all other nodes, where $q$ is the total number of containers to be stowed\\
    $\tau^L(i)$ & The length of container type node $i\in T$\\
    $\tau^W(i)$ & The weight of container type node $i\in T$\\
    $TEU_i$ & The TEU value of container type node $i\in T$\\
\end{tabularx}
\\\\
As the following is a network-flow formulation, the decision variables $x_{ij} \in \mathbb{R}^+$ represent the flow (the number) of containers from the source to the sink node, through the arcs $(i,j)\in A$.

\begin{align}
    & \text{\bf min} \sum_{p\in P}\mathcal{C}^Ty^T_p \label{chao:obj}\\
    & \text{\bf s.t.}\nonumber \\
    & \sum_{(i,j) \in \mathcal{A}^-_i}\!\! x_{ij} -\!\! \sum_{(i,j)\in \mathcal{A}^+_i}\!\! x_{ij} = q_i \quad \forall i \in N \label{chao:flow}\\
    & x_{ij} = K_t \quad \forall t \in OD, (i,j) \in A^{TR}_{t} \label{chao:od_pairs}\\
    & x_{ij} = E_{ij} \quad \forall t \in OD, (i,j) \in A^{\tau}_{t}\label{chao:types}\\
    & \sum_{t\in OD^{ON_p}}\sum_{(i,j)\in A^L_{tl}} TEU_i x_{ij} \leq C^{20}_l \quad \forall p\in P, l\in L\label{chao:cap_teu}\\
    & \sum_{t\in OD^{ON_p}}\sum_{(i,j)\in A^L_{tl},\tau(i)=\alpha}x_{ij} \leq C^\alpha_l \quad \forall p\in P, l\in L, \alpha \in \{20,40,R\}\label{chao:cap}\\
    & \sum_{t\in OD^{OD}_p}\sum_{(i,j)\in A^L_{tl}}\tau^W(i)x_{ij}\leq W^{Max}_l \quad \forall p\in P, l\in L \label{chao:cap_weight}\\
    &x_{ij} \leq M y_{ij} \quad \forall p\in P, t \in OD^A_p, (i,j)\in A^{LU}_{tl}\label{chao:delta}\\
    &\sum_{t\in OD^O_p}\sum_{(i,j)\in A^{LO}_{tl}} x_{ij} \leq M(1- y_{ij}) \quad \forall p\in P, t \in OD^A_p, (i,j)\in A^{LU}_{tl}\label{chao:hatch}\\
    &\sum_{b\in N}\sum_{l\in BL_b}\sum_{A^L_{tl}} x_{ij}\leq y^T_p \quad \forall p\in P, N \in Bin \label{chao:makespan}
\end{align}

The model's objective is the minimization of the makespan at each port~\eqref{chao:obj}. This objective differs from the original publication, but this change was necessary for the sake of the comparison. Constraint~\eqref{chao:flow} is the classic flow-conservation constraint between the source and the sink node. To ensure that cargo is correctly routed through the network, it is necessary to constraint containers within their respective origin-destination arcs~\eqref{chao:od_pairs}, and the correct container type arcs~\eqref{chao:types}. Constraints~\eqref{chao:cap_teu},~\eqref{chao:cap}, and~\eqref{chao:cap_weight} are the revised capacity and weight constraints. The absence of hatch=overstowage is ensured by constraints~\eqref{chao:delta}, which indicate the presents of container moves in below deck locations, and constraint~\eqref{chao:hatch} which imposes that no hatch-overstowage is allowed. Finally, constraint~\eqref{chao:makespan} calculates the makespan at each port.

%% file: main.bbl
\begin{thebibliography}{88}
\expandafter\ifx\csname natexlab\endcsname\relax\def\natexlab#1{#1}\fi
\providecommand{\url}[1]{\texttt{#1}}
\providecommand{\href}[2]{#2}
\providecommand{\path}[1]{#1}
\providecommand{\DOIprefix}{doi:}
\providecommand{\ArXivprefix}{arXiv:}
\providecommand{\URLprefix}{URL: }
\providecommand{\Pubmedprefix}{pmid:}
\providecommand{\doi}[1]{\href{http://dx.doi.org/#1}{\path{#1}}}
\providecommand{\Pubmed}[1]{\href{pmid:#1}{\path{#1}}}
\providecommand{\bibinfo}[2]{#2}
\ifx\xfnm\relax \def\xfnm[#1]{\unskip,\space#1}\fi
\bibitem[{Ajspur et~al.(2019)Ajspur, Jensen and Andersen}]{Ajspur2019AModels}
\bibinfo{author}{Ajspur, M.L.}, \bibinfo{author}{Jensen, R.M.},
  \bibinfo{author}{Andersen, K.H.}, \bibinfo{year}{2019}.
\newblock \bibinfo{title}{{A decomposed fourier-motzkin elimination framework
  to derive vessel capacity models}}. volume \bibinfo{volume}{11756 LNCS}.
\newblock \bibinfo{publisher}{Springer International Publishing}.
\newblock \URLprefix \url{http://dx.doi.org/10.1007/978-3-030-31140-7_6},
  \DOIprefix\doi{10.1007/978-3-030-31140-7{\_}6}.
\bibitem[{Ambrosino et~al.(2010)Ambrosino, Anghinolfi, Paolucci and
  Sciomachen}]{Ambrosino2010AnProblem}
\bibinfo{author}{Ambrosino, D.}, \bibinfo{author}{Anghinolfi, D.},
  \bibinfo{author}{Paolucci, M.}, \bibinfo{author}{Sciomachen, A.},
  \bibinfo{year}{2010}.
\newblock \bibinfo{title}{{An Experimental Comparison of Different Heuristics
  for the Master Bay Plan Problem}}, in: \bibinfo{booktitle}{Experimental
  Algorithms}. \bibinfo{publisher}{Springer Berlin Heidelberg}. volume
  \bibinfo{volume}{6049}, pp. \bibinfo{pages}{314--325}.
\newblock \URLprefix
  \url{http://link.springer.com/10.1007/978-3-642-13193-6_27},
  \DOIprefix\doi{10.1007/978-3-642-13193-6{\_}27}.
\bibitem[{Ambrosino et~al.(2015a)Ambrosino, Paolucci and
  Sciomachen}]{Ambrosino2015APlanning}
\bibinfo{author}{Ambrosino, D.}, \bibinfo{author}{Paolucci, M.},
  \bibinfo{author}{Sciomachen, A.}, \bibinfo{year}{2015}a.
\newblock \bibinfo{title}{{A MIP Heuristic for Multi Port Stowage Planning}}.
\newblock \bibinfo{journal}{Transportation Research Procedia}
  \bibinfo{volume}{10}, \bibinfo{pages}{725--734}.
\newblock \URLprefix
  \url{http://linkinghub.elsevier.com/retrieve/pii/S2352146515002136},
  \DOIprefix\doi{10.1016/j.trpro.2015.09.026}.
\bibitem[{Ambrosino et~al.(2015b)Ambrosino, Paolucci and
  Sciomachen}]{Ambrosino2015ComputationalProblems}
\bibinfo{author}{Ambrosino, D.}, \bibinfo{author}{Paolucci, M.},
  \bibinfo{author}{Sciomachen, A.}, \bibinfo{year}{2015}b.
\newblock \bibinfo{title}{{Computational evaluation of a MIP model for
  multi-port stowage planning problems}}.
\newblock \bibinfo{journal}{Soft Computing}
  \DOIprefix\doi{10.1007/s00500-015-1879-y}.
\bibitem[{Ambrosino et~al.(2015c)Ambrosino, Paolucci and
  Sciomachen}]{Ambrosino2015ExperimentalProblem}
\bibinfo{author}{Ambrosino, D.}, \bibinfo{author}{Paolucci, M.},
  \bibinfo{author}{Sciomachen, A.}, \bibinfo{year}{2015}c.
\newblock \bibinfo{title}{{Experimental evaluation of mixed integer programming
  models for the multi-port master bay plan problem}}.
\newblock \bibinfo{journal}{Flexible Services and Manufacturing Journal}
  \bibinfo{volume}{27}, \bibinfo{pages}{263--284}.
\newblock \URLprefix \url{http://link.springer.com/10.1007/s10696-013-9185-4},
  \DOIprefix\doi{10.1007/s10696-013-9185-4}.
\bibitem[{Ambrosino et~al.(2018)Ambrosino, Paolucci and
  Sciomachen}]{Ambrosino2018ShippingApproach}
\bibinfo{author}{Ambrosino, D.}, \bibinfo{author}{Paolucci, M.},
  \bibinfo{author}{Sciomachen, A.}, \bibinfo{year}{2018}.
\newblock \bibinfo{title}{{Shipping Liner Company Stowage Plans: An
  Optimization Approach}}, in: \bibinfo{editor}{{\.{Z}}ak, J.},
  \bibinfo{editor}{Hadas, Y.}, \bibinfo{editor}{Rossi, R.} (Eds.),
  \bibinfo{booktitle}{Advances in Intelligent Systems and Computing}.
  \bibinfo{publisher}{Springer International Publishing},
  \bibinfo{address}{Cham}. volume \bibinfo{volume}{572} of
  \textit{\bibinfo{series}{Advances in Intelligent Systems and Computing}}, pp.
  \bibinfo{pages}{405--420}.
\newblock \URLprefix \url{http://link.springer.com/10.1007/978-3-319-57105-8},
  \DOIprefix\doi{10.1007/978-3-319-57105-8{\_}20}.
\bibitem[{Ambrosino et~al.(2004)Ambrosino, Sciomachen and
  Tanfani}]{Ambrosino2004StowingProblem}
\bibinfo{author}{Ambrosino, D.}, \bibinfo{author}{Sciomachen, A.},
  \bibinfo{author}{Tanfani, E.}, \bibinfo{year}{2004}.
\newblock \bibinfo{title}{{Stowing a containership: the master bay plan
  problem}}.
\newblock \bibinfo{journal}{Transportation Research Part A}
  \bibinfo{volume}{38}, \bibinfo{pages}{81--99}.
\newblock \URLprefix
  \url{http://linkinghub.elsevier.com/retrieve/pii/S0965856403000892
  http://www.sciencedirect.com/science/article/pii/S0965856403000892},
  \DOIprefix\doi{10.1016/j.tra.2003.09.002}.
\bibitem[{Ambrosino et~al.(2006)Ambrosino, Sciomachen and
  Tanfani}]{Ambrosino2006AProblem}
\bibinfo{author}{Ambrosino, D.}, \bibinfo{author}{Sciomachen, A.},
  \bibinfo{author}{Tanfani, E.}, \bibinfo{year}{2006}.
\newblock \bibinfo{title}{{A decomposition heuristics for the container ship
  stowage problem}}.
\newblock \bibinfo{journal}{Journal of Heuristics} \bibinfo{volume}{12},
  \bibinfo{pages}{211--233}.
\newblock \URLprefix \url{http://link.springer.com/10.1007/s10732-006-5905-1
  http://link.springer.com/article/10.1007/s10732-006-5905-1},
  \DOIprefix\doi{10.1007/s10732-006-5905-1}.
\bibitem[{Aslidis(1989)}]{Aslidis1989CombinatorialProblems}
\bibinfo{author}{Aslidis, A.H.}, \bibinfo{year}{1989}.
\newblock \bibinfo{title}{{Combinatorial algorithms for stacking problems}}.
\newblock Ph.D. thesis.
\bibitem[{Avriel and Penn(1993)}]{AvrielMordecaiPenn1993}
\bibinfo{author}{Avriel, M.}, \bibinfo{author}{Penn, M.}, \bibinfo{year}{1993}.
\newblock \bibinfo{title}{{Exact and approximate solutions of the container
  ship stowage problem}}.
\newblock \bibinfo{journal}{Industrial Engineering} \bibinfo{volume}{25},
  \bibinfo{pages}{271--274}.
\bibitem[{Avriel et~al.(2000)Avriel, Penn and
  Shpirer}]{Avriel2000ContainerGraphs}
\bibinfo{author}{Avriel, M.}, \bibinfo{author}{Penn, M.},
  \bibinfo{author}{Shpirer, N.}, \bibinfo{year}{2000}.
\newblock \bibinfo{title}{{Container ship stowage problem: complexity and
  connection to the coloring of circle graphs}}.
\newblock \bibinfo{journal}{Discrete Applied Mathematics}
  \bibinfo{volume}{103}, \bibinfo{pages}{271--279}.
\newblock \URLprefix
  \url{https://linkinghub.elsevier.com/retrieve/pii/S0166218X99002450},
  \DOIprefix\doi{10.1016/S0166-218X(99)00245-0}.
\bibitem[{Avriel et~al.(1998)Avriel, Penn, Shpirer and
  Witteboon}]{Avriel1998StowageShifts}
\bibinfo{author}{Avriel, M.}, \bibinfo{author}{Penn, M.},
  \bibinfo{author}{Shpirer, N.}, \bibinfo{author}{Witteboon, S.},
  \bibinfo{year}{1998}.
\newblock \bibinfo{title}{{Stowage planning for container ships to reduce the
  number of shifts}}.
\newblock \bibinfo{journal}{Annals of Operations Research}
  \bibinfo{volume}{76}, \bibinfo{pages}{55--71}.
\newblock \URLprefix
  \url{http://search.ebscohost.com/login.aspx?direct=true&db=bth&AN=18925445&site=ehost-live
  http://link.springer.com/article/10.1023/A:1018956823693}.
\bibitem[{Aye et~al.(2010)Aye, Low, Ying, Jing, Fan and
  Min}]{Aye2010VisualizationSystem}
\bibinfo{author}{Aye, W.C.}, \bibinfo{author}{Low, M.Y.H.},
  \bibinfo{author}{Ying, H.S.}, \bibinfo{author}{Jing, H.W.},
  \bibinfo{author}{Fan, L.}, \bibinfo{author}{Min, Z.}, \bibinfo{year}{2010}.
\newblock \bibinfo{title}{{Visualization and simulation tool for automated
  stowage plan generation system}}.
\newblock \bibinfo{journal}{Proceedings of the International MultiConference of
  Engineers and Computer Scientists 2010, IMECS 2010} \bibinfo{volume}{II},
  \bibinfo{pages}{1013--1019}.
\bibitem[{Azevedo et~al.(2014)Azevedo, Cassilda~Maria, de~Sena, Chaves, Neto
  and Moretti}]{Azevedo2014SolvingMeta-heuristics}
\bibinfo{author}{Azevedo, A.}, \bibinfo{author}{Cassilda~Maria, R.},
  \bibinfo{author}{de~Sena, G.J.}, \bibinfo{author}{Chaves, A.A.},
  \bibinfo{author}{Neto, L.L.S.}, \bibinfo{author}{Moretti, A.C.},
  \bibinfo{year}{2014}.
\newblock \bibinfo{title}{{Solving the 3D container ship loading planning
  problem by representation by rules and meta-heuristics}}.
\newblock \bibinfo{journal}{International Journal of Data Analysis Techniques
  and Strategies} \bibinfo{volume}{6}, \bibinfo{pages}{228--260}.
\newblock \DOIprefix\doi{10.1504/IJDATS.2014.063060}.
\bibitem[{Azevedo et~al.(2018)Azevedo, de~Salles~Neto, Chaves and
  Moretti}]{Azevedo2018SolvingAlgorithm}
\bibinfo{author}{Azevedo, A.T.}, \bibinfo{author}{de~Salles~Neto, L.L.},
  \bibinfo{author}{Chaves, A.A.}, \bibinfo{author}{Moretti, A.C.},
  \bibinfo{year}{2018}.
\newblock \bibinfo{title}{{Solving the 3D stowage planning problem integrated
  with the quay crane scheduling problem by representation by rules and genetic
  algorithm}}.
\newblock \bibinfo{journal}{Applied Soft Computing Journal}
  \bibinfo{volume}{65}, \bibinfo{pages}{495--516}.
\newblock \DOIprefix\doi{10.1016/j.asoc.2018.01.006}.
\bibitem[{Bilican et~al.(2020)Bilican, Evren and
  Karatas}]{Bilican2020AParameters}
\bibinfo{author}{Bilican, M.S.}, \bibinfo{author}{Evren, R.},
  \bibinfo{author}{Karatas, M.}, \bibinfo{year}{2020}.
\newblock \bibinfo{title}{{A Mathematical Model and Two-Stage Heuristic for the
  Container Stowage Planning Problem with Stability Parameters}}.
\newblock \bibinfo{journal}{IEEE Access} \bibinfo{volume}{8},
  \bibinfo{pages}{113392--113413}.
\newblock \DOIprefix\doi{10.1109/ACCESS.2020.3003557}.
\bibitem[{Botter and Brinati(1992)}]{Botter1992StowageSolution}
\bibinfo{author}{Botter, R.}, \bibinfo{author}{Brinati, M.},
  \bibinfo{year}{1992}.
\newblock \bibinfo{title}{{Stowage container planning: a model for getting an
  optimal solution}}.
\newblock \bibinfo{journal}{Computer Applications in Automation of Shipyard
  Operation and Ship Design} \bibinfo{volume}{VII}, \bibinfo{pages}{217--228}.
\bibitem[{Chang et~al.(2022)Chang, Hamedi and Haghani}]{Chang2022SolvingMode}
\bibinfo{author}{Chang, Y.}, \bibinfo{author}{Hamedi, M.},
  \bibinfo{author}{Haghani, A.}, \bibinfo{year}{2022}.
\newblock \bibinfo{title}{{Solving integrated problem of stowage planning with
  crane split by an improved genetic algorithm based on novel encoding mode}}.
\newblock \bibinfo{journal}{Measurement and Control} \URLprefix
  \url{http://journals.sagepub.com/doi/10.1177/00202940221097981},
  \DOIprefix\doi{10.1177/00202940221097981}.
\bibitem[{Chao and Lin(2021)}]{Chao2021}
\bibinfo{author}{Chao, S.L.}, \bibinfo{author}{Lin, P.H.},
  \bibinfo{year}{2021}.
\newblock \bibinfo{title}{{Minimizing overstowage in master bay plans of large
  container ships}}.
\newblock \bibinfo{journal}{Maritime Economics and Logistics}
  \bibinfo{volume}{23}, \bibinfo{pages}{71--93}.
\newblock \URLprefix \url{https://doi.org/10.1057/s41278-019-00126-6},
  \DOIprefix\doi{10.1057/s41278-019-00126-6}.
\bibitem[{Cho(1981)}]{Cho1981DevelopmentPlanning}
\bibinfo{author}{Cho, D.}, \bibinfo{year}{1981}.
\newblock \bibinfo{title}{{Development of a methodology for containership load
  planning}}.
\newblock Ph.D. thesis. Oregon State University.
\bibitem[{Chou and Fang(2021)}]{Chou2021ApplyingStudy}
\bibinfo{author}{Chou, C.C.}, \bibinfo{author}{Fang, P.Y.},
  \bibinfo{year}{2021}.
\newblock \bibinfo{title}{{Applying expert knowledge to containership stowage
  planning: an empirical study}}.
\newblock \bibinfo{journal}{Maritime Economics and Logistics}
  \bibinfo{volume}{23}, \bibinfo{pages}{4--27}.
\newblock \URLprefix \url{https://doi.org/10.1057/s41278-018-0113-0},
  \DOIprefix\doi{10.1057/s41278-018-0113-0}.
\bibitem[{Christensen et~al.(2019)Christensen, Erera and
  Pacino}]{Christensen2019AProblem}
\bibinfo{author}{Christensen, J.}, \bibinfo{author}{Erera, A.},
  \bibinfo{author}{Pacino, D.}, \bibinfo{year}{2019}.
\newblock \bibinfo{title}{{A rolling horizon heuristic for the stochastic cargo
  mix problem}}.
\newblock \bibinfo{journal}{Transportation Research Part E: Logistics and
  Transportation Review} \bibinfo{volume}{123}.
\newblock \DOIprefix\doi{10.1016/j.tre.2018.10.010}.
\bibitem[{Christensen and Pacino(2017)}]{Christensen2017AStowage}
\bibinfo{author}{Christensen, J.}, \bibinfo{author}{Pacino, D.},
  \bibinfo{year}{2017}.
\newblock \bibinfo{title}{{A matheuristic for the Cargo Mix Problem with Block
  Stowage}}.
\newblock \bibinfo{journal}{Transportation Research Part E: Logistics and
  Transportation Review} \bibinfo{volume}{97}.
\newblock \DOIprefix\doi{10.1016/j.tre.2016.10.005}.
\bibitem[{Conca et~al.(2018)Conca, Febbraro, Giglio and
  Rebora}]{Conca2018AutomationPlanning}
\bibinfo{author}{Conca, A.}, \bibinfo{author}{Febbraro, A.D.},
  \bibinfo{author}{Giglio, D.}, \bibinfo{author}{Rebora, F.},
  \bibinfo{year}{2018}.
\newblock \bibinfo{title}{{Automation in freight port call process: Real time
  data sharing to improve the stowage planning}}, in:
  \bibinfo{booktitle}{Transportation Research Procedia}.
\newblock \DOIprefix\doi{10.1016/j.trpro.2018.09.009}.
\bibitem[{Cruz-reyes et~al.(2013)Cruz-reyes, H, Melin, H and
  O}]{Cruz-reyes2013ConstructivePlanning}
\bibinfo{author}{Cruz-reyes, L.}, \bibinfo{author}{H, P.H.},
  \bibinfo{author}{Melin, P.}, \bibinfo{author}{H, H.J.F.}, \bibinfo{author}{O,
  J.M.}, \bibinfo{year}{2013}.
\newblock \bibinfo{title}{{Constructive Algorithm for a Benchmark in Ship
  Stowage Planning}} , \bibinfo{pages}{393--408}.
\bibitem[{Cruz-Reyes et~al.(2015)Cruz-Reyes, Hern{\'{a}}ndez, Melin,
  Joaqu{\'{i}}n, Huacuja, Jos{\'{e}}, Soberanes, Javier, M{\'{e}}xico,
  Tecnol{\'{o}}gico, Madero, M{\'{e}}xico, Tijuana, California, Aut{\'{o}}noma,
  Tampico, M{\'{e}}xico and Le{\'{o}}n}]{Cruz-Reyes2015}
\bibinfo{author}{Cruz-Reyes, L.}, \bibinfo{author}{Hern{\'{a}}ndez, P.H.},
  \bibinfo{author}{Melin, P.}, \bibinfo{author}{Joaqu{\'{i}}n, H.},
  \bibinfo{author}{Huacuja, F.}, \bibinfo{author}{Jos{\'{e}}, H.},
  \bibinfo{author}{Soberanes, P.}, \bibinfo{author}{Javier, J.},
  \bibinfo{author}{M{\'{e}}xico, T.N.D.}, \bibinfo{author}{Tecnol{\'{o}}gico,
  I.}, \bibinfo{author}{Madero, D.C.}, \bibinfo{author}{M{\'{e}}xico, T.N.D.},
  \bibinfo{author}{Tijuana, I.T.D.}, \bibinfo{author}{California, B.},
  \bibinfo{author}{Aut{\'{o}}noma, U.}, \bibinfo{author}{Tampico, D.T.},
  \bibinfo{author}{M{\'{e}}xico, T.N.D.}, \bibinfo{author}{Le{\'{o}}n, I.T.D.},
  \bibinfo{year}{2015}.
\newblock \bibinfo{title}{{Lower and Upper Bounds for the Master Bay Planning
  Problem}}.
\newblock \bibinfo{journal}{International Journal of Combinatorial Optimization
  Problems and Informatics} \bibinfo{volume}{6}, \bibinfo{pages}{42--52}.
\bibitem[{Delgado et~al.(2012a)Delgado, Jensen and
  Guilbert}]{Delgado2012AStowage}
\bibinfo{author}{Delgado, A.}, \bibinfo{author}{Jensen, R.M.},
  \bibinfo{author}{Guilbert, N.}, \bibinfo{year}{2012}a.
\newblock \bibinfo{title}{{A placement heuristic for a commercial decision
  support system for container vessel stowage}}.
\newblock \bibinfo{journal}{38th Latin America Conference on Informatics, CLEI
  2012 - Conference Proceedings} \DOIprefix\doi{10.1109/CLEI.2012.6427181}.
\bibitem[{Delgado et~al.(2012b)Delgado, Jensen, Janstrup, Rose and
  Andersen}]{Delgado2012ABays}
\bibinfo{author}{Delgado, A.}, \bibinfo{author}{Jensen, R.M.},
  \bibinfo{author}{Janstrup, K.}, \bibinfo{author}{Rose, T.H.},
  \bibinfo{author}{Andersen, K.H.}, \bibinfo{year}{2012}b.
\newblock \bibinfo{title}{{A Constraint Programming model for fast optimal
  stowage of container vessel bays}}.
\newblock \bibinfo{journal}{European Journal of Operational Research}
  \bibinfo{volume}{220}, \bibinfo{pages}{251--261}.
\newblock \URLprefix
  \url{https://linkinghub.elsevier.com/retrieve/pii/S0377221712000483},
  \DOIprefix\doi{10.1016/j.ejor.2012.01.028}.
\bibitem[{Delgado et~al.(2009)Delgado, Jensen and Schulte}]{Delgado2009}
\bibinfo{author}{Delgado, A.}, \bibinfo{author}{Jensen, R.M.},
  \bibinfo{author}{Schulte, C.}, \bibinfo{year}{2009}.
\newblock \bibinfo{title}{{Generating optimal stowage plans for container
  vessel bays}}, in: \bibinfo{booktitle}{Lecture Notes in Computer Science
  (including subseries Lecture Notes in Artificial Intelligence and Lecture
  Notes in Bioinformatics)}, pp. \bibinfo{pages}{6--20}.
\bibitem[{Ding and Chou(2015)}]{Ding2015StowageShifts}
\bibinfo{author}{Ding, D.}, \bibinfo{author}{Chou, M.C.}, \bibinfo{year}{2015}.
\newblock \bibinfo{title}{{Stowage Planning for Container Ships: A Heuristic
  Algorithm to Reduce the Number of Shifts}}.
\newblock \bibinfo{journal}{European Journal of Operational Research}
  \URLprefix
  \url{http://linkinghub.elsevier.com/retrieve/pii/S0377221715002660},
  \DOIprefix\doi{10.1016/j.ejor.2015.03.044}.
\bibitem[{Dubrovsky et~al.(2002)Dubrovsky, Levitin and
  Penn}]{Dubrovsky2002AProblem}
\bibinfo{author}{Dubrovsky, O.}, \bibinfo{author}{Levitin, G.},
  \bibinfo{author}{Penn, M.}, \bibinfo{year}{2002}.
\newblock \bibinfo{title}{{A genetic algorithm with a compact solution encoding
  for the container ship stowage problem}}.
\newblock \bibinfo{journal}{Journal of Heuristics} \bibinfo{volume}{8},
  \bibinfo{pages}{585--599}.
\newblock \URLprefix
  \url{http://link.springer.com/article/10.1023/A:1020373709350},
  \DOIprefix\doi{10.1023/A:1020373709350}.
\bibitem[{El~Yaagoubi et~al.(2022)El~Yaagoubi, Charhbili, Boukachour and
  El~Hilali~Alaoui}]{ElYaagoubi2022Multi-objectiveSystem}
\bibinfo{author}{El~Yaagoubi, A.}, \bibinfo{author}{Charhbili, M.},
  \bibinfo{author}{Boukachour, J.}, \bibinfo{author}{El~Hilali~Alaoui, A.},
  \bibinfo{year}{2022}.
\newblock \bibinfo{title}{{Multi-objective optimization of the 3D container
  stowage planning problem in a barge convoy system}}.
\newblock \bibinfo{journal}{Computers {\&} Operations Research} ,
  \bibinfo{pages}{105796}\URLprefix
  \url{https://linkinghub.elsevier.com/retrieve/pii/S0305054822000855},
  \DOIprefix\doi{10.1016/J.COR.2022.105796}.
\bibitem[{Hamedi(2011)}]{Hamedi2011CONTAINERSHIPOPERATIONS}
\bibinfo{author}{Hamedi, M.}, \bibinfo{year}{2011}.
\newblock \bibinfo{title}{{CONTAINERSHIP LOAD PLANNING WITH CRANE OPERATIONS}}.
\newblock Ph.D. thesis. University of Maryland.
\bibitem[{Hsu et~al.(2021)Hsu, Wang, Fu and Dang}]{Hsu2021JointApproach}
\bibinfo{author}{Hsu, H.P.}, \bibinfo{author}{Wang, C.N.}, \bibinfo{author}{Fu,
  H.P.}, \bibinfo{author}{Dang, T.T.}, \bibinfo{year}{2021}.
\newblock \bibinfo{title}{{Joint scheduling of yard crane, yard truck, and quay
  crane for container terminal considering vessel stowage plan: An integrated
  simulation-based optimization approach}}.
\newblock \bibinfo{journal}{Mathematics} \bibinfo{volume}{9}.
\newblock \DOIprefix\doi{10.3390/math9182236}.
\bibitem[{Hu and Cai(2017)}]{Hu2017Multi-objectiveRoute}
\bibinfo{author}{Hu, M.}, \bibinfo{author}{Cai, W.}, \bibinfo{year}{2017}.
\newblock \bibinfo{title}{{Multi-objective optimization based on improved
  genetic algorithm for containership stowage on full route}}.
\newblock \bibinfo{journal}{2017 4th International Conference on Industrial
  Engineering and Applications, ICIEA 2017} ,
  \bibinfo{pages}{224--228}\DOIprefix\doi{10.1109/IEA.2017.7939211}.
\bibitem[{Hu et~al.(2012)Hu, Hu, Shi, Luo and
  Song}]{Hu2012CombinatorialTerminal}
\bibinfo{author}{Hu, W.}, \bibinfo{author}{Hu, Z.}, \bibinfo{author}{Shi, L.},
  \bibinfo{author}{Luo, P.}, \bibinfo{author}{Song, W.}, \bibinfo{year}{2012}.
\newblock \bibinfo{title}{{Combinatorial Optimization and Strategy for Ship
  Stowage and Loading Schedule of Container Terminal}}.
\newblock \bibinfo{journal}{Journal of Computers} \bibinfo{volume}{7},
  \bibinfo{pages}{2078--2092}.
\newblock \URLprefix
  \url{http://ojs.academypublisher.com/index.php/jcp/article/view/8182},
  \DOIprefix\doi{10.4304/jcp.7.8.2078-2092}.
\bibitem[{{International Chamber of
  Shipping}(2023)}]{InternationalChamberofShipping2023EnvironmentalTransport}
\bibinfo{author}{{International Chamber of Shipping}}, \bibinfo{year}{2023}.
\newblock \bibinfo{title}{{Environmental Performance: Comparison of CO2
  Emissions by Different Modes of Transport}}.
\newblock \URLprefix
  \url{https://www.ics-shipping.org/shipping-fact/environmental-performance-environmental-performance/}.
\bibitem[{Iris et~al.(2018)Iris, Christensen, Pacino and
  Ropke}]{Iris2018FlexibleScheduling}
\bibinfo{author}{Iris, C.}, \bibinfo{author}{Christensen, J.},
  \bibinfo{author}{Pacino, D.}, \bibinfo{author}{Ropke, S.},
  \bibinfo{year}{2018}.
\newblock \bibinfo{title}{{Flexible ship loading problem with transfer vehicle
  assignment and scheduling}}.
\newblock \bibinfo{journal}{Transportation Research Part B: Methodological}
  \bibinfo{volume}{111}, \bibinfo{pages}{113--134}.
\newblock \DOIprefix\doi{10.1016/j.trb.2018.03.009}.
\bibitem[{Jensen and Ajspur(2018)}]{Jensen2018TheCapacity}
\bibinfo{author}{Jensen, R.M.}, \bibinfo{author}{Ajspur, M.L.},
  \bibinfo{year}{2018}.
\newblock \bibinfo{title}{{The Standard Capacity Model: Towards a Polyhedron
  Representation of Container Vessel Capacity}}, in:
  \bibinfo{booktitle}{Computational Logistics}. \bibinfo{publisher}{Springer
  International Publishing}. volume \bibinfo{volume}{8197}, pp.
  \bibinfo{pages}{175--190}.
\newblock \URLprefix \url{http://link.springer.com/10.1007/978-3-642-41019-2},
  \DOIprefix\doi{10.1007/978-3-030-00898-7}.
\bibitem[{Jensen and Ajspur(2022)}]{Jensen2022RevenueChallenge}
\bibinfo{author}{Jensen, R.M.}, \bibinfo{author}{Ajspur, M.L.},
  \bibinfo{year}{2022}.
\newblock \bibinfo{title}{{Revenue management in liner shipping: Addressing the
  vessel capacity challenge}}.
\newblock \bibinfo{journal}{Maritime Transport Research} \bibinfo{volume}{3}.
\newblock \DOIprefix\doi{10.1016/j.martra.2022.100069}.
\bibitem[{Jensen et~al.(2012)Jensen, Leknes and
  Bebbington}]{Jensen2012FastDiagrams}
\bibinfo{author}{Jensen, R.M.}, \bibinfo{author}{Leknes, E.},
  \bibinfo{author}{Bebbington, T.}, \bibinfo{year}{2012}.
\newblock \bibinfo{title}{{Fast interactive decision support for modifying
  stowage plans using binary decision diagrams}}, in:
  \bibinfo{booktitle}{Lecture Notes in Engineering and Computer Science}, pp.
  \bibinfo{pages}{1555--1561}.
\bibitem[{Jensen et~al.(2018)Jensen, Pacino, Ajspur and
  Vesterdal}]{Jensen2018ContainerPlanning}
\bibinfo{author}{Jensen, R.M.}, \bibinfo{author}{Pacino, D.},
  \bibinfo{author}{Ajspur, M.L.}, \bibinfo{author}{Vesterdal, C.},
  \bibinfo{year}{2018}.
\newblock \bibinfo{title}{{Container Vessel Stowage Planning}}.
\newblock \bibinfo{publisher}{Weilbach}.
\bibitem[{Jin and Mi(2019)}]{Jin2019AnBay}
\bibinfo{author}{Jin, J.}, \bibinfo{author}{Mi, W.}, \bibinfo{year}{2019}.
\newblock \bibinfo{title}{{An AIMMS-based decision-making model for optimizing
  the intelligent stowage of export containers in a single bay}}.
\newblock \bibinfo{journal}{Discrete and Continuous Dynamical Systems - Series
  S} \bibinfo{volume}{12}, \bibinfo{pages}{1101--1115}.
\newblock \DOIprefix\doi{10.3934/dcdss.2019076}.
\bibitem[{Kaisar(2006)}]{Kaisar2006ATRANSPORTATION}
\bibinfo{author}{Kaisar, E.I.}, \bibinfo{year}{2006}.
\newblock \bibinfo{title}{{A STOWAGE PLANNING MODEL FOR MULTIPORT CONTAINER
  TRANSPORTATION}}.
\newblock Ph.D. thesis.
\newblock \URLprefix \url{http://drum.lib.umd.edu/handle/1903/9139}.
\bibitem[{Kang and Kim(2002)}]{Kang2002StowageTransportation}
\bibinfo{author}{Kang, J.G.}, \bibinfo{author}{Kim, Y.D.},
  \bibinfo{year}{2002}.
\newblock \bibinfo{title}{{Stowage planning in maritime container
  transportation}}.
\newblock \bibinfo{journal}{Journal of the Operational Research Society}
  \bibinfo{volume}{53}, \bibinfo{pages}{415--426}.
\newblock \URLprefix
  \url{http://www.ingentaconnect.com/content/pal/01605682/2002/00000053/00000004/2601322}.
\bibitem[{Kebedow and Oppen(2018)}]{Kebedow2018IncludingProblem}
\bibinfo{author}{Kebedow, K.G.}, \bibinfo{author}{Oppen, J.},
  \bibinfo{year}{2018}.
\newblock \bibinfo{title}{{Including containers with dangerous goods in the
  multi-port master bay planning problem}}.
\newblock \bibinfo{journal}{Mendel} \bibinfo{volume}{24},
  \bibinfo{pages}{23--36}.
\newblock \DOIprefix\doi{10.13164/mendel.2018.2.023}.
\bibitem[{Kebedow and Oppen(2019a)}]{Kebedow2019IncludingStowage}
\bibinfo{author}{Kebedow, K.G.}, \bibinfo{author}{Oppen, J.},
  \bibinfo{year}{2019}a.
\newblock \bibinfo{title}{{Including containers with dangerous goods in the
  cargo mix problem for container vessel stowage}}.
\newblock \bibinfo{journal}{Communications - Scientific Letters of the
  University of Zilina} \bibinfo{volume}{21}, \bibinfo{pages}{100--113}.
\newblock \DOIprefix\doi{10.26552/com.c.2019.2.100-113}.
\bibitem[{Kebedow and Oppen(2019b)}]{Kebedow2019IncludingProblem}
\bibinfo{author}{Kebedow, K.G.}, \bibinfo{author}{Oppen, J.},
  \bibinfo{year}{2019}b.
\newblock \bibinfo{title}{{Including containers with dangerous goods in the
  slot planning problem}}.
\newblock \bibinfo{journal}{Proceedings of the International Conference on
  Industrial Engineering and Operations Management} \bibinfo{volume}{2019},
  \bibinfo{pages}{225--232}.
\bibitem[{Korach et~al.(2020)Korach, Brouer and
  Jensen}]{Korach2020MatheuristicsBays}
\bibinfo{author}{Korach, A.}, \bibinfo{author}{Brouer, B.D.},
  \bibinfo{author}{Jensen, R.M.}, \bibinfo{year}{2020}.
\newblock \bibinfo{title}{{Matheuristics for slot planning of container vessel
  bays}}.
\newblock \bibinfo{journal}{European Journal of Operational Research}
  \bibinfo{volume}{282}, \bibinfo{pages}{873--885}.
\newblock \URLprefix \url{https://doi.org/10.1016/j.ejor.2019.09.042},
  \DOIprefix\doi{10.1016/j.ejor.2019.09.042}.
\bibitem[{Kroer et~al.(2012)Kroer, Kjaer~Svendsen, M{\o}ller~Jensen and
  Kiniry}]{Kroer2012SATPlanning}
\bibinfo{author}{Kroer, C.}, \bibinfo{author}{Kjaer~Svendsen, M.},
  \bibinfo{author}{M{\o}ller~Jensen, R.}, \bibinfo{author}{Kiniry, J.R.},
  \bibinfo{year}{2012}.
\newblock \bibinfo{title}{{SAT and SMT-based Interactive Configuration for
  Container Vessel Stowage Planning}}.
\newblock \bibinfo{type}{Technical Report}. IT University of Copenhagen.
\bibitem[{Kroer et~al.(2016)Kroer, Svendsen, Jensen, Kiniry and
  Leknes}]{Kroer2016SymbolicPlanning}
\bibinfo{author}{Kroer, C.}, \bibinfo{author}{Svendsen, M.K.},
  \bibinfo{author}{Jensen, R.M.}, \bibinfo{author}{Kiniry, J.},
  \bibinfo{author}{Leknes, E.}, \bibinfo{year}{2016}.
\newblock \bibinfo{title}{{Symbolic configuration for interactive container
  ship stowage planning}}.
\newblock \bibinfo{journal}{Computational Intelligence} \bibinfo{volume}{32},
  \bibinfo{pages}{259--283}.
\newblock \DOIprefix\doi{10.1111/coin.12051}.
\bibitem[{Larsen and Pacino(2021)}]{Larsen2021AProblem}
\bibinfo{author}{Larsen, R.}, \bibinfo{author}{Pacino, D.},
  \bibinfo{year}{2021}.
\newblock \bibinfo{title}{{A heuristic and a benchmark for the stowage planning
  problem}}.
\newblock \bibinfo{journal}{Maritime Economics and Logistics}
  \bibinfo{volume}{23}, \bibinfo{pages}{94--122}.
\newblock \DOIprefix\doi{10.1057/s41278-020-00172-5}.
\bibitem[{Lee et~al.(2020)Lee, Lee and Shin}]{Lee2020LashingContainerships}
\bibinfo{author}{Lee, C.}, \bibinfo{author}{Lee, M.K.}, \bibinfo{author}{Shin,
  J.Y.}, \bibinfo{year}{2020}.
\newblock \bibinfo{title}{{Lashing Force Prediction Model with Multimodal Deep
  Learning and AutoML for Stowage Planning Automation in Containerships}}.
\newblock \bibinfo{journal}{Logistics} \bibinfo{volume}{5}, \bibinfo{pages}{1}.
\newblock \DOIprefix\doi{10.3390/logistics5010001}.
\bibitem[{Li et~al.(2020)Li, Zhang, Liu and Liang}]{Li2020OptimizingShipping}
\bibinfo{author}{Li, J.}, \bibinfo{author}{Zhang, Y.}, \bibinfo{author}{Liu,
  Z.}, \bibinfo{author}{Liang, X.}, \bibinfo{year}{2020}.
\newblock \bibinfo{title}{{Optimizing the Stowage Planning and Block Relocation
  Problem in Inland Container Shipping}}.
\newblock \bibinfo{journal}{IEEE Access} \bibinfo{volume}{8},
  \bibinfo{pages}{207499--207514}.
\newblock \DOIprefix\doi{10.1109/ACCESS.2020.3037675}.
\bibitem[{Li et~al.(2018)Li, Zhang, Ma and Ji}]{Li2018Multi-PortUncertainties}
\bibinfo{author}{Li, J.}, \bibinfo{author}{Zhang, Y.}, \bibinfo{author}{Ma,
  J.}, \bibinfo{author}{Ji, S.}, \bibinfo{year}{2018}.
\newblock \bibinfo{title}{{Multi-Port Stowage Planning for Inland Container
  Liner Shipping Considering Weight Uncertainties}}.
\newblock \bibinfo{journal}{IEEE Access} \bibinfo{volume}{6},
  \bibinfo{pages}{66468--66480}.
\newblock \DOIprefix\doi{10.1109/ACCESS.2018.2878308}.
\bibitem[{Liu et~al.(2011)Liu, Low, Hsu, Huang, Zeng and
  Win}]{Liu2011RandomizedPlans}
\bibinfo{author}{Liu, F.}, \bibinfo{author}{Low, M.Y.H.}, \bibinfo{author}{Hsu,
  W.J.}, \bibinfo{author}{Huang, S.Y.}, \bibinfo{author}{Zeng, M.},
  \bibinfo{author}{Win, C.A.}, \bibinfo{year}{2011}.
\newblock \bibinfo{title}{{Randomized algorithm with tabu search for
  multi-objective optimization of large containership stowage plans}}.
\newblock \bibinfo{journal}{Lecture Notes in Computer Science (including
  subseries Lecture Notes in Artificial Intelligence and Lecture Notes in
  Bioinformatics)} \bibinfo{volume}{6971 LNCS}, \bibinfo{pages}{256--272}.
\newblock \DOIprefix\doi{10.1007/978-3-642-24264-9{\_}20}.
\bibitem[{Martin et~al.(1988)Martin, Randhawa and
  McDowell}]{Martin1988ComputerizedEvaluation}
\bibinfo{author}{Martin, G.L.}, \bibinfo{author}{Randhawa, S.U.},
  \bibinfo{author}{McDowell, E.D.}, \bibinfo{year}{1988}.
\newblock \bibinfo{title}{{Computerized container-ship load planning: A
  methodology and evaluation}}.
\newblock \bibinfo{journal}{Computers and Industrial Engineering}
  \bibinfo{volume}{14}, \bibinfo{pages}{429--440}.
\newblock \DOIprefix\doi{10.1016/0360-8352(88)90045-9}.
\bibitem[{Monaco et~al.(2014)Monaco, Sammarra and
  Sorrentino}]{Monaco2014TheProblem}
\bibinfo{author}{Monaco, M.F.}, \bibinfo{author}{Sammarra, M.},
  \bibinfo{author}{Sorrentino, G.}, \bibinfo{year}{2014}.
\newblock \bibinfo{title}{{The Terminal-Oriented Ship Stowage Planning
  Problem}}.
\newblock \bibinfo{journal}{European Journal of Operational Research}
  \URLprefix
  \url{http://linkinghub.elsevier.com/retrieve/pii/S0377221714004536},
  \DOIprefix\doi{10.1016/j.ejor.2014.05.030}.
\bibitem[{Nugroho et~al.(2021)Nugroho, Djatmiko, {Murdjito}, Ardhi, Supomo and
  Buana}]{Nugroho2021RegulatoryConsiderations}
\bibinfo{author}{Nugroho, S.}, \bibinfo{author}{Djatmiko, E.B.},
  \bibinfo{author}{{Murdjito}}, \bibinfo{author}{Ardhi, E.W.},
  \bibinfo{author}{Supomo, H.}, \bibinfo{author}{Buana, I.G.N.S.},
  \bibinfo{year}{2021}.
\newblock \bibinfo{title}{{Regulatory framework of a computer-based stowage
  planning: safety and efficiency considerations}}.
\newblock \bibinfo{journal}{IOP Conference Series: Materials Science and
  Engineering} \bibinfo{volume}{1052}, \bibinfo{pages}{012065}.
\newblock \DOIprefix\doi{10.1088/1757-899x/1052/1/012065}.
\bibitem[{Pacino(2013)}]{Pacino2013AnPlanning}
\bibinfo{author}{Pacino, D.}, \bibinfo{year}{2013}.
\newblock \bibinfo{title}{{An LNS Approach for Container Stowage Multi-port
  Master Planning}}, in: \bibinfo{booktitle}{Computational Logistics}, pp.
  \bibinfo{pages}{35--44}.
\newblock \URLprefix
  \url{http://link.springer.com/10.1007/978-3-642-41019-2_3},
  \DOIprefix\doi{10.1007/978-3-642-41019-2{\_}3}.
\bibitem[{Pacino(2018)}]{Pacino2018CranePlanning}
\bibinfo{author}{Pacino, D.}, \bibinfo{year}{2018}.
\newblock \bibinfo{title}{{Crane Intensity and Block Stowage Strategies in
  Stowage Planning}}, in: \bibinfo{booktitle}{Lecture Notes in Computer Science
  (including subseries Lecture Notes in Artificial Intelligence and Lecture
  Notes in Bioinformatics)}. \bibinfo{publisher}{Springer International
  Publishing}. volume \bibinfo{volume}{11184 LNCS}, pp.
  \bibinfo{pages}{191--206}.
\newblock \URLprefix
  \url{http://link.springer.com/10.1007/978-3-030-00898-7_12},
  \DOIprefix\doi{10.1007/978-3-030-00898-7{\_}12}.
\bibitem[{Pacino et~al.(2011)Pacino, Delgado, Jensen and
  Bebbington}]{Pacino2011FastVessels}
\bibinfo{author}{Pacino, D.}, \bibinfo{author}{Delgado, A.},
  \bibinfo{author}{Jensen, R.}, \bibinfo{author}{Bebbington, T.},
  \bibinfo{year}{2011}.
\newblock \bibinfo{title}{{Fast generation of near-optimal plans for
  eco-efficient stowage of large container vessels}}.
\newblock \bibinfo{journal}{Computational Logistics} ,
  \bibinfo{pages}{286--301}\URLprefix
  \url{http://link.springer.com/chapter/10.1007/978-3-642-24264-9_22}.
\bibitem[{Pacino et~al.(2012)Pacino, Delgado, Jensen and
  Bebbington}]{Pacino2012AnTanks}
\bibinfo{author}{Pacino, D.}, \bibinfo{author}{Delgado, A.},
  \bibinfo{author}{Jensen, R.M.}, \bibinfo{author}{Bebbington, T.},
  \bibinfo{year}{2012}.
\newblock \bibinfo{title}{{An accurate model for seaworthy container vessel
  stowage planning with ballast tanks}}.
\newblock \bibinfo{journal}{Lecture Notes in Computer Science (including
  subseries Lecture Notes in Artificial Intelligence and Lecture Notes in
  Bioinformatics)} \bibinfo{volume}{7555 LNCS}, \bibinfo{pages}{17--32}.
\newblock \DOIprefix\doi{10.1007/978-3-642-33587-7{\_}2}.
\bibitem[{Pacino and Jensen(2010)}]{Pacino2010ABays}
\bibinfo{author}{Pacino, D.}, \bibinfo{author}{Jensen, R.M.},
  \bibinfo{year}{2010}.
\newblock \bibinfo{title}{{A 3-Phase Randomized Constraint Based Local Search
  Algorithm for Stowing Under Deck Locations of Container Vessel Bays}}.
\newblock \URLprefix
  \url{https://en.itu.dk/-/media/EN/Research/About-ITU-Research/Technical-Reports/2010/ITU-TR-2010-123-pdf}.
\bibitem[{Pacino and Jensen(2013)}]{Pacino2013FastSearch}
\bibinfo{author}{Pacino, D.}, \bibinfo{author}{Jensen, R.M.},
  \bibinfo{year}{2013}.
\newblock \bibinfo{title}{{Fast slot planning using constraint-based local
  search}}.
\newblock \bibinfo{journal}{Lecture Notes in Electrical Engineering}
  \bibinfo{volume}{186 LNEE}, \bibinfo{pages}{49--63}.
\newblock \DOIprefix\doi{10.1007/978-94-007-5651-9-4}.
\bibitem[{Parre{\~{n}}o et~al.(2016)Parre{\~{n}}o, Pacino and
  Alvarez-Valdes}]{Parreno2016AProblem}
\bibinfo{author}{Parre{\~{n}}o, F.}, \bibinfo{author}{Pacino, D.},
  \bibinfo{author}{Alvarez-Valdes, R.}, \bibinfo{year}{2016}.
\newblock \bibinfo{title}{{A GRASP algorithm for the container stowage slot
  planning problem}}.
\newblock \bibinfo{journal}{Transportation Research Part E: Logistics and
  Transportation Review} \bibinfo{volume}{94}, \bibinfo{pages}{141--157}.
\newblock \DOIprefix\doi{10.1016/j.tre.2016.07.011}.
\bibitem[{Parre{\~{n}}o-Torres(2020)}]{Parreno-Torres2020ImprovingProblems}
\bibinfo{author}{Parre{\~{n}}o-Torres, C.}, \bibinfo{year}{2020}.
\newblock \bibinfo{title}{{Improving container terminal efficiency: New models
  and algorithms for Premarshalling and Stowage Problems}} \URLprefix
  \url{http://roderic.uv.es/handle/10550/75440}.
\bibitem[{Parre{\~{n}}o-Torres et~al.(2019)Parre{\~{n}}o-Torres, Alvarez-Valdes
  and Parre{\~{n}}o}]{Parreno-Torres2019SolutionProblem}
\bibinfo{author}{Parre{\~{n}}o-Torres, C.}, \bibinfo{author}{Alvarez-Valdes,
  R.}, \bibinfo{author}{Parre{\~{n}}o, F.}, \bibinfo{year}{2019}.
\newblock \bibinfo{title}{{Solution strategies for a multiport container ship
  stowage problem}}.
\newblock \bibinfo{journal}{Mathematical Problems in Engineering}
  \bibinfo{volume}{2019}.
\newblock \DOIprefix\doi{10.1155/2019/9029267}.
\bibitem[{Parre{\~{n}}o-Torres et~al.(2021)Parre{\~{n}}o-Torres, {\c{C}}alık,
  Alvarez-Valdes and Ruiz}]{Parreno-Torres2021SolvingAlgorithm}
\bibinfo{author}{Parre{\~{n}}o-Torres, C.}, \bibinfo{author}{{\c{C}}alık, H.},
  \bibinfo{author}{Alvarez-Valdes, R.}, \bibinfo{author}{Ruiz, R.},
  \bibinfo{year}{2021}.
\newblock \bibinfo{title}{{Solving the generalized multi-port container stowage
  planning problem by a matheuristic algorithm}}.
\newblock \bibinfo{journal}{Computers and Operations Research}
  \bibinfo{volume}{133}, \bibinfo{pages}{105383}.
\newblock \URLprefix \url{https://doi.org/10.1016/j.cor.2021.105383},
  \DOIprefix\doi{10.1016/j.cor.2021.105383}.
\bibitem[{Rashed et~al.(2021)Rashed, Eltawil and Gheith}]{Rashed2021AVessels}
\bibinfo{author}{Rashed, D.}, \bibinfo{author}{Eltawil, A.},
  \bibinfo{author}{Gheith, M.}, \bibinfo{year}{2021}.
\newblock \bibinfo{title}{{A Fuzzy Logic-Based Algorithm to Solve the Slot
  Planning Problem in Container Vessels}}.
\newblock \bibinfo{journal}{Logistics} \bibinfo{volume}{5},
  \bibinfo{pages}{67}.
\newblock \DOIprefix\doi{10.3390/logistics5040067}.
\bibitem[{Roberti and Pacino(2018)}]{Roberti2018APlans}
\bibinfo{author}{Roberti, R.}, \bibinfo{author}{Pacino, D.},
  \bibinfo{year}{2018}.
\newblock \bibinfo{title}{{A decomposition method for finding optimal container
  stowage plans}}.
\newblock \bibinfo{journal}{Transportation Science} \bibinfo{volume}{52},
  \bibinfo{pages}{1444--1462}.
\newblock \DOIprefix\doi{10.1287/trsc.2017.0795}.
\bibitem[{Saginaw and Perakis(1989)}]{Saginaw1989DecisionPlanning}
\bibinfo{author}{Saginaw, D.J.}, \bibinfo{author}{Perakis, A.N.},
  \bibinfo{year}{1989}.
\newblock \bibinfo{title}{{Decision support system for containership stowage
  planning}}.
\newblock \bibinfo{journal}{Marine Technology and SNAME News}
  \bibinfo{volume}{26}, \bibinfo{pages}{47--61}.
\newblock \DOIprefix\doi{10.5957/mt1.1989.26.1.47}.
\bibitem[{Sciomachen and Tanfani(2003)}]{Sciomachen2003TheProblem}
\bibinfo{author}{Sciomachen, A.}, \bibinfo{author}{Tanfani, E.},
  \bibinfo{year}{2003}.
\newblock \bibinfo{title}{{The master bay plan problem: A solution method based
  on its connection to the three-dimensional bin packing problem}}.
\newblock \bibinfo{journal}{IMA Journal of Management Mathematics}
  \bibinfo{volume}{14}, \bibinfo{pages}{251--269}.
\newblock \DOIprefix\doi{10.1093/imaman/14.3.251}.
\bibitem[{Sciomachen and Tanfani(2007)}]{Sciomachen2007AProductivity}
\bibinfo{author}{Sciomachen, A.}, \bibinfo{author}{Tanfani, E.},
  \bibinfo{year}{2007}.
\newblock \bibinfo{title}{{A 3D-BPP approach for optimising stowage plans and
  terminal productivity}}.
\newblock \bibinfo{journal}{European Journal of Operational Research}
  \bibinfo{volume}{183}, \bibinfo{pages}{1433--1446}.
\newblock \DOIprefix\doi{10.1016/j.ejor.2005.11.067}.
\bibitem[{Serban and Carp(2017)}]{Serban2017AStrategy}
\bibinfo{author}{Serban, C.}, \bibinfo{author}{Carp, D.}, \bibinfo{year}{2017}.
\newblock \bibinfo{title}{{A genetic algorithm for solving a container storage
  problem using a residence time strategy}}.
\newblock \bibinfo{journal}{Studies in Informatics and Control}
  \bibinfo{volume}{26}, \bibinfo{pages}{59--66}.
\newblock \DOIprefix\doi{10.24846/v26i1y201707}.
\bibitem[{Shen et~al.(2017)Shen, Zhao, Xia and Du}]{Shen2017AProblem}
\bibinfo{author}{Shen, Y.}, \bibinfo{author}{Zhao, N.}, \bibinfo{author}{Xia,
  M.}, \bibinfo{author}{Du, X.}, \bibinfo{year}{2017}.
\newblock \bibinfo{title}{{A deep Q-learning network for ship stowage planning
  problem}}.
\newblock \bibinfo{journal}{Polish Maritime Research} \bibinfo{volume}{24},
  \bibinfo{pages}{102--109}.
\newblock \DOIprefix\doi{10.1515/pomr-2017-0111}.
\bibitem[{Shields(1984)}]{Shields1984}
\bibinfo{author}{Shields, J.}, \bibinfo{year}{1984}.
\newblock \bibinfo{title}{{Containership Stowage: A Computer-Aided Preplanning
  System}} \URLprefix \url{http://trid.trb.org/view.aspx?id=419881}.
\bibitem[{Song et~al.(2010)Song, Dou, Ren and Liu}]{Song2010ResearchLogistics}
\bibinfo{author}{Song, X.y.}, \bibinfo{author}{Dou, X.c.},
  \bibinfo{author}{Ren, Y.}, \bibinfo{author}{Liu, X.}, \bibinfo{year}{2010}.
\newblock \bibinfo{title}{{Research on application of simulation technology in
  container ship stowage problem of port logistics}}.
\newblock \bibinfo{journal}{Proceedings IE {\&} EM 2010 : 2010 IEEE 17th
  International Conference on Industrial Engineering and Engineering
  Management} , \bibinfo{pages}{29--31}.
\bibitem[{{The Economist}(2013)}]{TheEconomist2013FreeHero}
\bibinfo{author}{{The Economist}}, \bibinfo{year}{2013}.
\newblock \bibinfo{title}{{Free exchange - The humble hero}}.
\newblock \bibinfo{journal}{The Economist} \URLprefix
  \url{https://www.economist.com/finance-and-economics/2013/05/18/the-humble-hero}.
\bibitem[{Tierney et~al.(2014)Tierney, Pacino and
  Jensen}]{Tierney2014OnProblems}
\bibinfo{author}{Tierney, K.}, \bibinfo{author}{Pacino, D.},
  \bibinfo{author}{Jensen, R.M.}, \bibinfo{year}{2014}.
\newblock \bibinfo{title}{{On the complexity of container stowage planning
  problems}}.
\newblock \bibinfo{journal}{Discrete Applied Mathematics}
  \bibinfo{volume}{169}, \bibinfo{pages}{225--230}.
\newblock \URLprefix \url{http://dx.doi.org/10.1016/j.dam.2014.01.005},
  \DOIprefix\doi{10.1016/j.dam.2014.01.005}.
\bibitem[{Wilson et~al.(2001)Wilson, Roach and Ware}]{Wilson2001ContainerStudy}
\bibinfo{author}{Wilson, I.}, \bibinfo{author}{Roach, P.},
  \bibinfo{author}{Ware, J.}, \bibinfo{year}{2001}.
\newblock \bibinfo{title}{{Container stowage pre-planning: using search to
  generate solutions, a case study}}.
\newblock \bibinfo{journal}{Knowledge-Based Systems} \bibinfo{volume}{14},
  \bibinfo{pages}{137--145}.
\newblock \URLprefix
  \url{http://www.sciencedirect.com/science/article/pii/S0950705101000909}.
\bibitem[{Wilson and Roach(2000)}]{Wilson2000ContainerSolutions}
\bibinfo{author}{Wilson, I.D.}, \bibinfo{author}{Roach, P.A.},
  \bibinfo{year}{2000}.
\newblock \bibinfo{title}{{Container stowage planning: a methodology for
  generating computerised solutions}}.
\newblock \bibinfo{journal}{Journal of the Operational Research Society}
  \bibinfo{volume}{51}, \bibinfo{pages}{1248--1255}.
\newblock \URLprefix
  \url{https://www.tandfonline.com/doi/full/10.1057/palgrave.jors.2601022},
  \DOIprefix\doi{10.1057/palgrave.jors.2601022}.
\bibitem[{Wu et~al.(2021)Wu, Xia and Wu}]{Wu2021ResearchShips}
\bibinfo{author}{Wu, Q.}, \bibinfo{author}{Xia, Q.}, \bibinfo{author}{Wu, M.},
  \bibinfo{year}{2021}.
\newblock \bibinfo{title}{{Research on intelligent loading system for container
  ships}}.
\newblock \bibinfo{journal}{IOP Conference Series: Earth and Environmental
  Science} \bibinfo{volume}{632}.
\newblock \DOIprefix\doi{10.1088/1755-1315/632/2/022074}.
\bibitem[{Yifan et~al.(2016)Yifan, Ning and Weijian}]{Yifan2016Group-BayShip}
\bibinfo{author}{Yifan, S.}, \bibinfo{author}{Ning, Z.},
  \bibinfo{author}{Weijian, M.}, \bibinfo{year}{2016}.
\newblock \bibinfo{title}{{Group-Bay Stowage Planning Problem for Container
  Ship}}.
\newblock \bibinfo{journal}{Polish Maritime Research} \bibinfo{volume}{23},
  \bibinfo{pages}{152--159}.
\newblock \DOIprefix\doi{10.1515/pomr-2016-0060}.
\bibitem[{Zhang et~al.(2008)Zhang, Lin, Ji and Zhang}]{Zhang2008}
\bibinfo{author}{Zhang, W.Y.}, \bibinfo{author}{Lin, Y.}, \bibinfo{author}{Ji,
  Z.S.}, \bibinfo{author}{Zhang, G.F.}, \bibinfo{year}{2008}.
\newblock \bibinfo{title}{{Review of containership stowage plans for full
  routes}}.
\newblock \bibinfo{journal}{Journal of Marine Science and Application}
  \bibinfo{volume}{7}, \bibinfo{pages}{278--285}.
\newblock \DOIprefix\doi{10.1007/s11804-008-7087-8}.
\bibitem[{Zhao et~al.(2018)Zhao, Guo, Xiang, Xia, Shen and
  Mi}]{Zhao2018ContainerSearch}
\bibinfo{author}{Zhao, N.}, \bibinfo{author}{Guo, Y.}, \bibinfo{author}{Xiang,
  T.}, \bibinfo{author}{Xia, M.}, \bibinfo{author}{Shen, Y.},
  \bibinfo{author}{Mi, C.}, \bibinfo{year}{2018}.
\newblock \bibinfo{title}{{Container Ship Stowage Based on Monte Carlo Tree
  Search}}.
\newblock \bibinfo{journal}{Journal of Coastal Research} \bibinfo{volume}{83},
  \bibinfo{pages}{540--547}.
\newblock \DOIprefix\doi{10.2112/SI83-090.1}.
\bibitem[{Zhu et~al.(2020)Zhu, Ji and Guo}]{Zhu2020IntegerProblem}
\bibinfo{author}{Zhu, H.}, \bibinfo{author}{Ji, M.}, \bibinfo{author}{Guo, W.},
  \bibinfo{year}{2020}.
\newblock \bibinfo{title}{{Integer Linear Programming Models for the
  Containership Stowage Problem}}.
\newblock \bibinfo{journal}{Mathematical Problems in Engineering}
  \bibinfo{volume}{2020}.
\newblock \DOIprefix\doi{10.1155/2020/4382745}.
\bibitem[{Zurheide and Fischer(2015)}]{Zurheide2015RevenueIndustry}
\bibinfo{author}{Zurheide, S.}, \bibinfo{author}{Fischer, K.},
  \bibinfo{year}{2015}.
\newblock \bibinfo{title}{{Revenue management methods for the liner shipping
  industry}}.
\newblock \bibinfo{journal}{Flexible Services and Manufacturing Journal}
  \bibinfo{volume}{27}, \bibinfo{pages}{200--223}.
\newblock \DOIprefix\doi{10.1007/s10696-014-9192-0}.

\end{thebibliography}
